\documentclass{article}

\usepackage[utf8]{inputenc}
\usepackage[T1]{fontenc}

\usepackage{geometry}
\usepackage{ragged2e}

\usepackage[pagebackref, hidelinks]{hyperref}
\hypersetup{
colorlinks=true, 
citecolor = blue, 
linkcolor = blue, 
urlcolor = blue, 
anchorcolor = blue
}

\usepackage{amsmath, amsthm, amsfonts}
\usepackage{amscd, amssymb, latexsym}

\usepackage[all]{xy}
\usepackage{mathtools} 

\usepackage{comment, color}
\usepackage{enumitem}

\usepackage{svg}
\usepackage{import}

\newtheorem{Theorem}{Theorem}[section]
\newtheorem{Lemma}[Theorem]{Lemma}
\newtheorem{Proposition}[Theorem]{Proposition}
\newtheorem{Corollary}[Theorem]{Corollary}
\newtheorem{Definition}[Theorem]{Definition}
\newtheorem{Example}[Theorem]{Example}
\newtheorem{Remark}[Theorem]{Remark}

\newtheorem{Question}[Theorem]{Question}

\newcommand{\asrt}[1]{\ensuremath{\left[\!\left[ {#1} \right]\!\right]}}

\newcommand{\N}{\mathbb{N}}
\newcommand{\Z}{\mathbb{Z}}

\newcommand{\Field}{\mathbb{K}}

\renewcommand{\S}{\mathbb{S}}

\newcommand{\TT}{\mathcal{T}}

\DeclareMathOperator{\Card}{Card}

\DeclareMathOperator{\rank}{rank}

\DeclareMathOperator{\len}{len}

\DeclareMathOperator{\DL}{\mathfrak{del}}
\DeclareMathOperator{\Lap}{\mathfrak{d}}
\DeclareMathOperator{\E}{\mathfrak{e}}
\DeclareMathOperator{\W}{\mathfrak{w}}
\DeclareMathOperator{\I}{\mathfrak{i}}
\DeclareMathOperator{\CR}{\mathfrak{r}}
\DeclareMathOperator{\CB}{\mathfrak{c}}

\DeclareMathOperator{\ming}{mg}

\title{Loops in surfaces, chord diagrams, interlace graphs:\\ operad factorisations and generating grammars}

\author{Christopher-Lloyd Simon}
\date{\today} 

\begin{document}

\maketitle

\begin{abstract}
A filoop is a generic immersion of a circle in a closed oriented surface, whose complement is a disjoint union of discs, considered up to orientation preserving diffeomorphisms. 
It gives rise to a chord diagram $C$ which has an interlace graph $G$, called a chordiagraph.

For a graph $G$ with even degrees, we compute a quantity $\ming(G)$ which yields, for every chord diagram $C$ with interlace graph $G$, the minimal genus of filoops with chord diagram $C$. If $\ming(G)=0$ then $C$ admits exactly two framings of genus $0$, corresponding to spheriloops.

After recalling the Cunningham factorisation of connected graphs, we describe a canonical factorisation of filoops into spheric sums followed by toric sums, for which the genus is additive. This is analogous to the factorisation of compact connected 3-manifolds along spheres and tori. 

We describe unambiguous context-sensitive grammars generating the set of all graphs and with $\ming=0$ and deduce stability properties with respect to spheric and toric factorisations.
Similar results hold for chordiagraphs with $\ming = 0$ and their corresponding spheriloops.

\end{abstract}

\section*{Acknowledgements}
I wish to thank Etienne Ghys for introducing me to these questions through his book \cite{Ghys_promenade_2017}, and Sergei Tabachnikov for encouraging me to revisit my unpublished results about them: this led to several new results in the present and forthcoming papers.
I am also grateful to Patrick Popescu-Pampu for his careful remarks, and pointing out relevant references towards previous results recovered in this work.

\setcounter{tocdepth}{1}
\tableofcontents

\section{Introduction}

Let $F$ be a smooth surface which is connected oriented and closed, of Euler characteristic $2-2g$.

A \emph{filoop} is a generic immersion $\gamma \colon \S^1 \looparrowright F$ considered up to orientation preserving diffeomorphisms of source and target such that $F\setminus \gamma$ is a disjoint union of discs.


A \emph{chord diagram} is a finite cyclic word $C$ in which each letter appears exactly twice, and a \emph{framing} is a map $\varphi \colon C \to \{\infty,0\}$ which is bijective in restriction to every double occurrence letter.

We recall in Proposition \ref{Prop:framed-chordiag-loop-S} how every framed chord diagram $C_\varphi$ arises from a unique filoop $\gamma$.

To a chord diagram $C$ is associated the \emph{interlace graph} $G$ of its chords.
Such interlace graphs of chord diagrams will be called \emph{chordiagraphs}.

\begin{figure}[h]
    \centering
    \includegraphics[width=\textwidth]{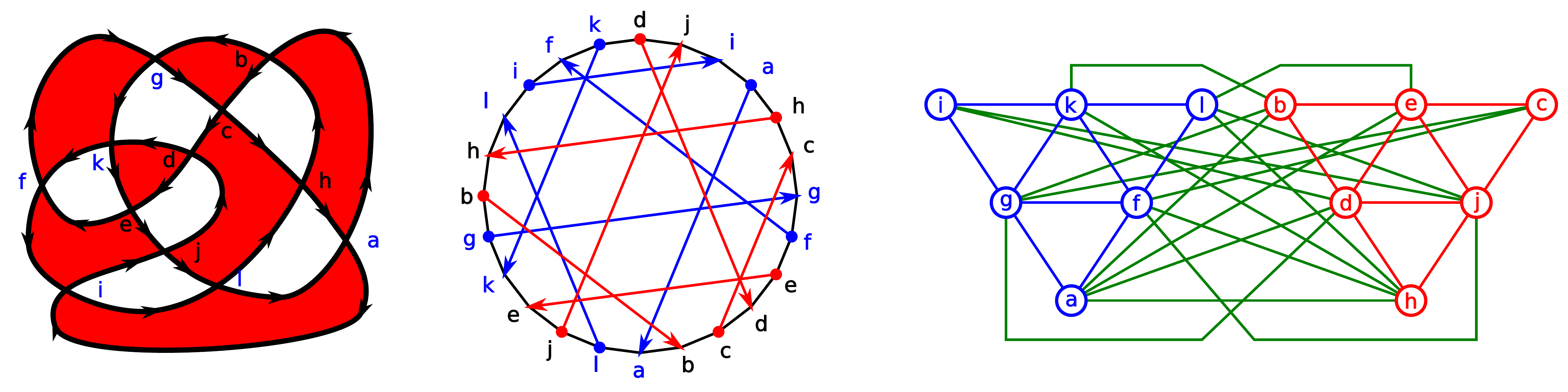}\caption{\label{fig:Intro:AbcDeFbagHdCiGfEhI_chordialoop} A spheriloop, its corresponding framed chord diagram, and its interlace graph.}
\end{figure}

We wish to relate the topology of filoops with the combinatorics of their chord diagrams and interlace graphs, for individual objects and families of objects characterised by certain properties.
Of special interest are the filoops in the sphere, which we call \emph{spheriloops}, and their chordiagraphs.

\subsection{Homological invariants of filoops and their interlace graphs}

The study of spheriloops and their chord diagrams was initiated by Gauss in \cite{Gauss_words-plane-loops_1807}, observing that their interlace graphs must satisfy parity conditions \ref{EveN1} \& \ref{EveN2} but that these were not sufficient. This description was completed by Rosenstiehl in \cite{Rosenstiehl_Gauss-interlace_1999} by adding the homological condition \ref{RoCo} and showing these three conditions characterise chordiagraphs of spheriloops.

\begin{Definition}[Rosenstiehl form \& Gaussian graphs]
\label{Def:Intro:Graph-Conditions}
Consider a simple graph $G$ with adjacency matrix $\E\colon V_G\times V_G \to \{0,1\}$. We define its \emph{Rosenstiehl form} by $\CR=\E+\E^2$, and its restriction to $E_G$ defines the \emph{Rosenstiehl cocycle} $\CR \in Z^1(G;\Z/2)$ (\textcolor{black!50!green}{edges} in Figure~\ref{fig:Intro:AbcDeFbagHdCiGfEhI_chordialoop}).
We say that $G$ satisfies:
\begin{itemize}[noitemsep]
    \item[\ref{EveN1}] when every vertex $x\in V_G$ has an even degree $\CR(x,x) \equiv 0$
    \item[\ref{EveN2}] when distinct non-adjacent vertices $x,y \in V_G$ share an even number of neighbours 
    $\CR(x,y) \equiv 0$
    \item[\ref{RoCo}] when the Rosenstiehl cocycle is null in cohomology, namely $\CR \in B^1(G;\Z/2)$.
\end{itemize}
A graph satisfying \ref{EveN1} and \ref{EveN2} and \ref{RoCo} will be called \emph{Gaussian} (even if it is not a chordiagraph).
\end{Definition}

A \emph{bicolouring} of $\gamma$ is a map $F\setminus \gamma \to \{\bullet, \circ\}$ changing colour when crossing $\gamma$ at smooth points. Such a bicolouring yields a map $\chi^\pm \colon V_G \to \{\bullet, \circ\}$ and a colour form $\CB_\chi \colon V_G\times V_G \to \{0,1\}$.

The filoop is \emph{bicolourable} if and only if $G$ satisfies \ref{EveN1}, in which case it admits two bicolourings.
For a chord diagram $C$ whose interlace graph $G$ satisfies \ref{EveN1}, we show in Proposition \ref{Prop:Bicolor=Framing} that pairs of bicolourings $\chi^{\pm} \colon V_G \to \{\bullet, \circ\}$ correspond to pairs of framings $\varphi^{\pm} \colon C \to \{\infty,0\}$.

Our main Theorems \ref{Thm:intersection-cocycle} and \ref{Thm:intercolorgraphorms}, relating the topology of filoops to the combinatorics of their interlace graphs and chord diagrams will recover, refine and generalise the results of Rosenstiehl.

\begin{Corollary}[Minimal genus]
    \label{Cor:Intro:minimal-genus-bicolor}
    Fix a chordiagraph $G$ satisfying \ref{EveN1} with Rosenstiehl form $\CR$.
    For every chord diagram $C$ with interlace graph $G$, the minimal genus of its framings equals the minimum of the half-rank over the space of bicolour forms $c_\chi\colon V_G\times V_G\to \{0,1\}$ translated by $\CR$:
    \begin{equation*}
        \ming(\E)= \tfrac{1}{2}\min \{\rank(\CR+\CB_\chi) \mid \chi \colon V_G\to \{\bullet, \circ\} \}
    \end{equation*}

    Fix a chordiagraph $G$ satisfying \ref{EveN1} and \ref{EveN2} with Rosenstiehl cocycle $\CR\in Z^1(G;\Z/2)$.    
    For every chord diagram $C$ with interlace graph $G$, the minimal genus of its framings equals the minimum of the half-rank over the space of coboundaries $B^1(G;\Z/2)$ translated by $\CR\in Z^1(G;\Z/2)$. 
\end{Corollary}


\begin{Corollary}[Spheriloops with prescribed chord diagram]
    \label{Cor:Intro:Spherical-framings}
    A chord diagram with a Gaussian interlace graph $G$ admits $\Card H_0(G; \Z/2)=2^{b_0}$ framings of genus $0$. 
    
    Two of them are obtained by integrating the Rosenstiehl form as in Proposition \ref{Prop:Bicolor=Framing} so that consecutive letters have the same framing if and only if they have different colours.
    The others differ by locally-constant change of framings on chords associated to connected components of $G$.
\end{Corollary}

\subsection{Unique factorisations of graphs, chord diagrams and filoops}

We first recall Cunningham's factorisation of graphs and Bouchet's factorisation of chord diagrams.

\begin{Definition}

A \emph{split} of a graph $G$ is a bipartition of its vertices $V=V_0\sqcup V_1$ which is non-trivial (both $V_i$ contain at least $2$ vertices) such that the edges $E(V_0,V_1)$ connecting $V_0$ and $V_1$ are precisely those of the bipartite complete graph $K(U_0,U_1)$ on subsets of vertices $U_0\subset V_0$ and $U_1\subset V_1$.

A connected graph is called \emph{prime} when it has no splits.
A connected graph is called \emph{degenerate} when all its non-trivial bipartitions are splits: those are cliques $K_n$ or stars $S_{n} = K_{1,n-1}$.
\end{Definition}

\begin{Definition}
    A \emph{graph-labelled-tree} (or GLT) is a set of graphs $G_x$ indexed by the nodes of a tree $T$, with bijections $\rho_x \colon V_{x} \to N(x)$ from the vertices of $G_x$ to the neighbours of $x$ in $T$.
    
    The \emph{accessibility graph} $G$ of such a GLT has vertices the leaves of $T$, and distinct vertices $u,v\in V_G$ are connected when for every node $x_i$ of the injective sequence $x_1,\dots,x_l\in V_T$ joining $u=x_0$ to $v=x_{l+1}$, the vertices of $G_{x_i}$ mapped by $\rho_{x_i}$ towards $x_{i-1}$ and $x_{i+1}$ are connected.

    The GLT is \emph{reduced} when: the nodes have degree at least $3$, every $G_x$ is prime or degenerate, and for connected nodes $x,y\in V_G$ we exclude the possibility that $G_x$ and $G_y$ are both cliques, or stars linked by a center and an extremity.
\end{Definition}

\begin{figure}[h]
 \centering
 \def\svgwidth{3.8cm}
 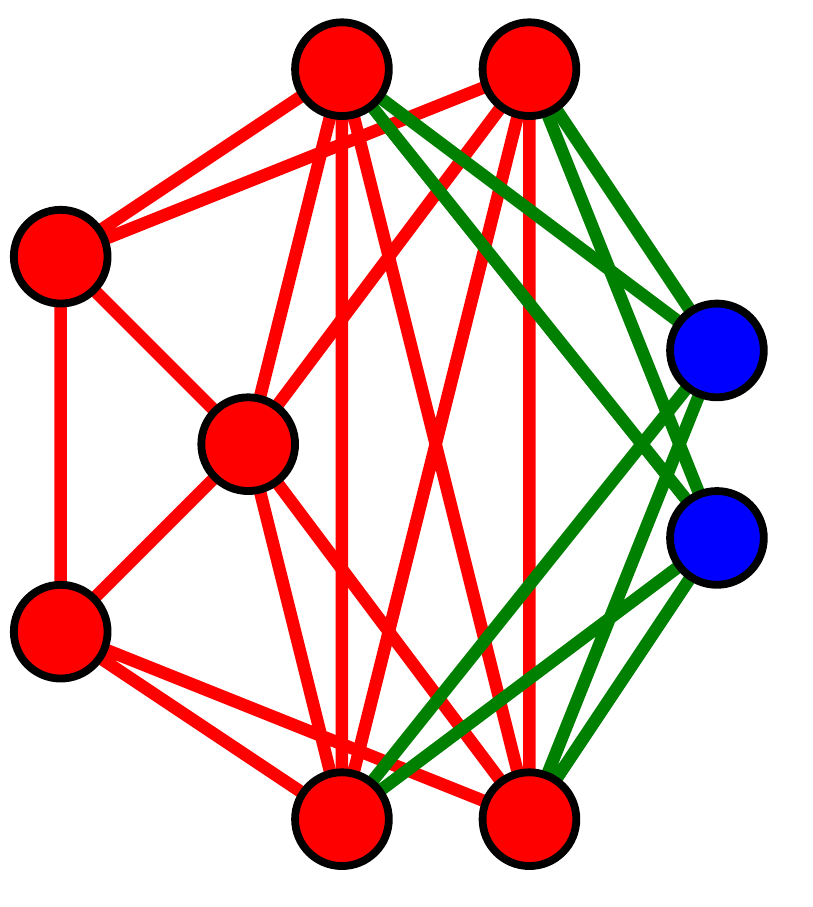
 \hfill
 \def\svgwidth{7cm}
 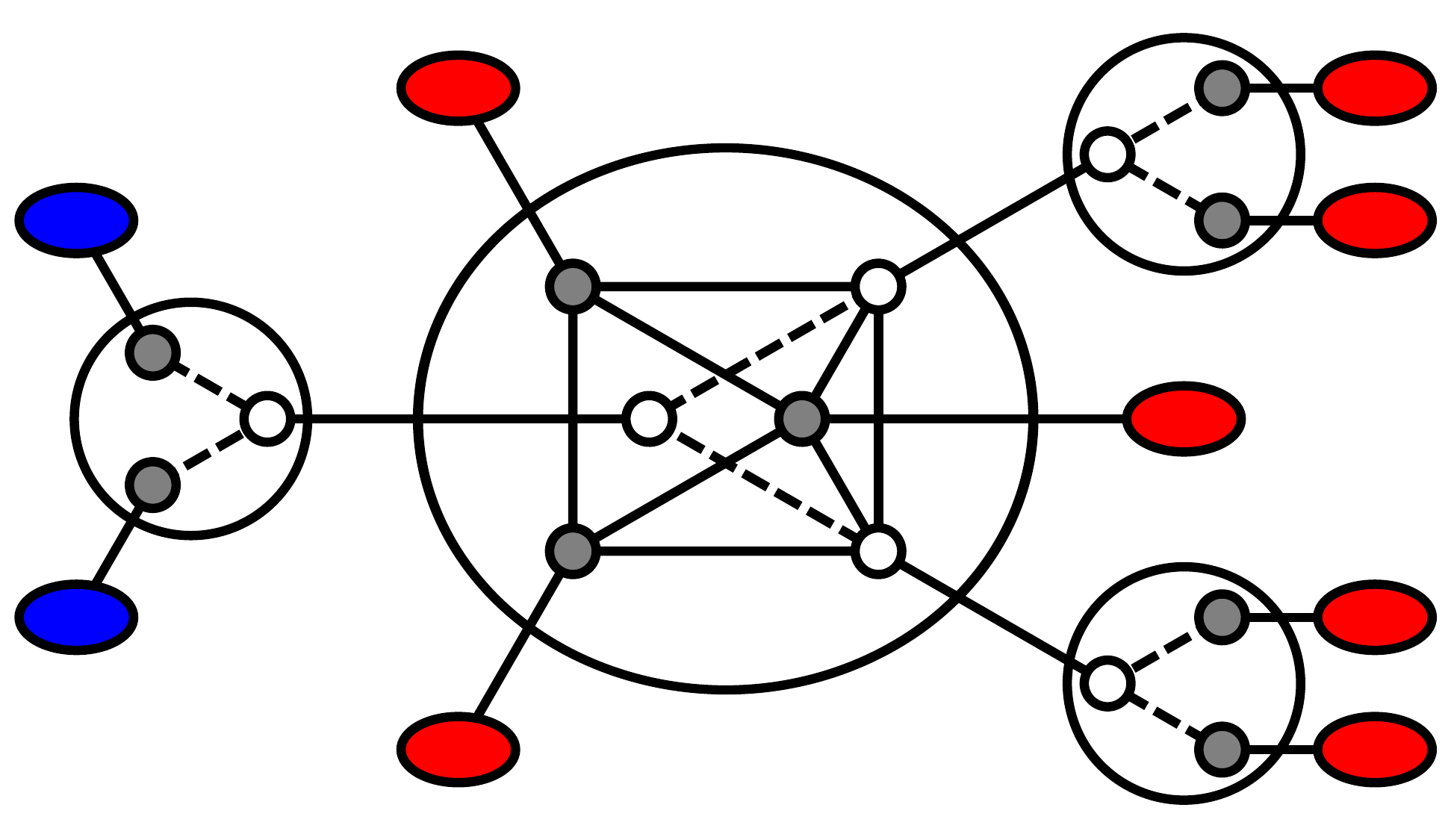
 \caption{\label{fig:Intro:AbcDeFbagHdCiGfEhI_graph} 
 A connected (Gaussian chordia)graph and its reduced GLT factorisation.}
\end{figure}

The \cite[Theorem 3]{Cunningham_decomp-di-graphs_1982} of Cunningham revisited as \cite[Theorem 2.9]{Gio-Paul_split-decomp_2012} by Gioan--Paul states that a connected graph is the accessibility graph of a unique reduced graph-labelled-tree.

Bouchet lifted this unique factorisation to chord diagrams by showing in \cite{Bouchet_prime-graphs-circle-graphs_1987} that a connected chordiagraph, either prime or degenerate, is the interlace graph of a unique chord diagram up to orientation reversal. 
Thus starting from the Cunningham factorisation of a connected graph, one may reconstruct all chord diagrams with that interlace graph by inserting chord diagrams in the place of chords. The only ambiguity arising during this reconstruction appears while replacing a chord by a sub-chord diagram: the choices form a torsor under the rectangular group $(\Z/2)^2$.
In particular, two (unoriented) chord diagrams have the same interlace graphs if and only if they are related by a sequence of \emph{mutations}, which consist in applying a symmetry of the rectangular group $(\Z/2\Z)^2$ to a sub-chord diagram consisting of (at most) two intervals of the circle.

Our main results in section \ref{Sec:FactOperad} consist in lifting these unique factorisations to filoops into \emph{spheric-sums} and \emph{plumbings}, whose precise definitions \ref{Def:spheric-sum-loops} and \ref{Def:plumbing-loops} can be sketched as follows.

\begin{Definition}[Sketch]
A \emph{spheric-sum} of filoops consists in performing a connected sum along discs intersecting each loop in an interval. It corresponds to disjoint unions of interlace graphs.

A \emph{plumbing} of filoops consists in performing a connected sum along discs intersecting each filoop in a pair of transverse intervals.
Essential plumbings correspond to splits of interlace graphs.
\end{Definition}

\begin{figure}[h]
    \centering
    \def\svgwidth{12cm}
    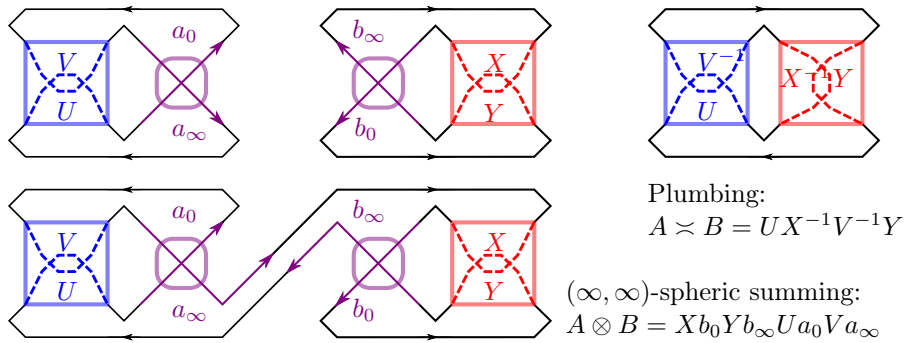
    \caption{Based filoops $\alpha$ and $\beta$, their spheric sum $\alpha \otimes \beta$ and plumbing $\alpha \asymp \beta$.}
    \label{fig:Intro:sum-plumb-loops}
\end{figure}

\begin{Corollary}[Spheric-sum factorisation]
    \label{Cor:Intro:Spheric-sum-facto}
    A filoop admits a unique spheric-sum factorisation into filoops with connected interlace graphs: it has the structure of a plane tree, which can be read off its chord diagram, whose vertices correspond to the connected components of its interlace graph.
    
    The genus is additive under spheric sums.
\end{Corollary}

\begin{Theorem}[Essential-plumbing factorisation]
    \label{Thm:Intro:unifacto-filoops}
    Fix a filoop $\gamma$ with connected interlace graph $G$.
    
    A split factorisation of $G$ into a graph-labelled-tree $(T,G_x)$ yields an essential-plumbing factorisation of $\gamma$ into a filoop-labelled-tree $(T,\gamma_x)$.

    The interlace graph of $\gamma_x$ is the local-complement of $G_x$ at all its vertices in any order, which will prescribe the order in which one must perform the plumbings associated to its incident tree-edges.
    
    The genus of $\gamma$ is the sum of genera of the $\gamma_x$.
\end{Theorem}

\subsection{Generating grammars for Gaussian graphs and spheriloops}

The following algorithm, combining the results in the previous section with Corollary \ref{Cor:Intro:Spherical-framings}, generates all spheriloops with a fixed connected Gaussian chordiagraph in linear time on its number of vertices.

\begin{Corollary} 
\label{Cor:Intro:all-spheriloops-with-graph-G}
Consider a connected Gaussian chordiagraph $G$. The set of all spheriloops whose chord diagrams have interlace graph $G$ are obtained as follows.

\begin{enumerate}[noitemsep]
    \item Compute the Cunningham decomposition $(T,G_x)$ of $G$, and realise every (prime or degenerate) factor by a (unique) unoriented chord diagram.
    
    \item Compose the chord diagrams according to the tree $T$ in all possible ways (at most $4^{k-1}$ choices where $k$ is the number of internal nodes of the tree).
    
    \item Compute the two framings of genus $0$, by integrating the Rosenstiehl form as in Corollary \ref{Cor:Intro:Spherical-framings}.

    \item Construct the loop from the framed chord diagram as in Proposition \ref{Prop:framed-chordiag-loop-S}.
\end{enumerate}

In particular, if a connected Gaussian chordiagraph is prime or degenerate, then it corresponds to a unique loop in the sphere up to homeomorphism.
\end{Corollary}

\begin{figure}[h]
 \centering
 \def\svgwidth{4cm}
 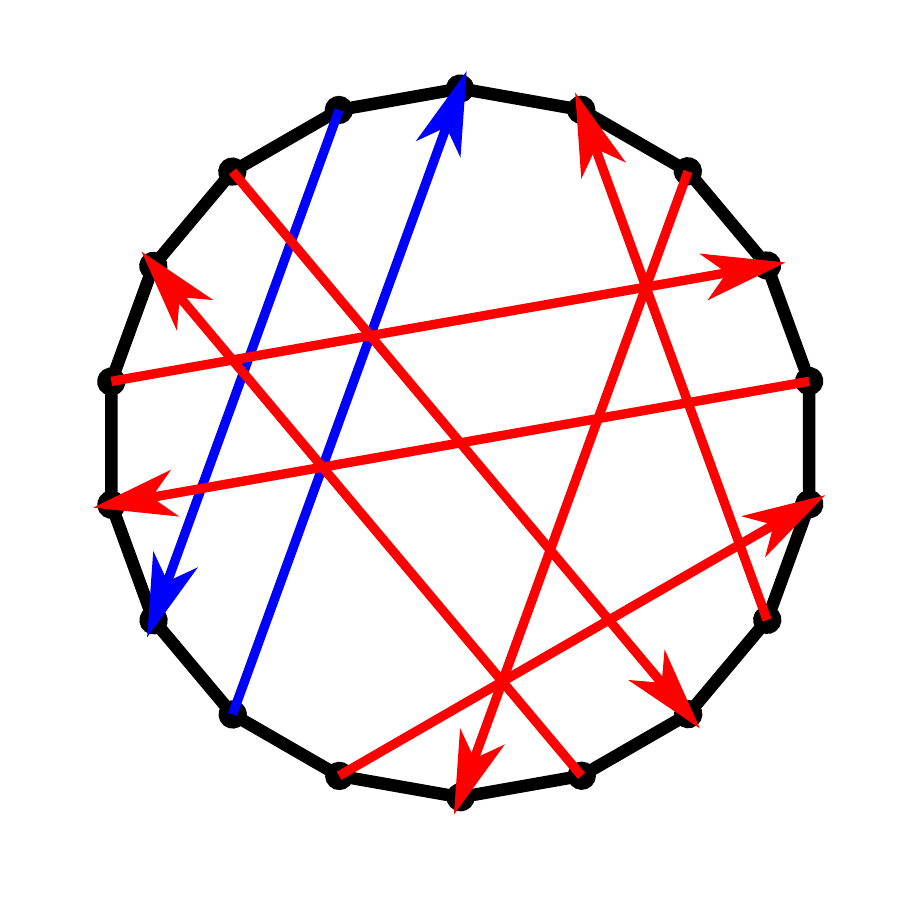
 \hfill
 \def\svgwidth{8cm}
 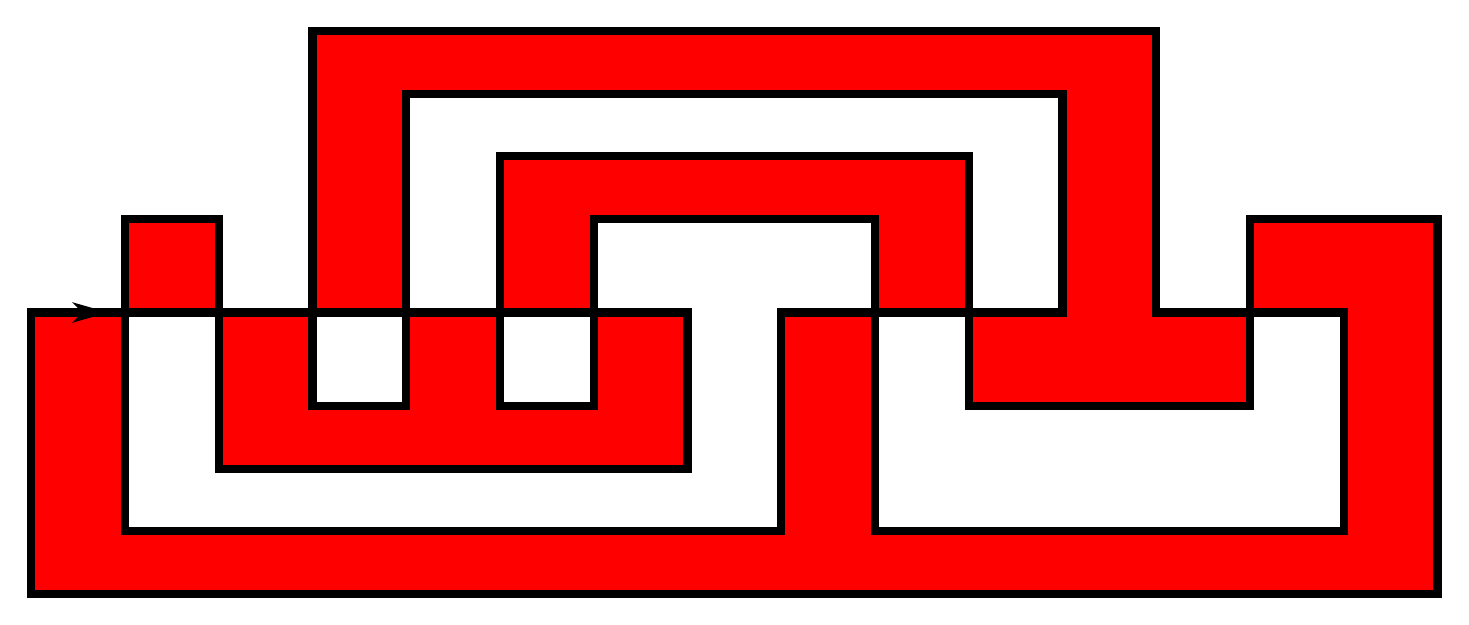\caption{\label{fig:AbcDeFbagHdCiGfEhI_chordialoop} The unique framing of genus $0$ and its corresponding spheriloop.}
\end{figure}

Finally we characterise the GLT factorisations of connected Gaussian graphs in \ref{Prop:GLT-Gauss} after introducing weighted versions of the conditions \ref{EveN1}, \ref{EveN2}, \ref{RoCo}, of which a combination are called \emph{CL2-graphs}.
We therefore deduce that connected Gaussian graphs form a language generated by an arborescent grammar which is context-sensitive and unambiguous.

\begin{Theorem}[Grammatical generation of Gaussian GLT]
\label{Thm:Intro:Construct-Gauss-GLT}
    The following dynamical generation algorithm yields all and only connected Gaussian graphs:
    \begin{enumerate}[noitemsep]
        \item[Atoms] For every (prime or degenerate) graph $G$, compute its CL2-weightings $\bmod{2}$ and construct a weighted GLT whose tree has only one node decorated by $G$.

        \item[Bonds] Choose weighted GLT $A$ and $B$ from those already constructed with leaves $f_a$ and $f_b$ connected to control vertices $u_a\in A_x$ and $u_b\in B_y$ such that the weight of $u_a$ is equal to the degree of $f_b$ and the weight of $f_a$ is equal to the degree of $u_b$, and form the grafting $C = K_2((A,f_a),(B,f_b))$.
        
        \item[Molecules] Among the weighted GLT so constructed, those for which all leaves are connected to control vertices of odd weight and even degree are Gaussian.
    \end{enumerate}
    When restricting to prime and degenerate graphs in the first step, every Gaussian graph is generated only once (possibly up to automorphisms depending on the implementation).
\end{Theorem}

\begin{Remark}[Generating spheriloops]
Adapting the first step in algorithm \ref{Thm:Intro:Construct-Gauss-GLT} by restricting to chordiagraphs (using Bouchet's algorithmic description \cite{Bouchet_Circle-graph-obstructions_1994} recalled in Theorem \ref{Thm:Bouchet-chordiagraph}), yields all Gaussian chordiagraphs.
Combining this with Corollary \ref{Cor:all-spheriloops-with-graph-G} provides an algorithm generating all spheriloops exactly once, which is linear on the number of intersections.
\end{Remark}


\subsection*{Originality and relations with previous works}

A number of ideas and results appearing here have predecessors in the literature.

\paragraph{Section \ref{Sec:Genus-Framings-chordiagraphs}.}
Most of the results in section \ref{Sec:Genus-Framings-chordiagraphs} were known, but its originality lies in the topological methods and the concise self-contained presentation.

The motivating result of Rosenstiehl in Corollary \ref{Cor:Rosenstiehl} was obtained by alternative topological methods in \cite{Chaves-Weber_plombages-rubans-mots-Gauss_1994} and in \cite{Fraysseix-Mendez_Gauss-codes_1999}.

Similar results to Theorem \ref{Thm:intersection-cocycle}, Proposition \ref{Prop:Bicolor=Framing} and Corollary \ref{Cor:Intro:minimal-genus-bicolor} were obtained by combinatorial methods in \cite{Crapo-Rosenstiehl_bicoloops_2001}.
The Corollary \ref{Cor:Intro:minimal-genus-bicolor} was also generalised by Turaev to non-generic loops as well as other related objects, in particular in \cite[Theorem 8.5.4]{Turaev_CurveSurfacesChartsWords_2005}.

However our Corollary \ref{Cor:Intro:Spherical-framings} appears to be new (even though it may be possible to recover it from the more general framework developed by Turaev, and seems related to the \cite[Section 6.2]{Turaev_CurveSurfacesChartsWords_2005} counting certain filoops weighted by their automorphism groups).

\paragraph{Section \ref{Sec:FactOperad}.}
The factorisation Theorem \ref{Thm:Intro:unifacto-filoops} achieved in section \ref{Sec:FactOperad} relies on the combinatorial results of \cite{Gio-Paul_split-decomp_2012, Bouchet_prime-graphs-circle-graphs_1987}, and generalises them to filoops to provide a new insight on their operad structure.

Following Vassiliev's approach to knot theory (whose invariants explained in \cite{ChDuMo_Vassiliev_2012} boil down to framed chord diagrams modulo 1T \& 4T relations), Arnold explored in \cite{Arnold_Topo-plane-curves-caustics_2000} the cohomology of the space of spheriloops (which he considered as a non-commutative version of the former) and believed that a complete understanding is out of reach.
Several investigations pursued these ideas such as \cite{Ito_finite-type-inv-curveSurfaces_2009, Cahn-Levi_VassilInv-VirtuaLegendKnots_2015}.

Our factorisation provides a fathomable description of filoops which leads to the definition of invariants reflecting their topology and complexity (such as the number of prime nodes in the factorisation-tree).
We hope this may lead to a better understanding of their classical invariants.


\paragraph{Section \ref{Sec:GLT-Gaussian}.}
The unambiguous grammatical generation of all Gaussian graphs and spheriloops presented in section \ref{Sec:GLT-Gaussian} are new (they rely on the novel Proposition \ref{Prop:GLT-Gauss} describing their GLT factorisations and Corollary \ref{Cor:Intro:Spherical-framings} yielding the spheric framings of Gaussian chord diagrams).

These results are of interest both from the computational and enumerative perspective, and may also yield random sampling models for Gaussian graphs and spheriloops.

Indeed, they provide an efficient method for tabulating spheriloops according to various refinements on their topological complexity. This problem has gathered much interest for the classification of knot projections (see \cite{Dowker-Thistlethwaite_classification-knot-projections_1983, Weber_classical-knot-theory_2001, Valette_classification-spheriloops-6-gauss-diagram_2016}).

Besides, the unambiguous grammatical generation of much simpler classes of graphs and chord diagrams (yielding context-free grammars) have led to their generation and enumeration, for instance in \cite{CLS_TopoDenoCourbAlgReel_2018} in relation with real algebraic singular plane curves.

The enumeration of spheriloops remains an open question, of interest to physical models of quantum field theory which have been used to conjecture their asymptotic behaviour \cite{Schaeffer-Justin_asymptotic-plane-curves_2004}.

The ideas and results in this paper may contribute to such enumeration, and improve the tabulations of the related combinatorial objects introduced by Turaev initiated in \cite{Gibson_Tabulating-Virtual-strings_2008}.

\subsection*{Further directions of research and relations with other works}

We aim to apply those results in future work to the asymptotic enumeration of spheriloops, the spectral analysis of Gaussian graphs, and the study of their polynomial invariants. In particular...

\paragraph{Enumeration.}

The Remarks \ref{Rem:Vari-LaGraphs} and \ref{Rem:Vari-CL2-Graphs} discuss the algebraic varieties defined by subsets of the conditions \ref{EveN1}, \ref{EveN2}, \ref{RoCo} or the weighted versions thereof.
A further study could make them amenable to the polynomial method explained in \cite{Tao_polynomial-method_2014} to deduce enumerative properties.

Note that a bicolourable linear chord diagram on $[1,2n]$ corresponds to a permutation $\sigma$ pairing letters $2i$ and $2\sigma(i)+1$. Hence there are $(2^n)n!$ bicolourable based filoops with $n$ intersections.
One may try to apply the results in this paper to analyse the generating function $\sum_{n,g} \tfrac{b(g,n)}{n!2^n}x^n y^g$ for the number $b_{n,g}$ of bicolourable filoops with $n$ marked intersections and genus $g$, and deduce asymptotic properties for the representation theory of symmetric groups.
We refer to \cite[Chapter 3]{Lando-Zvonkin_graphs-on-surfaces_2004} and \cite[Section 6.2]{Turaev_CurveSurfacesChartsWords_2005} for some approaches to related enumeration problems.

\paragraph{Polynomial invariants.}

In future work, we will define a ``skein algebra'' of filoops, and apply it to recover the criterion of Lovasz--Marx \cite{Lovasz-Marx_forbidden-substructure-Gauss-codes_1976} telling when a chordiagram admits a framing of genus $0$.
It yields polynomial invariants closely related to the Bollob\'as--Riordan Polynomial of filoops \cite{Chmutov-Pak_Kauffman-virtualink-Tutte_2007} and the Arratia-Bollob\'as polynomial \cite{Arratia-Bollobas-Sorkin_interlace-polynomial_2004} their interlace graphs.
These are also related to the $u$-polynomial and the skein algebra defined and studied in \cite{Turaev_Virtual-strings_2015}.

We believe that the structural results in this paper may help to understand their properties.

\paragraph{Spectral properties.}
Remark \ref{Rem:Laplacian-Lagraphs} expresses \ref{EveN1}, \ref{EveN2}, \ref{RoCo} in terms of the graph Laplacian. How to understand its spectrum in terms of Cunningham's factorisation?
Would Proposition \ref{Prop:GLT-Gauss} help to deduce properties on the spectrum of Gaussian graphs?

After defining \ref{Def:locomp-graph} the local complementation of a graph, it is compelling to consider the action of the group generated by local complementations involutions (indexed by the fixed set of vertices) on the set of graphs (where it satisfies partial commutation and braid relations in  Remark \ref{Rem:locomp-commute-braid}). We wonder what the Schreier graphs associated to the orbits of this action look like. Do they yield interesting families of expander graphs?
Lemma \ref{Lem:GLT-Local-Comp} shows how the reduced GLT decomposition evolves under such operations: this should help in experimenting and solving these problems.

\paragraph{Factorisations of $3$-manifolds}

Theorem \ref{Thm:Intro:unifacto-filoops} is analogous to the unique factorisation of compact connected 3-manifolds into connected sums, first along spheres (by Kneser-Milnor) into prime manifolds, and then along tori (by Jaco--Shalen and Johannson) into geometric manifolds (which are either atoroidal or Seifert-fibered).
We refer to \cite{Thurston_geometrisation_1982, Scott_geometry-3-manifolds_1983} and references therein.
The prime manifolds would correspond to the filoops with connected interlace graphs, among which the atoroidal manifolds would correspond to filoops whose interlace graphs are prime and the Seifert-fibered manifolds would correspond to filoops whose interlace graphs are degenerate.

\paragraph{Knots and singularities.}

The plumbing pictured in Figure~\ref{fig:Intro:sum-plumb-loops} is reminiscent of the tangle addition introduced by Conway in \cite{Conway_enumeration-knots-links-tangles_1970} as an attempt to classify knots according to their diagrams.

We hope that Theorem \ref{Thm:Intro:unifacto-filoops} may help to understand the topology of the lifted knot $\Vec{\gamma}$ in the unit tangent bundle $SF$.
This may have application to the work \cite{FPST_morsifications-mutations_2022} on the topological study of complex analytic singularites from their real morsifications.

\newpage

\section{Homological invariants of filoops and their interlace graphs}

\label{Sec:Genus-Framings-chordiagraphs}

\subsection{Filoops correspond to framed chord diagrams}


Let $F$ be a smooth surface which is connected oriented and closed, of Euler characteristic $2-2g$.

\begin{Definition}
A \emph{loop} is a generic immersion $\gamma \colon \S^1 \looparrowright F$ (all multiple points are transverse double points), considered up to orientation preserving diffeomorphisms of the source and the target.

This amounts to a $4$-valent graph embedded in $F$ up to orientation preserving diffeomorphisms, with a choice of one of the two Eulerian circuits that self-intersect transversely at every vertex.

A \emph{filoop} is a loop which is filling, in the sense that $F\setminus \gamma$ is a disjoint union of discs.
\end{Definition}

Each intersection $x\in\gamma(\S^1)$ has two preimages $x_\infty,x_0\in \S^1$ ordered so that $(\Vec{\gamma}(x_\infty), \Vec{\gamma}(x_0))$ forms a positive basis of the tangent plane $T_{x}F$.
Hence a loop yields a framed chord diagram...


\begin{Definition}

A chord diagram is a finite cyclic word $C$ in which every letter appears twice, considered up to relabelling of the letters.
A framing of $C$ is a function $\varphi\colon C \to \{0,\infty\}$ assigning different values to the two occurrences of each letter.
\end{Definition}

Observe that inverting the orientation of the source $\S^1$ of $\gamma$ corresponds to inverting the cyclic order of source $C$ of $\varphi$, whereas inverting the orientation of the target $F$ of $\gamma$ corresponds to exchanging the values of the target $\{0,\infty\}$ of $\varphi$.

\begin{figure}[h]
    \centering
    \def\svgwidth{11cm}
    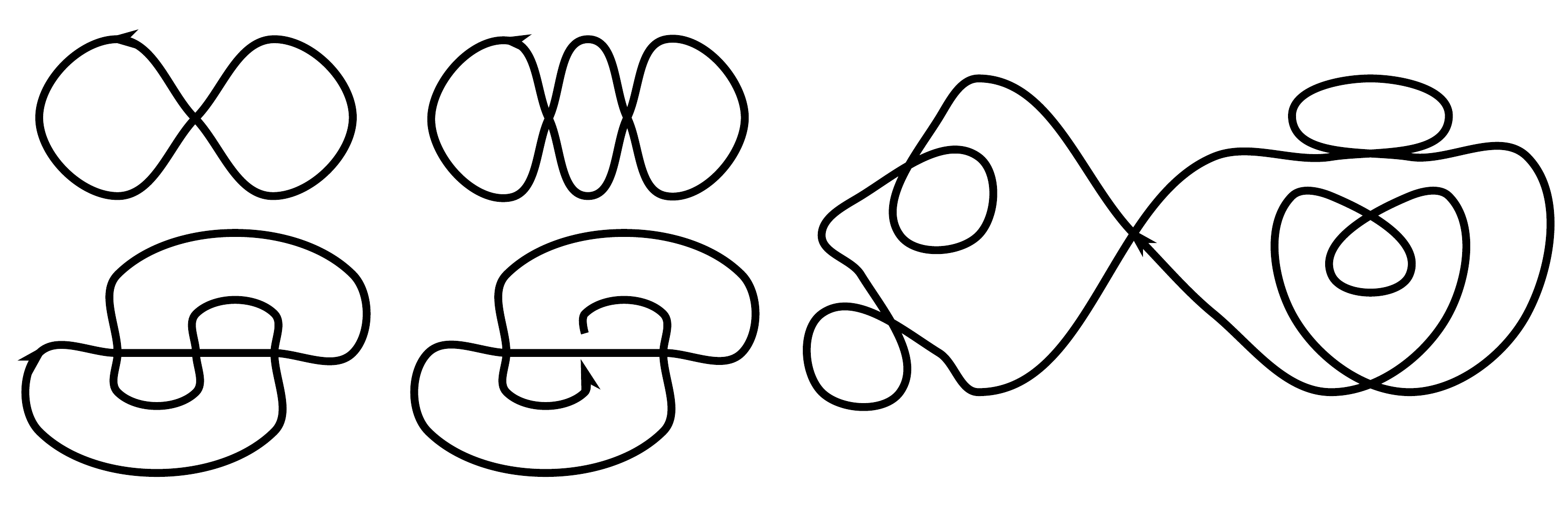
    \def\svgwidth{3cm}
\begingroup%
  \makeatletter%
  \providecommand\color[2][]{%
    \errmessage{(Inkscape) Color is used for the text in Inkscape, but the package 'color.sty' is not loaded}%
    \renewcommand\color[2][]{}%
  }%
  \providecommand\transparent[1]{%
    \errmessage{(Inkscape) Transparency is used (non-zero) for the text in Inkscape, but the package 'transparent.sty' is not loaded}%
    \renewcommand\transparent[1]{}%
  }%
  \providecommand\rotatebox[2]{#2}%
  \newcommand*\fsize{\dimexpr\f@size pt\relax}%
  \newcommand*\lineheight[1]{\fontsize{\fsize}{#1\fsize}\selectfont}%
  \ifx\svgwidth\undefined%
    \setlength{\unitlength}{375bp}%
    \ifx\svgscale\undefined%
      \relax%
    \else%
      \setlength{\unitlength}{\unitlength * \real{\svgscale}}%
    \fi%
  \else%
    \setlength{\unitlength}{\svgwidth}%
  \fi%
  \global\let\svgwidth\undefined%
  \global\let\svgscale\undefined%
  \makeatother%
  \begin{picture}(1,1)%
    \lineheight{1}%
    \setlength\tabcolsep{0pt}%
    \put(0,0){\includegraphics[width=\unitlength,page=1]{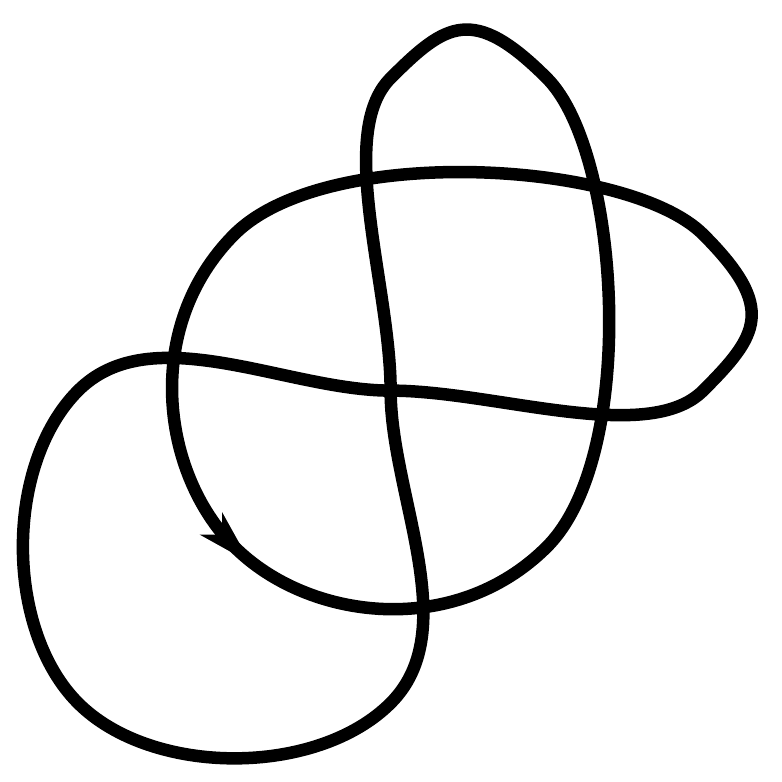}}%
    \put(0.56000001,0.38){\color[rgb]{0,0,0}\makebox(0,0)[lt]{\lineheight{1.25}\smash{\begin{tabular}[t]{l}e\end{tabular}}}}%
    \put(0.80000006,0.49999999){\color[rgb]{0,0,0}\makebox(0,0)[lt]{\lineheight{1.25}\smash{\begin{tabular}[t]{l}b\end{tabular}}}}%
    \put(0.65,0.8){\color[rgb]{0,0,0}\makebox(0,0)[lt]{\lineheight{1.25}\smash{\begin{tabular}[t]{l}c\end{tabular}}}}%
    \put(0.38000001,0.66){\color[rgb]{0,0,0}\makebox(0,0)[lt]{\lineheight{1.25}\smash{\begin{tabular}[t]{l}d\end{tabular}}}}%
    \put(0.25,0.44999998){\color[rgb]{0,0,0}\makebox(0,0)[lt]{\lineheight{1.25}\smash{\begin{tabular}[t]{l}f\end{tabular}}}}%
    \put(0.58000002,0.11999998){\color[rgb]{0,0,0}\makebox(0,0)[lt]{\lineheight{1.25}\smash{\begin{tabular}[t]{l}a\end{tabular}}}}%
  \end{picture}%
\endgroup%

    \caption{The filoops with framed chord diagrams $a_\infty a_0$ and $a_\infty b_0 b_\infty a_0$ and $a_\infty b_0 c_\infty a_0 b_\infty c_0$ and $b_0 a_\infty b_\infty a_0$ and $a_0 b_0 b_\infty c_\infty c_0 a_\infty d_0 d_\infty e_\infty f_\infty f_0 e_0$ and $a_\infty b_0 c_\infty d_\infty e_0 f_\infty b_\infty a_0 f_0 c_0 d_0 e_\infty$.}
    \label{fig:loop_framed-chordiag}
\end{figure}

\begin{Proposition}[filoops = framed chord diagrams]
\label{Prop:framed-chordiag-loop-S}
Every framed chord diagram $\varphi \colon C \to \{0,\infty\}$ arises from a unique filoop $\gamma \colon \S^1 \to F$ up to orientation preserving homeomorphism.
\end{Proposition}

\begin{proof}
The unique framed chord diagram with no chords corresponds to the unique simple loop in the sphere up to orientation preserving homeomorphism. Now assume that $C$ is not empty.

Let us first associate to $C_\varphi$ an oriented $4$-valent connected ribbon-graph $M$, namely a $4$-valent connected graph whose half-edges are cyclically ordered around each vertex.
Each letter $x$ of $C$ corresponds to a ribbon-vertex of $M$, that is an oriented disc with $4$ disjoint closed intervals on its boundary labelled by $x_\infty^+, x_0^+, x_\infty^-,x_0^-$, in that cyclic order, as in Figure~\ref{fig:ribbon-graph}.
We connect these oriented ribbon-vertices by oriented ribbon-edges as prescribed by the framed chord diagram.

\begin{figure}[h]
    \centering
    \def\svgwidth{6cm}
\begingroup%
  \makeatletter%
  \providecommand\color[2][]{%
    \errmessage{(Inkscape) Color is used for the text in Inkscape, but the package 'color.sty' is not loaded}%
    \renewcommand\color[2][]{}%
  }%
  \providecommand\transparent[1]{%
    \errmessage{(Inkscape) Transparency is used (non-zero) for the text in Inkscape, but the package 'transparent.sty' is not loaded}%
    \renewcommand\transparent[1]{}%
  }%
  \providecommand\rotatebox[2]{#2}%
  \newcommand*\fsize{\dimexpr\f@size pt\relax}%
  \newcommand*\lineheight[1]{\fontsize{\fsize}{#1\fsize}\selectfont}%
  \ifx\svgwidth\undefined%
    \setlength{\unitlength}{375bp}%
    \ifx\svgscale\undefined%
      \relax%
    \else%
      \setlength{\unitlength}{\unitlength * \real{\svgscale}}%
    \fi%
  \else%
    \setlength{\unitlength}{\svgwidth}%
  \fi%
  \global\let\svgwidth\undefined%
  \global\let\svgscale\undefined%
  \makeatother%
  \begin{picture}(1,0.8)%
    \lineheight{1}%
    \setlength\tabcolsep{0pt}%
    \put(0,0){\includegraphics[width=\unitlength,page=1]{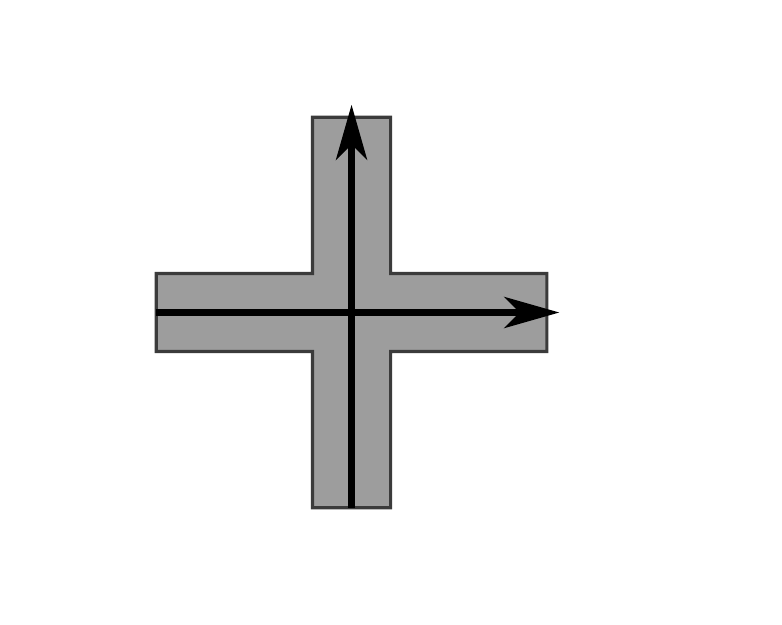}}%
    \put(0.40199986,0.04999995){\color[rgb]{0,0,0}\makebox(0,0)[lt]{\lineheight{1.25}\smash{\begin{tabular}[t]{l}$x_0^-$\end{tabular}}}}%
    \put(0.75000017,0.37799986){\color[rgb]{0,0,0}\makebox(0,0)[lt]{\lineheight{1.25}\smash{\begin{tabular}[t]{l}$x_\infty^+$\end{tabular}}}}%
    \put(0.40599998,0.70199996){\color[rgb]{0,0,0}\makebox(0,0)[lt]{\lineheight{1.25}\smash{\begin{tabular}[t]{l}$x_0^+$\end{tabular}}}}%
    \put(0.05800011,0.37799986){\color[rgb]{0,0,0}\makebox(0,0)[lt]{\lineheight{1.25}\smash{\begin{tabular}[t]{l}$x_\infty^-$\end{tabular}}}}%
  \end{picture}%
\endgroup%

    \def\svgwidth{6cm}
\begingroup%
  \makeatletter%
  \providecommand\color[2][]{%
    \errmessage{(Inkscape) Color is used for the text in Inkscape, but the package 'color.sty' is not loaded}%
    \renewcommand\color[2][]{}%
  }%
  \providecommand\transparent[1]{%
    \errmessage{(Inkscape) Transparency is used (non-zero) for the text in Inkscape, but the package 'transparent.sty' is not loaded}%
    \renewcommand\transparent[1]{}%
  }%
  \providecommand\rotatebox[2]{#2}%
  \newcommand*\fsize{\dimexpr\f@size pt\relax}%
  \newcommand*\lineheight[1]{\fontsize{\fsize}{#1\fsize}\selectfont}%
  \ifx\svgwidth\undefined%
    \setlength{\unitlength}{750bp}%
    \ifx\svgscale\undefined%
      \relax%
    \else%
      \setlength{\unitlength}{\unitlength * \real{\svgscale}}%
    \fi%
  \else%
    \setlength{\unitlength}{\svgwidth}%
  \fi%
  \global\let\svgwidth\undefined%
  \global\let\svgscale\undefined%
  \makeatother%
  \begin{picture}(1,0.7)%
    \lineheight{1}%
    \setlength\tabcolsep{0pt}%
    \put(0,0){\includegraphics[width=\unitlength,page=1]{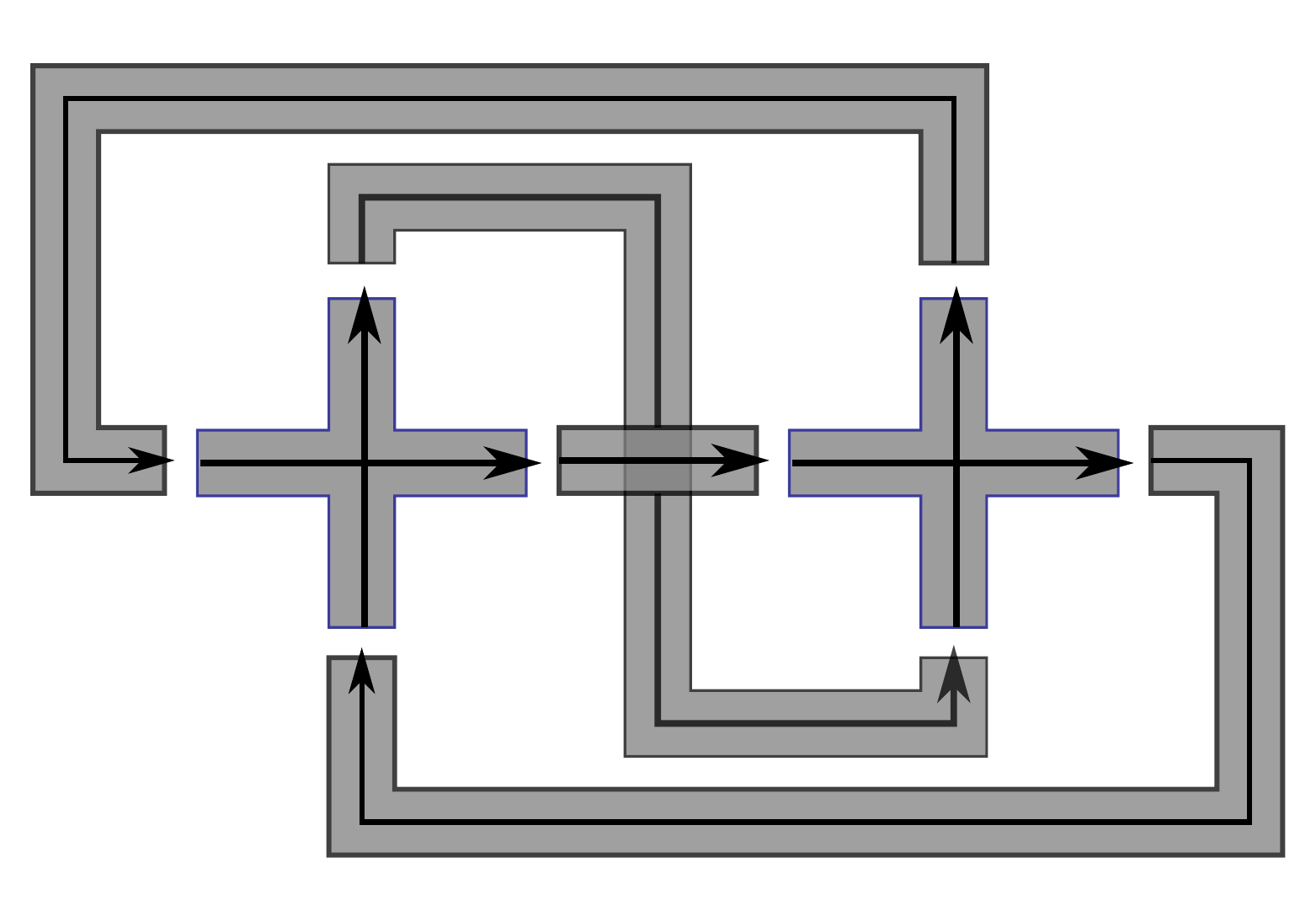}}%
    \put(0.35000002,0.27499999){\color[rgb]{0,0,0}\makebox(0,0)[lt]{\lineheight{1.25}\smash{\begin{tabular}[t]{l}$a_\infty$\end{tabular}}}}%
    \put(0.32500001,0.44999999){\color[rgb]{0,0,0}\makebox(0,0)[lt]{\lineheight{1.25}\smash{\begin{tabular}[t]{l}$a_0$\end{tabular}}}}%
    \put(0.80000001,0.27499999){\color[rgb]{0,0,0}\makebox(0,0)[lt]{\lineheight{1.25}\smash{\begin{tabular}[t]{l}$b_\infty$\end{tabular}}}}%
    \put(0.77500003,0.44999999){\color[rgb]{0,0,0}\makebox(0,0)[lt]{\lineheight{1.25}\smash{\begin{tabular}[t]{l}$b_0$\end{tabular}}}}%
  \end{picture}%
\endgroup%

    \caption{Constructing the ribbon-graph associated to a framed chord diagram}
    \label{fig:ribbon-graph}
\end{figure}
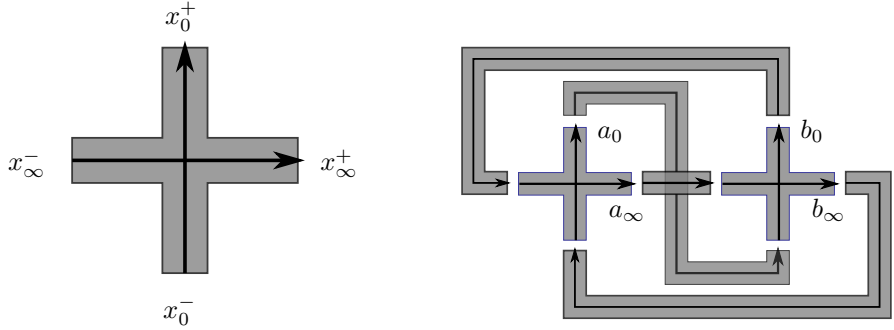

Attach discs to the boundary components of $M$ to obtain a closed connected oriented surface $F$.
The pair $M\subset F$ is uniquely determined up to orientation preserving homeomorphisms by $C_\varphi$.
Moreover the orientation of $C$ selects one of the two Eulerian circuits of $M$ which are transverse at every vertex, hence a filoop $\gamma \colon \S^1\to F$.
\end{proof}

\begin{Remark}[Homology]
    Let $n\in \N^*$ be the number of chords in $C$.
    The $4$-valent map $M$ with $n$ vertices and $2n$ edges has Euler characteristic $-n$ so it is homotopic to a wedge of $1+n$ circles.
    
    The topological maps $M \cap D \hookrightarrow M \sqcup D \twoheadrightarrow M\cup D$ obtained while identifying the boundary of $f$ discs $D = \sqcup D_i$ to the boundary components of $M$, or equivalently $\partial M \hookrightarrow M \sqcup D \twoheadrightarrow F$, yields the Mayer-Vietoris long exact sequence in homology:
    \begin{align*}
        0 &\to H_2(\partial M) \to H_2(M)\oplus H_2(D) \to H_2(F) \to \\
          &\to H_1(\partial M) \to H_1(M)\oplus H_1(D) \to H_1(F) \to \\
          &\to H_0(\partial M) \to H_0(M)\oplus H_0(D) \to H_0(F) \to 0 \\
    \end{align*}
    Since $H_2(M)\oplus H_2(D)=0=H_1(D)$ and $H_0(\partial M) \to H_0(M)\oplus H_0(D)$ is injective, we deduce the short exact sequence $0 \to H_2(F) \to H_1(\partial M) \to H_1(M) \to H_1(F) \to 0$.
    
    Over any field $\Field$ (we will mostly be working over the field $2$ elements), this may be written $0 \to \Field \to \Field^{f} \to \Field^{1+n} \to \Field^{2g} \to 0$. In particular we recover the relation $-1+f-(1+n)+2g=0$.
\end{Remark}

\begin{Remark}[Genus of $C_\varphi$]
	\label{Rem:genus-chordiag}
    If $C$ has $n$ chords and $M$ has $f$ boundary components, then $F$ has Euler characteristic $2-2g=f-n$.
	The boundary components of $M$ can be read off $C_\varphi$ by travelling along the circle, jumping across each chord that is met along the way, and pursuing along the circle in the same or opposite direction according to the framing.

\begin{figure}[h]
    \centering
    \def\svgwidth{6.5cm}
\begingroup%
  \makeatletter%
  \providecommand\color[2][]{%
    \errmessage{(Inkscape) Color is used for the text in Inkscape, but the package 'color.sty' is not loaded}%
    \renewcommand\color[2][]{}%
  }%
  \providecommand\transparent[1]{%
    \errmessage{(Inkscape) Transparency is used (non-zero) for the text in Inkscape, but the package 'transparent.sty' is not loaded}%
    \renewcommand\transparent[1]{}%
  }%
  \providecommand\rotatebox[2]{#2}%
  \newcommand*\fsize{\dimexpr\f@size pt\relax}%
  \newcommand*\lineheight[1]{\fontsize{\fsize}{#1\fsize}\selectfont}%
  \ifx\svgwidth\undefined%
    \setlength{\unitlength}{562.5bp}%
    \ifx\svgscale\undefined%
      \relax%
    \else%
      \setlength{\unitlength}{\unitlength * \real{\svgscale}}%
    \fi%
  \else%
    \setlength{\unitlength}{\svgwidth}%
  \fi%
  \global\let\svgwidth\undefined%
  \global\let\svgscale\undefined%
  \makeatother%
  \begin{picture}(1,0.53333333)%
    \lineheight{1}%
    \setlength\tabcolsep{0pt}%
    \put(0,0){\includegraphics[width=\unitlength,page=1]{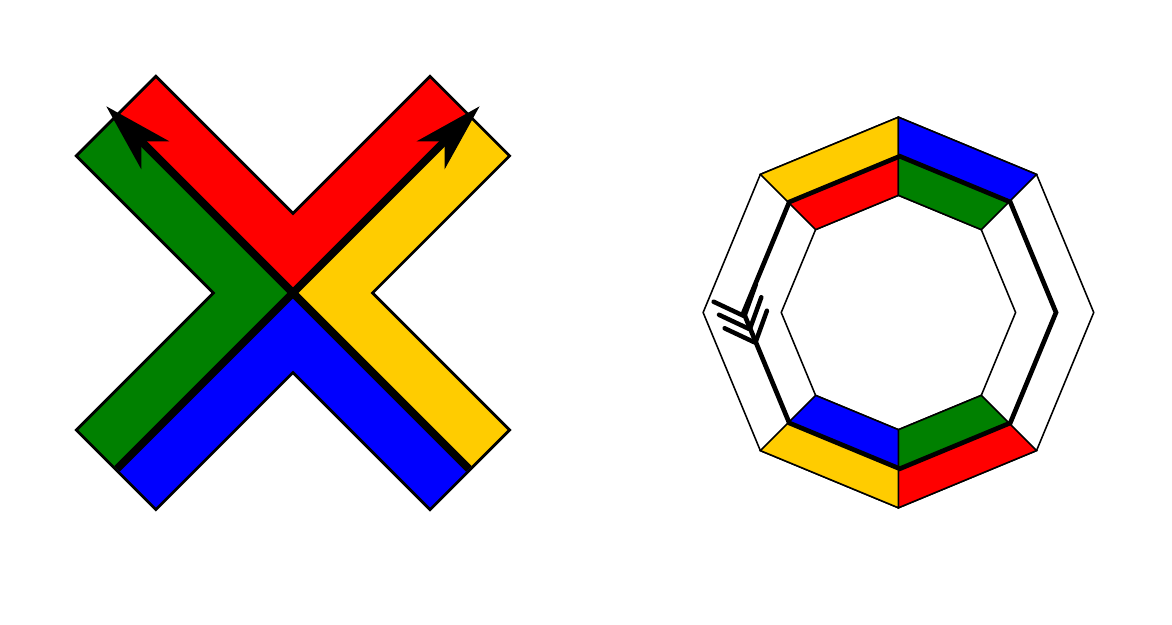}}%
    \put(0.19999999,0.43333334){\color[rgb]{0,0,0}\makebox(0,0)[lt]{\lineheight{1.25}\smash{\begin{tabular}[t]{l}$F_{++}$\end{tabular}}}}%
    \put(0.36666668,0.26666667){\color[rgb]{0,0,0}\makebox(0,0)[lt]{\lineheight{1.25}\smash{\begin{tabular}[t]{l}$F_{+-}$\end{tabular}}}}%
    \put(0.03333333,0.26666667){\color[rgb]{0,0,0}\makebox(0,0)[lt]{\lineheight{1.25}\smash{\begin{tabular}[t]{l}$F_{-+}$\end{tabular}}}}%
    \put(0.19999999,0.1){\color[rgb]{0,0,0}\makebox(0,0)[lt]{\lineheight{1.25}\smash{\begin{tabular}[t]{l}$F_{--}$\end{tabular}}}}%
    \put(0.69999998,0.46666667){\color[rgb]{0,0,0}\makebox(0,0)[lt]{\lineheight{1.25}\smash{\begin{tabular}[t]{l}$x_\infty$\end{tabular}}}}%
    \put(0.69999998,0.03333335){\color[rgb]{0,0,0}\makebox(0,0)[lt]{\lineheight{1.25}\smash{\begin{tabular}[t]{l}$x_0$\end{tabular}}}}%
    \put(0,0){\includegraphics[width=\unitlength,page=2]{loops_intersection_regions.pdf}}%
  \end{picture}%
\endgroup%

    \hfill
    \def\svgwidth{3.5cm}
    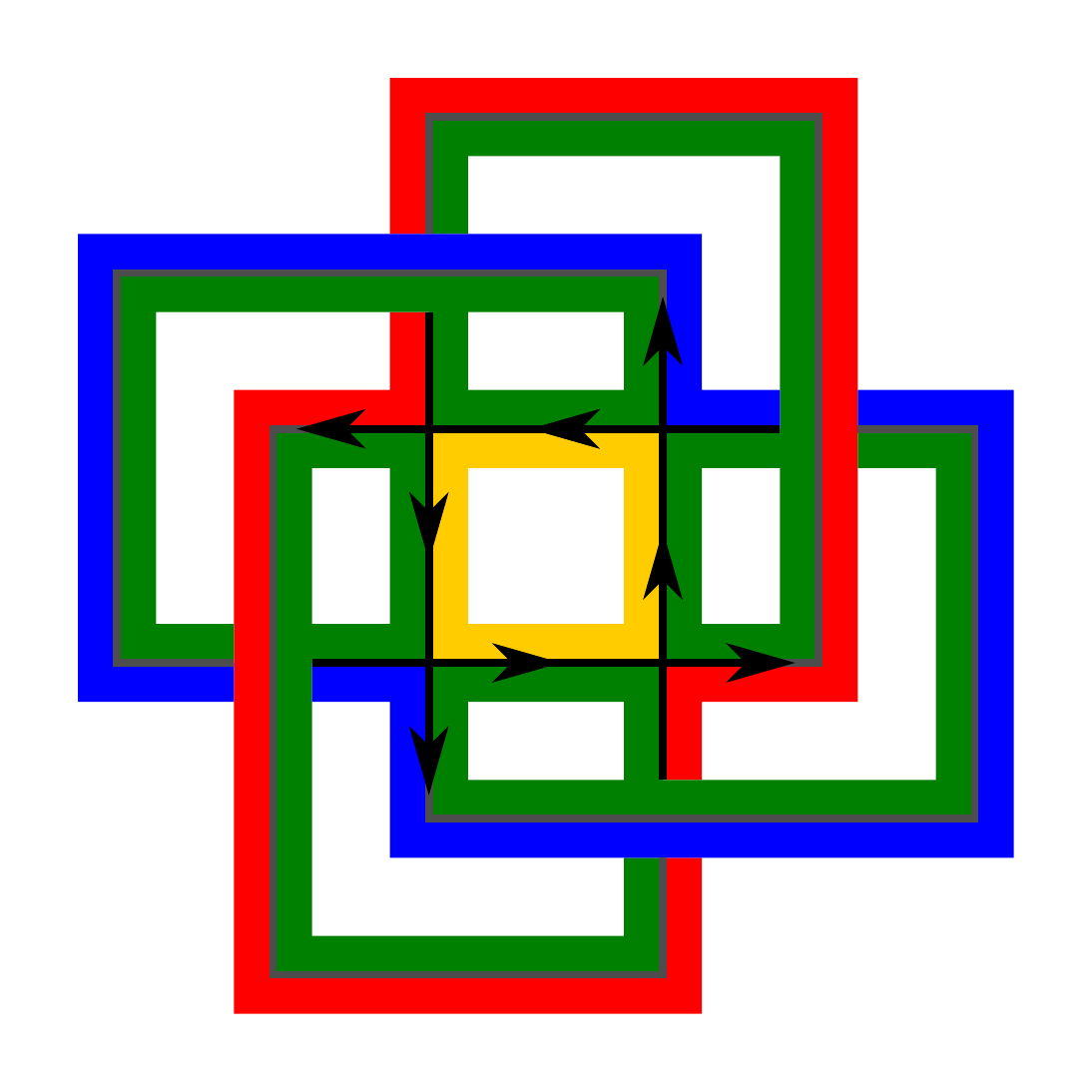
    \def\svgwidth{3.5cm}
    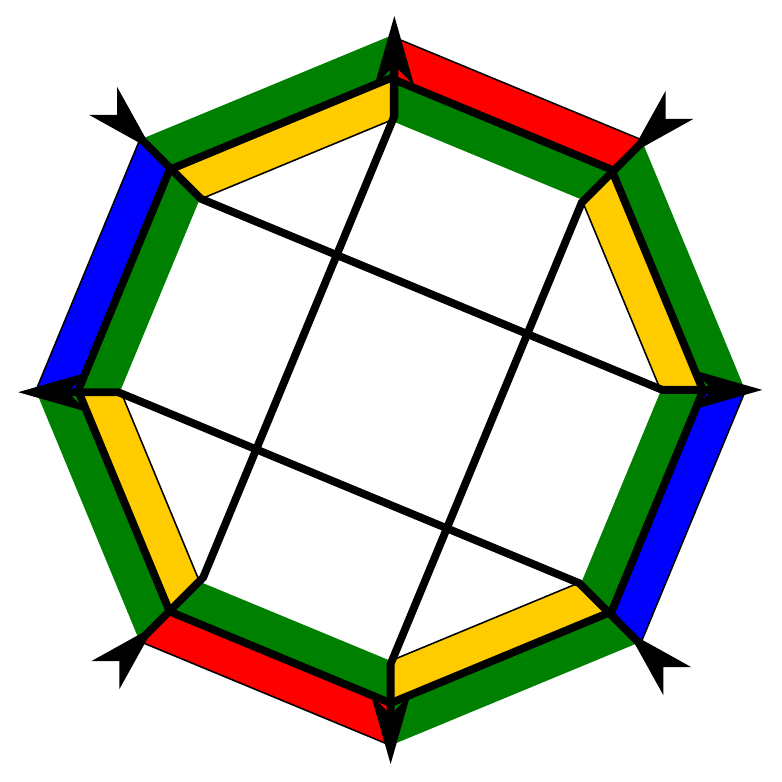
    \caption{Reading the boundary components of $F\setminus \gamma$, for example on $C_\varphi=a_\infty c_0 d_\infty b_0 c_\infty a_0 d_\infty b_0$.}
    \label{fig:faces-from-chordiag}
\end{figure}
\end{Remark}

\subsection{Intersection form of the ribbon-graph}

    
Let us define a couple of functions which will help us understand the relationship between the topology of a filoop $\gamma \colon \S^1 \to F$ and the combinatorics of its associated framed chord diagram $C_\varphi$.
The goal is to compute the symplectic intersection form on $H_1(M;\Z/2)$, whose rank is equal to that of $H_1(F;\Z/2)$ namely twice the genus of $F$.

\begin{Definition}[Graphs and Rosenstiehl form]
\label{Def:Graph-Conditions}
A simple graph $G$ on a set of vertices $V_G$ is encoded by its adjacency matrix $\E\colon V_G\times V_G \to \{0,1\}$, which is symmetric and vanishes on the diagonal.

We define its \emph{Rosenstiehl} symmetric bilinear form $\CR \colon V_G\times V_G \to \Z/2$ by $\CR=\E+\E^2$, so that: 
\begin{align*}
    \tag{RC0}
    \forall (x,y)\notin E_G, \quad
    \CR(x,y) = 0+\Card N(x)\cap N(y) \\
    \tag{RC1}
    \forall (x,y)\in E_G, \quad
    \CR(x,y) = 1+\Card N(x)\cap N(y)
\end{align*}
A symmetric function $V_G\times V_G\to \Z/2$ restricted to $E_G$ yields a $1$-cochain, which is a $1$-cocycle as $G$ has no $2$-cells: we call it the associated cocycle in $Z^1(G;\Z/2)$ and denote it by the same letter.
\end{Definition}

\begin{Definition}[Even and Gaussian conditions]
We say that a graph $G$ satisfies the property:
\begin{itemize}[align=left, noitemsep]
    \item[\ref{EveN1}] when every vertex has an even degree: 
    \begin{equation} \label{EveN1} \tag{EN1}
        \forall x\in V_G, 
        \qquad 
        \CR(x,x) \equiv \Card N(x) \equiv 0
    \end{equation}
    \item[\ref{EveN2}] when distinct non-adjacent vertices share an even number of neighbours: 
    \begin{equation} \label{EveN2}\tag{EN2}
        \forall x,y \in V_G,
        \quad x \ne y,
        \qquad
        \E(x,y) = 0 \implies 
        \CR(x,y) \equiv \Card N(x) \cap N(y) \equiv 0
    \end{equation}
    \item[\ref{RoCo}] when the Rosenstiehl cocycle is null in cohomology, namely $\CR \in B^1(G;\Z/2)$, or when: 
    \begin{equation} \label{RoCo}\tag{RC}
        \forall (x_1,\dots,x_l=x_0) \in V_G^{l}, \qquad \textstyle\prod_1^l \E(x_{i-1},x_{i})\ne 0 \implies \textstyle\sum_1^l \CR(x_{i-1},x_{i}) \equiv 0
    \end{equation}
\end{itemize}

A graph will be called \emph{Gaussian} when it satisfies \ref{EveN1} and \ref{EveN2} and \ref{RoCo}.
\end{Definition}

\begin{Definition}[Interlace graph of a chordiag]
	To a chord diagram $C$ we associate its \emph{interlace graph} $G$, with vertices $V_G$ the set of chords and edges $E_G$ the set of pairs of intersecting chords.

    A \emph{chordiagraph} is the interlace graph of some chord diagram (those will be described in \ref{SubSec:chordiagraphs}).
\end{Definition}

\begin{Definition}[Intersection form of a framed chordiag]
    To a framed chord diagram $C_\varphi$ with interlace graph $G$, we associate its \emph{intersection form} $\I_\varphi \colon V_G \times V_G \to \Z/2$ defined by:
	\begin{equation} \label{GhyCo} \tag{GC}
		\I_\varphi(x,y)= \E(x,y)+\Card \{z_j \in C_\varphi \mid z_j\in (x_\infty,x_0) \wedge z_{1/j} \in (y_\infty,y_0)\}
	\end{equation}
	
	Note that a chord $z\in C$ contributes $1\bmod{2}$ to $\I_{\varphi}(x,y)-\E(x,y)$ when its endpoints both lie in $(x_\infty,x_0)\cup (y_\infty,y_0)$ but do not both lie in $(x_\infty,x_0)\cap (y_\infty,y_0)$, or equivalently when it intersects at least one of the chords $x$ or $y$ and its endpoints both lie in the union of intervals $(x_\infty,x_0)\cup (y_\infty,y_0)$.

    The restriction of $\I_\varphi$ to $E_G$ yields the \emph{intersection cocycle} still denoted $\I_\varphi \in Z^1(G;\Z/2)$.
\end{Definition}

\begin{figure}[h]
    \centering
    \def\svgwidth{3.5cm}
\begingroup%
  \makeatletter%
  \providecommand\color[2][]{%
    \errmessage{(Inkscape) Color is used for the text in Inkscape, but the package 'color.sty' is not loaded}%
    \renewcommand\color[2][]{}%
  }%
  \providecommand\transparent[1]{%
    \errmessage{(Inkscape) Transparency is used (non-zero) for the text in Inkscape, but the package 'transparent.sty' is not loaded}%
    \renewcommand\transparent[1]{}%
  }%
  \providecommand\rotatebox[2]{#2}%
  \newcommand*\fsize{\dimexpr\f@size pt\relax}%
  \newcommand*\lineheight[1]{\fontsize{\fsize}{#1\fsize}\selectfont}%
  \ifx\svgwidth\undefined%
    \setlength{\unitlength}{375bp}%
    \ifx\svgscale\undefined%
      \relax%
    \else%
      \setlength{\unitlength}{\unitlength * \real{\svgscale}}%
    \fi%
  \else%
    \setlength{\unitlength}{\svgwidth}%
  \fi%
  \global\let\svgwidth\undefined%
  \global\let\svgscale\undefined%
  \makeatother%
  \begin{picture}(1,1)%
    \lineheight{1}%
    \setlength\tabcolsep{0pt}%
    \put(0,0){\includegraphics[width=\unitlength,page=1]{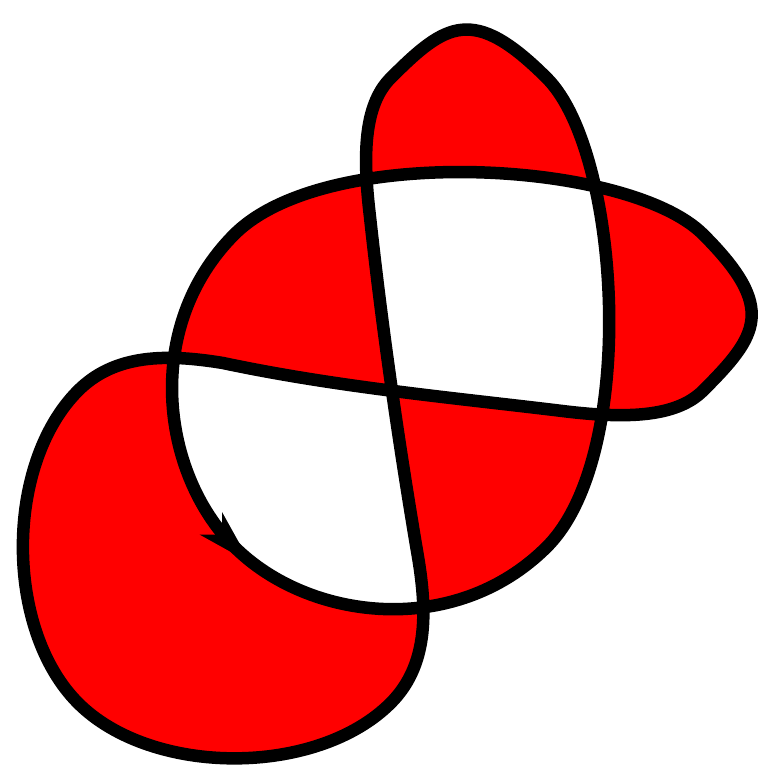}}%
    \put(0.53999999,0.40400002){\color[rgb]{0,0,0}\makebox(0,0)[lt]{\lineheight{1.25}\smash{\begin{tabular}[t]{l}e\end{tabular}}}}%
    \put(0.80000006,0.49999999){\color[rgb]{0,0,0}\makebox(0,0)[lt]{\lineheight{1.25}\smash{\begin{tabular}[t]{l}b\end{tabular}}}}%
    \put(0.65400001,0.80800001){\color[rgb]{0,0,0}\makebox(0,0)[lt]{\lineheight{1.25}\smash{\begin{tabular}[t]{l}c\end{tabular}}}}%
    \put(0.4,0.65399999){\color[rgb]{0,0,0}\makebox(0,0)[lt]{\lineheight{1.25}\smash{\begin{tabular}[t]{l}d\end{tabular}}}}%
    \put(0.25,0.44999998){\color[rgb]{0,0,1}\makebox(0,0)[lt]{\lineheight{1.25}\smash{\begin{tabular}[t]{l}f\end{tabular}}}}%
    \put(0.57200002,0.14399995){\color[rgb]{0,0,1}\makebox(0,0)[lt]{\lineheight{1.25}\smash{\begin{tabular}[t]{l}a\end{tabular}}}}%
  \end{picture}%
\endgroup%

    \hfill
    \def\svgwidth{4.2cm}
    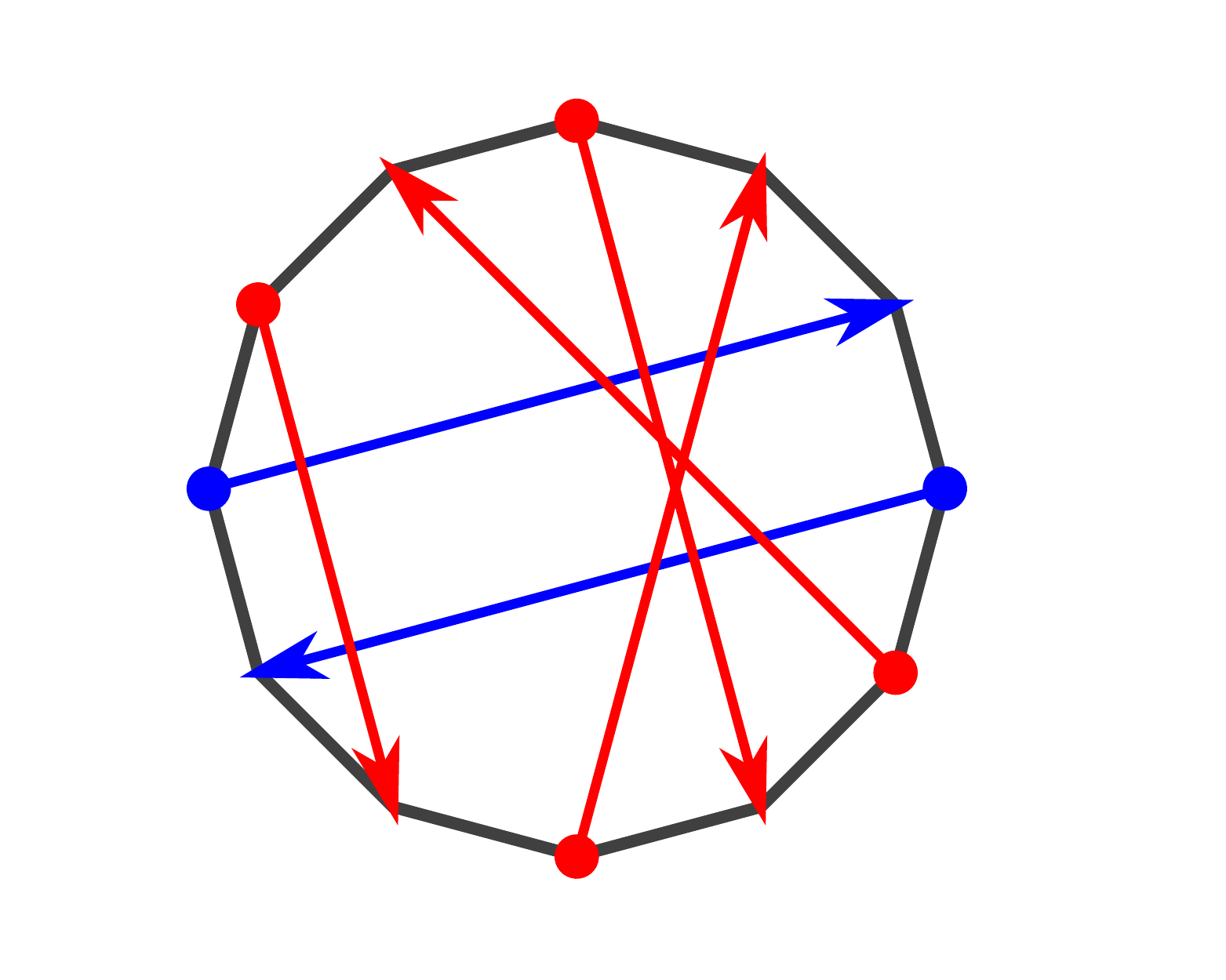
    \hfill
    \def\svgwidth{5.6cm}
\begingroup%
  \makeatletter%
  \providecommand\color[2][]{%
    \errmessage{(Inkscape) Color is used for the text in Inkscape, but the package 'color.sty' is not loaded}%
    \renewcommand\color[2][]{}%
  }%
  \providecommand\transparent[1]{%
    \errmessage{(Inkscape) Transparency is used (non-zero) for the text in Inkscape, but the package 'transparent.sty' is not loaded}%
    \renewcommand\transparent[1]{}%
  }%
  \providecommand\rotatebox[2]{#2}%
  \newcommand*\fsize{\dimexpr\f@size pt\relax}%
  \newcommand*\lineheight[1]{\fontsize{\fsize}{#1\fsize}\selectfont}%
  \ifx\svgwidth\undefined%
    \setlength{\unitlength}{525bp}%
    \ifx\svgscale\undefined%
      \relax%
    \else%
      \setlength{\unitlength}{\unitlength * \real{\svgscale}}%
    \fi%
  \else%
    \setlength{\unitlength}{\svgwidth}%
  \fi%
  \global\let\svgwidth\undefined%
  \global\let\svgscale\undefined%
  \makeatother%
  \begin{picture}(1,0.57142857)%
    \lineheight{1}%
    \setlength\tabcolsep{0pt}%
    \put(0,0){\includegraphics[width=\unitlength,page=1]{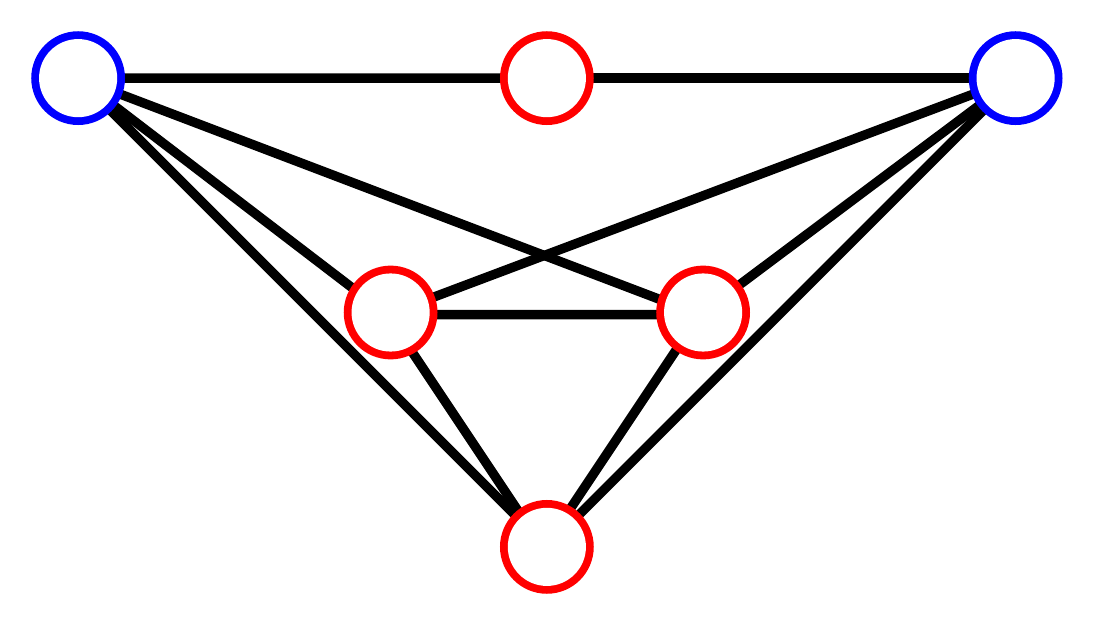}}%
    \put(0.05479917,0.48440295){\color[rgb]{0,0,0}\makebox(0,0)[lt]{\lineheight{1.25}\smash{\begin{tabular}[t]{l}a\end{tabular}}}}%
    \put(0.48162951,0.05583151){\color[rgb]{0,0,0}\makebox(0,0)[lt]{\lineheight{1.25}\smash{\begin{tabular}[t]{l}c\end{tabular}}}}%
    \put(0.33797439,0.26441133){\color[rgb]{0,0,0}\makebox(0,0)[lt]{\lineheight{1.25}\smash{\begin{tabular}[t]{l}b\end{tabular}}}}%
    \put(0.62573946,0.26441133){\color[rgb]{0,0,0}\makebox(0,0)[lt]{\lineheight{1.25}\smash{\begin{tabular}[t]{l}d\end{tabular}}}}%
    \put(0.9173131,0.47829248){\color[rgb]{0,0,0}\makebox(0,0)[lt]{\lineheight{1.25}\smash{\begin{tabular}[t]{l}f\end{tabular}}}}%
    \put(0.48236613,0.48440295){\color[rgb]{0,0,0}\makebox(0,0)[lt]{\lineheight{1.25}\smash{\begin{tabular}[t]{l}e\end{tabular}}}}%
  \end{picture}%
\endgroup%

    \caption{A filoop, its framed chord diagram, and its interlace graph with the intersection cocycle.}
    \label{fig:Interlace-graph-Intersection-cocyle}
\end{figure}

\begin{Definition}
\label{Def:gamma=alpha_x.beta_x}
In the surface $M\subset F$, we define for every intersection $x\in \gamma$ the loops $\alpha_x$ and $\beta_x$ running along $\gamma$ from $x_\infty$ to $x_0$ and $x_0$ to $x_\infty$ respectively.

Since $\gamma = \alpha_x\beta_x$ in $\pi_1(M,x)$, we have $[\gamma]=[\alpha_x]+[\beta_x]$ in $H_1(M;\Z/2)$.
\end{Definition}

\begin{Proposition}[Basis and intersection of $H_1(M;\Z/2)$]
\label{Proposition:BasInter-H1M}
The collections of $n+1$ cycles $[\alpha_x],[\gamma]$ and $[\beta_x],[\gamma]$ form two bases of $H_1(M;\Z/2)$.
We will use the former as our preferred basis.

In that basis, the symplectic form on $H_1(M;\Z/2)$ can be expressed in terms of $C_\varphi$ as:
\begin{equation*}
    [\gamma]\cdot[\alpha_x]= \Card N(x)
    \qquad \qquad
	[\alpha_x]\cdot[\alpha_y]= \I_\varphi(x,y)
\end{equation*}

\end{Proposition}

\begin{proof}
	The proof can be found in \cite[Page 232 to 236]{Ghys_promenade_2017}. 
	Let us briefly explain how to compute the intersection form in this basis using $C_\varphi$.
	Observe that the algebraic and geometric intersection numbers coincide modulo two.
	
	The chords of $C$ that interlace $x\in C$ correspond to the other intersections of $\gamma$ which appear once along $\alpha_x$, and their parity is $[\gamma] \cdot [\alpha_x]$.
	For $x,y\in C$, consider slight parallel translations of $\alpha_x$ to the left and of $\alpha_y$ to the right of $\gamma$ and observe (depending on whether $x$ and $y$ are interlaced) that their intersection number has the same parity as the quantity defining $\I_\varphi(x,y)$.
\end{proof}

\begin{Lemma} \label{Lemma:framing-cobound}
    If two framings $\varphi,\varphi'$ of a chord diagram $C$ differ only at $z\in V_G$, then their associated intersection forms are related by $\I_{\varphi'}(x,y)=\I_\varphi(x,y)$ when $x=z=y$ or $x\ne z \ne y$, and:
\begin{equation*}
	\forall x\in V_G\setminus\{z\}, \quad
	\I_{\varphi'}(x,z)-\I_\varphi(x,z)
    = \Card N(x)\setminus \{z\} 
    = \Card N(x)-\E(x,z)
    = [\alpha_x]\cdot [\beta_z]
\end{equation*}
\end{Lemma}
\begin{proof}
The first two equalities follow from Proposition \ref{Proposition:BasInter-H1M} and the last is represented in Figure~\ref{fig:Ghys_framing_coboundary}.
In particular, the last intersection takes the same value with respect to $\varphi$ or $\varphi'$.
\begin{figure}[h]
    \centering
    \def\svgwidth{3.5cm}
\begingroup%
  \makeatletter%
  \providecommand\color[2][]{%
    \errmessage{(Inkscape) Color is used for the text in Inkscape, but the package 'color.sty' is not loaded}%
    \renewcommand\color[2][]{}%
  }%
  \providecommand\transparent[1]{%
    \errmessage{(Inkscape) Transparency is used (non-zero) for the text in Inkscape, but the package 'transparent.sty' is not loaded}%
    \renewcommand\transparent[1]{}%
  }%
  \providecommand\rotatebox[2]{#2}%
  \newcommand*\fsize{\dimexpr\f@size pt\relax}%
  \newcommand*\lineheight[1]{\fontsize{\fsize}{#1\fsize}\selectfont}%
  \ifx\svgwidth\undefined%
    \setlength{\unitlength}{487.5bp}%
    \ifx\svgscale\undefined%
      \relax%
    \else%
      \setlength{\unitlength}{\unitlength * \real{\svgscale}}%
    \fi%
  \else%
    \setlength{\unitlength}{\svgwidth}%
  \fi%
  \global\let\svgwidth\undefined%
  \global\let\svgscale\undefined%
  \makeatother%
  \begin{picture}(1,1)%
    \lineheight{1}%
    \setlength\tabcolsep{0pt}%
    \put(0,0){\includegraphics[width=\unitlength,page=1]{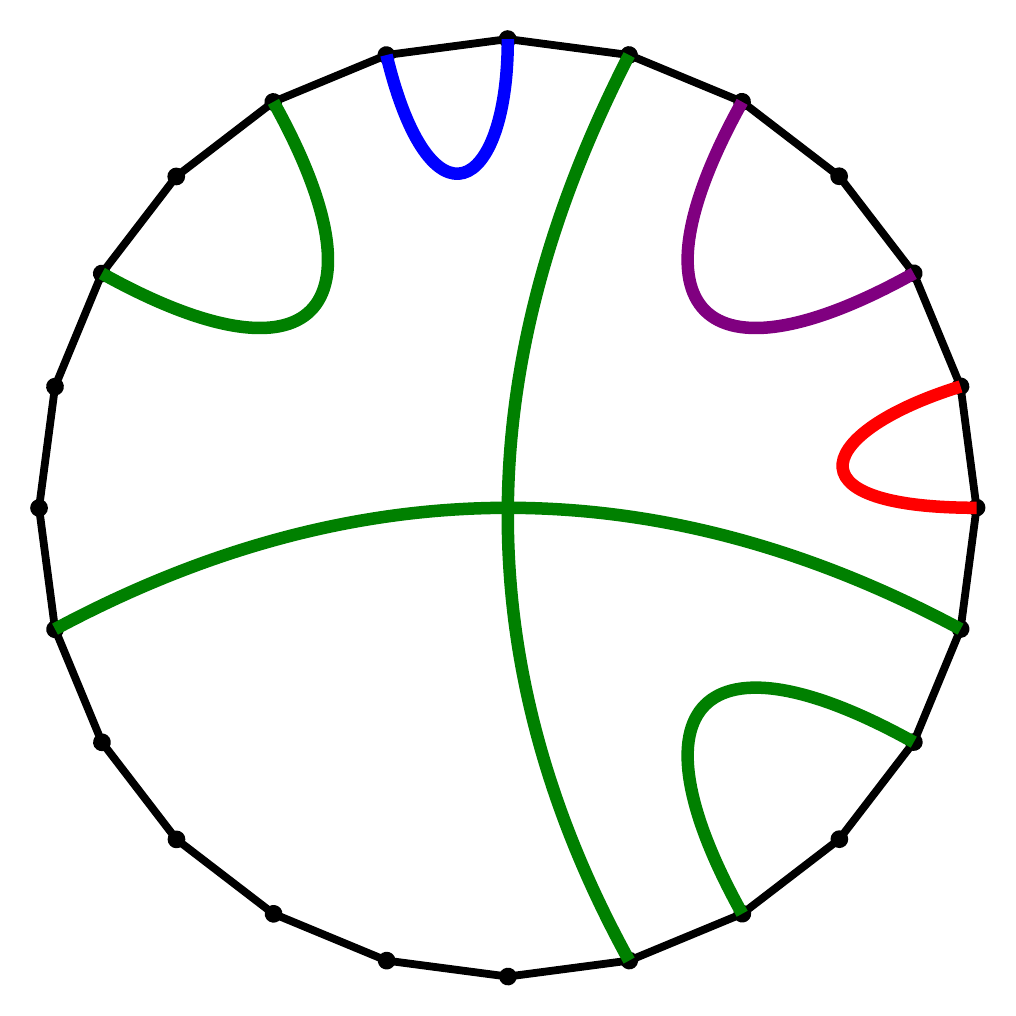}}%
    \put(-0.01241878,0.77756918){\color[rgb]{0,0,0}\makebox(0,0)[lt]{\lineheight{1.25}\smash{\begin{tabular}[t]{l}$-1$\end{tabular}}}}%
    \put(0.13593755,0.09858779){\color[rgb]{0,0,0}\makebox(0,0)[lt]{\lineheight{1.25}\smash{\begin{tabular}[t]{l}$z$\end{tabular}}}}%
    \put(0.84166256,0.1015212){\color[rgb]{0,0,0}\makebox(0,0)[lt]{\lineheight{1.25}\smash{\begin{tabular}[t]{l}$x$\end{tabular}}}}%
    \put(0,0){\includegraphics[width=\unitlength,page=2]{chordiag_frame-change_cobound.pdf}}%
    \put(0.70923079,0.02615397){\color[rgb]{0,0,0}\makebox(0,0)[lt]{\lineheight{1.25}\smash{\begin{tabular}[t]{l}$+1$\end{tabular}}}}%
    \put(0.92307691,0.19230771){\color[rgb]{0,0,0}\makebox(0,0)[lt]{\lineheight{1.25}\smash{\begin{tabular}[t]{l}$0$\end{tabular}}}}%
    \put(0.2076923,0.93230768){\color[rgb]{0,0,0}\makebox(0,0)[lt]{\lineheight{1.25}\smash{\begin{tabular}[t]{l}$0$\end{tabular}}}}%
    \put(0,0){\includegraphics[width=\unitlength,page=3]{chordiag_frame-change_cobound.pdf}}%
  \end{picture}%
\endgroup%

    \caption{
    $\textcolor{red}{\{y_j\in (z_\infty,z_0) \wedge y_{1/j} \in (x_\infty,x_0)\}}
    \triangle
    \textcolor{blue}{\{y_j\in (x_\infty,x_0) \wedge y_{1/j} \in (z_0,z_\infty)\}}
    = \textcolor{green!50!black}{N(x)\setminus \{z\}}$}
    \label{fig:Ghys_framing_coboundary}
\end{figure}
\end{proof}

\subsection{Bicolourings and bilinear forms}



\begin{Corollary}[Bicolourable]
Proposition \ref{Proposition:BasInter-H1M} implies that for every filoop $\gamma \subset M \subset F$ corresponding to a framed chord diagram $C_\varphi$ with interlace graph $G$, the following are equivalent:
\begin{itemize}[noitemsep]
    \item[$\iff$] The graph $G$ satisfies condition \ref{EveN1}: $\forall x\in V_G, \: \Card N(x)\equiv 0$.
    \item[$\iff$] In $H_1(M;\Z/2)$, the cycle $[\gamma]$ is symplectic-orthogonal to all $[\alpha_x]$.
    \item[$\iff$] In $H_1(F;\Z/2)$ the class of the filoop $\gamma$ is trivial.
    \item[$\iff$] The filoop $\gamma\subset F$ admits a \emph{bicolouring}, that is a colour map $F\setminus \gamma \to \{\bullet, \circ\}$ that changes colour when crossing $\gamma$ at smooth points, and it is unique up to global change of colours.
\end{itemize}
A graph, chord diagram or filoop satisfying these conditions will be called \emph{bicolourable}.
\end{Corollary}


\begin{Theorem}
    \label{Thm:intersection-cocycle}
    Consider a bicolourable framed chord diagram $C_\varphi$. 
    
    If $x,y\in V_G$ are not connected by an edge then $\I_\varphi(x,y)\equiv \E^2(x,y)\equiv \Card N(x)\cap N(y)$. 
    
    The cohomology class of $\I_\varphi \in Z^1(G;\Z/2)$ is independent of the framing $\varphi$. 
    
    Each cocycle in this cohomology class arises from $2^{b_0}$ framings, where $b_0=\dim H^0(G;\Z/2)$. 

    (In particular $\I_\varphi$ determines $\varphi$ up to frame inversions on connected components of $G$.)
\end{Theorem}

\begin{proof}

On the one hand, for distinct $x,y\in C$ with $\E(x,y)=0$, one may change the framing at $x$ or $y$ so that $(x_\infty,x_0)$ and $(y_\infty,x_0)$ are disjoint, and thus $\I_\varphi(x,y)=\Card N(x)\cap N(y)=\CR(x,y)$.

On the other hand Lemma \ref{Lemma:framing-cobound} implies that if two framings $\varphi,\varphi'$ differ only at $z$ then the cocycles $\I_\varphi, \I_{\varphi'} \in Z^1(G;\Z/2)$ differ by the coboundary $d\delta_z$ of the indicator function $\delta_z \colon V_G \to \{0,1\}$ of $x$.
Hence the cohomology classes $[\I_\varphi] \in H^1(G;\Z/2)$ do not depend on the framings $\varphi$.

Moreover the $\Z/2$-vector space $B^1(G;\Z/2)$ of $1$-coboundaries is generated by the $d\delta_x$ for $x\in V_G$ with one relation of the form $\sum d\delta_y$ for each connected component of $G$.
Hence in the cohomology class $[\I_\varphi]$, each cocycle arises from exactly $2^{b_0}$ framings 
(related to one another by changing orientations of all chords in unions of connected components of $G$).
\end{proof}


\begin{Proposition}[Bicolourings = framings] 
    \label{Prop:Bicolor=Framing}
    Consider a bicolourable chord diagram $C$.
    
    To a framing $\varphi \colon C\to \{\infty, 0\}$ we associate two bicolourings $\chi \colon V_G \to \{\bullet, \circ\}$ by choosing one for the corresponding filoop $F\setminus \gamma$ and colouring each intersection $x\in \gamma$ according to the region in the direction of the sum of the tangent vectors $\Vec{\gamma}(x_\infty)+\Vec{\gamma}(x_0)$. They differ by a global change in colour.

    To a bicolouring $\chi \colon V_G \to \{\bullet,\circ\}$ we associate two functions $\varphi \colon C\to \{\infty , 0\}$ by choosing the value on some letter, and extending so that consecutive letters have the same value if and only if they have different colours.
    These functions are framings and they differ by a global frame inversion.

    Seen as functions along the circle with letters of $C$, the bicolouring $\chi$ is the skew-derivative of the framing $\varphi$, and we just integrated the two skew-primitives.
    This yields a correspondence between bicolourings $\chi \colon V_G \to \{\bullet, \circ\}$ and framings $\varphi \colon C \to \{\infty,0\}$ up to global inversion.
\end{Proposition}

\begin{proof}
    We first show that the integration $\chi \mapsto \varphi$ yields framings.
    This is because between two occurrences of any letter $x\in C$, there is an even number of letters in total, hence an odd number of consecutive pairs; and an even number of alternations in colours, hence an odd number of repetitions in colours. Consequently the two occurrences of $x$ have different framings.

    Now we show that it is indeed inverse to the derivation $\varphi \mapsto \chi$.
    For this observe that while traversing the intersections of the filoop associated to $C_\varphi$, the variation of framing is related to the variation of colours as claimed. 
    In particular while constructing the filoop from the framing $(C,\varphi)$ as in \ref{Prop:framed-chordiag-loop-S}, one could consistently use only the bicolouring $(C,\chi)$ by attaching the bicoloured crosses so as to match the colours of their regions.
\end{proof}

\begin{Definition}[Bicolouring form]
Consider a graph $G$ with a bicolouring $\chi \colon V_G \to \{\bullet, \circ\}$.

The \emph{colour form} $\CB_\chi \colon V_G \times V_G \to \{0,1\}$ is defined as the differential of the indicator functions on $V_{\bullet}$ or $V_\circ$, and its restriction to the edges is called the \emph{colour coboundary} $d\chi \in B^1(G,\Z/2)$. 
\end{Definition}

\begin{Theorem}[Bilinear forms]
    \label{Thm:intercolorgraphorms}
    Consider a bicolourable framed chord diagram $C_\varphi$. 
    
    Its associated intersection form $\I_\varphi$, colour form $\CB_\chi$ and Rosenstiehl form $\CR$ are related by:
    \begin{equation*}
        \I_\varphi - \CB_\chi \equiv \E+\E^2 
    \end{equation*}
    Hence the intersection form can be computed from the bicoloured graph.
    We also recover the fact that $\I_\varphi$ determines $\CB_\chi$, whence $\varphi$ up to inversions on connected components of the interlace graph.

    In particular when $G$ satisfies \ref{EveN2}, namely $\E^2$ vanishes on the non-edges, the intersection form can be computed from the graph and its colour coboundary $d\chi \in B^1(G,\Z/2)$.
\end{Theorem}
\begin{proof}
    Note that changing the framing $\varphi$ at a chord $z$ corresponds to changing colouring $\chi$ of this chord $z$, hence to changing the colour form $d\chi$ by the coboundary $d\delta_z$.
    Hence the quantity $\I_\varphi - \CB_\chi$ is invariant by a change of framing, and we may thus choose them to alternate as reading the letters of $C$, so that all colours are equal and $\CB_\chi=0$.
    
    We must now prove that such an alternating framing satisfies $\I_\varphi = \E+\E^2$.
    For chords that are not interlaced this follows from Proposition \ref{Thm:intersection-cocycle}.
    Consider interlaced chords $x,y$, 
    and recall that:
    \begin{equation*}
    	\I_\varphi(x,y) \equiv 1+\Card \{z_j \in C_\varphi \mid z_j\in (x_\infty,x_0) \wedge z_{1/j} \in (y_0, y_\infty)\}
    \end{equation*}
    Hence must show that $\Card \{z_j \in C_\varphi \mid z_j\in (x_\infty,x_0) \wedge z_{1/j} \in (y_0, y_\infty)\} \equiv \Card N(x)\cap N(y)$, and by partitioning those sets, cancelling those in common and those with even cardinal, we want:
    \begin{align*}
        \Card \{z_j \in C_\varphi \mid z_j\in (x_\infty,y_\infty) \wedge z_{1/j} \in (y_\infty, x_0)\} + \Card \{z_j \in C_\varphi \mid z_j\in (y_\infty, x_0) \wedge z_{1/j} \in (x_0,y_0)\} 
        \\ \equiv 
        \Card \{z_j \in C_\varphi \mid z_j\in (y_\infty, x_0) \wedge z_{1/j} \in (y_0, x_\infty)\}
    \end{align*}
    These three sets yield a partition of the chords having one endpoint in the open interval $(y_\infty,x_0)$.
    Since $\varphi$ is alternating, this open interval $(y_\infty,x_0)$ contains an even number of letters, so the sum of these three cardinals is $\equiv 0$, and this concludes the proof. 
\end{proof}

\begin{Corollary}[Cohomology]
    \label{Cor:IC-cohomolog-RC}
    Consider a bicolourable chord diagram $C$ with interlace graph $G$. 
    
    The alternating framings $\varphi_\pm$ of $C$ have intersection form $\I_\varphi$ equal to the Rosenstiehl form.

    Thus all intersection cocycles $\I_\varphi$ are cohomologous to the Rosenstiehl cocycle $\CR$ in $Z^1(G;\Z/2)$. 
\end{Corollary}

\begin{Corollary}[Minimal genus]
    \label{Cor:minimal-genus-bicolor}
    Fix a bicolourable chordiagraph $G$ with Rosenstiehl form $\CR$.

    For every chord diagram $C$ with interlace graph $G$, the minimal genus of its framings equals the minimum of the half-rank over the space of bicolour forms $c_\chi\colon V_G\times V_G\to \{0,1\}$ translated by $\CR$:
    \begin{equation*}
        \ming(\E)= \tfrac{1}{2}\min \{\rank(\CR+\CB_\chi) \mid \chi \colon V_G\to \{\bullet, \circ\} \}
    \end{equation*}
\end{Corollary}


\begin{Remark}
    A bicolourable linear chord diagram $C$ corresponds to a permutation $\sigma$ pairing letters $2i$ and $2\sigma(i)+1$. Hence, there are $(2^n)n!$ bicolourable based filoops with $n$ intersections.
    
    In a future work we will relate the algebraic combinatorics of permutations (factorisation into disjoint cycles, composition law) with the algebraic topology of bicolourable filoops (connected sums).
\end{Remark}

\subsection{Lagrangian loops and spheriloops}

\begin{Corollary}[Lagrangian]
    Theorem \ref{Thm:intercolorgraphorms} implies that the following conditions are equivalent.
    
    The first and last only involve the graph $G$, while the others are valid for every framed chord diagram $C_\varphi$ with interlace $G$, corresponding to a filoop $\gamma \subset M \subset F$ with intersection form $\I_\varphi$.

    \begin{itemize}[noitemsep]
        \item[$\iff$] The graph $G$ satisfies conditions \ref{EveN1} and \ref{EveN2}: $\forall x,y\in V_G, \: \Card N(x)\cap N(y)\equiv 0$.
        \item[$\iff$] In $H_1(M;\Z/2)$, the cycles $[\gamma], [\alpha_x]$ are two by two symplectic-orthogonal.
        \item[$\iff$] In $H_1(F;\Z/2)$ the class of the loops $\gamma, \alpha_x$ are all trivial.
        \item[$\iff$] $F\setminus \gamma$ admits a bicolouring, such that non-interlaced intersections have the same colour.
        \item[$\iff$] The Resenstiehl form $\CR=\E+\E^2$ of the graph $G$ vanishes on non-edges.
    \end{itemize}

    A graph, chord diagram or loop satisfying these conditions will be called \emph{Lagrangian}.
\end{Corollary}

\begin{Corollary}[Minimal genus]
    \label{Cor:minimal-genus-Lagrange}
    Fix a Lagrangian chordiagraph $G$ with Rosenstiehl cocycle $\CR$.
    
    For every chord diagram with interlace graph $G$, the minimal genus of its framings equals the minimum of the half-rank over the space of coboundaries $B^1(G;\Z/2)$ translated by $\CR\in Z^1(G;\Z/2)$:
    \begin{equation*}
        \ming(\E)= \tfrac{1}{2}\min \{\rank(\CR+d\chi \mid \chi \in Z^0(G;\Z/2) \}
    \end{equation*} 
\end{Corollary}

\begin{Corollary}[Rosenstiehl]
    \label{Cor:Rosenstiehl}
    A chord diagram admits a framing of genus zero if and only if its interlace graph is  Gaussian.
    Such chord diagrams will thus be called Gaussian.
\end{Corollary}

\begin{Corollary}[Spheriloops]
    \label{Cor:Spherical-framings}
    A chord diagram $C$ with a Gaussian interlace graph $G$ admits $\Card H_0(G; \Z/2)=2^{b_0}$ framings of genus $0$. 
    Two of them are obtained by integrating the Rosenstiehl form as in Proposition \ref{Prop:Bicolor=Framing}, and the others differ by locally constant change of framings.
\end{Corollary}

We will construct all spheriloops with a given interlace graph in Corollary \ref{Cor:all-spheriloops-with-graph-G}.

\begin{Remark}[Variety of Lagrangian graphs]
\label{Rem:Vari-LaGraphs}
Let us work over the field $\Z/2$ with two elements.

The set $\mathcal{SE}(n)$ of simple graphs $\E$ on $n$ labelled vertices is a vector space of dimension $\binom{n}{2}$.

The subset of bicolourable graphs $\mathcal{BE}(n)$ is described by $n$ linear equations $\sum_k \E_{ik}=0$ sharing only $1$ non-trivial relation $\sum_i \sum_k \E_{ik}=0$, so it is a vector space of dimension $\binom{n-1}{2}$. 

The space of bicolourings forms a trivial vector bundle $\mathcal{BC}(n)\to \mathcal{BE}(n)$ with fibers $B^1(K_n)$ of dimension $n-1$, and of total dimension $\binom{n}{2}$. 

The subset $\mathcal{LE}(n)\subset \mathcal{BE}(n)$ of Lagrangian graphs is an algebraic subvariety described by the $\binom{n}{2}$ degenerate cubic equations $(1-\E_{ij})\sum_k \E_{ik}\E_{kj}$.

We may restrict the trivial vector bundle of colourings over Lagrangian graphs to obtain a trivial vector bundle $\mathcal{LC}(n)\to \mathcal{LE}(n)$ with fibers $B^1(G)$ of dimension $n-1$.

\end{Remark}

\begin{Remark}[Laplacian of Lagrangian graphs]
\label{Rem:Laplacian-Lagraphs}
A simple graph $G$ can be encoded by its Laplacian matrix $\Lap$, obtained by subtracting the adjacency matrix $\E$ from the diagonal matrix of degrees $\E^2(x,x)$.

The symmetric bilinear form $\DL=-\Lap+\Lap^2\colon V_G\times V_G\to \Z$ satisfies for all distinct $x,y\in V_G$:
\begin{equation*}
    \DL(x,x) = \Card N(x)^2 
    \quad \mathrm{and} \quad
    \DL(x,y) = \E(x,y)(\Card N(x) + 1 +\Card N(y)) - \Card N(x)\cap N(y)
\end{equation*}
so one may reformulate the conditions \ref{EveN1}, \ref{EveN2} and \ref{RoCo} graph classes in terms of $\DL$, in particular for bicolourable graphs we have $\DL = -\Lap+\Lap^2 \equiv \E+\E^2 = \CR$.

\end{Remark}

\section{Unique factorisations of graphs, chord diagrams and filoops}
\label{Sec:FactOperad}

For us, a graph $G$ is a relation on a finite set of vertices $V_G$ which is symmetric and nowhere reflexive.
It is represented by its adjacency matrix $\E \colon V\times V \to \{0,1\}$ which is symmetric and restricts to zero on the diagonal.
Those are often called undirected simple graphs (no loops or multiple edges).

For $U\subset V_G$, we denote by $N(U)=\{x\in V_G \setminus U\mid \exists y \in U,\; \E(x,y)=1\}$ its set of neighbours and by $G[U]$ the corresponding induced graph whose adjacency matrix is the minor $\E_{\mid U\times U}$ of $\E$.

\subsection{Split decomposition of graphs}

Let us briefly recall Cunningham's split decomposition of graphs revisited by Gioan--Paul \cite{Gio-Paul_split-decomp_2012}.

\begin{figure}[h]
 \centering
 \def\svgwidth{5cm}
 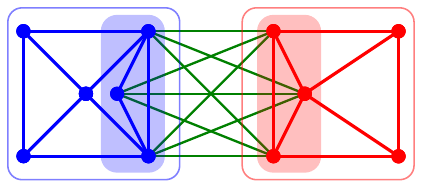
 \def\svgwidth{9cm}
 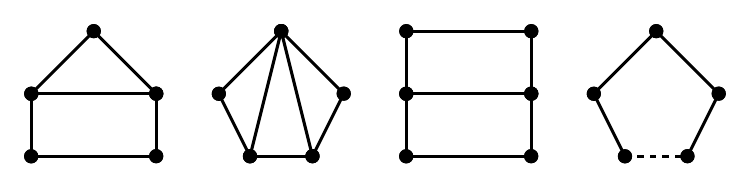
 \caption{\label{fig:split&house-gem-domino} A split of a graph. The house, gem, domino and $(n>4)$-cycles are prime.}
\end{figure}

\begin{Definition}
The \emph{toric sum} of two graphs $G_0$ and $G_1$ along non-empty subsets of their vertices $U_0\subset V_0$ and $U_1\subset V_1$, is the graph $(G_0,U_0) \boxtimes (G_1,U_1)$ obtained from the union $G_0\sqcup G_1$ by adjoining all edges of the complete bipartite graph $K(U_0,U_1)$.


Conversely, a \emph{split} of a graph $G$ is a bipartition of its vertices $V=V_0\sqcup V_1$ which is non-trivial (both $V_i$ contain at least $2$ vertices) such that the edges $E(V_0,V_1)$ connecting $V_0$ and $V_1$ are precisely those of the bipartite complete graph $K(U_0,U_1)$ on subsets of vertices $U_0\subset V_0$ and $U_1\subset V_1$.

A connected graph is called \emph{prime} when it has no splits.
In particular, it must have $>4$ vertices.
For example the house, gem, domino and $(n>4)$-cycles in Figure~\ref{fig:split&house-gem-domino} are prime.

A connected graph is called \emph{degenerate} when all its non-trivial bipartitions are splits: those are cliques $K_n$ or stars $S_{n} = K_{1,n-1}$.
In particular every graph with $<4$ vertices is degenerate.
\end{Definition}

The toric-sum factorisations (or splits) do not overlap, and form the (internal) edges of a tree: this provides a unique factorisation into a reduced graph-labelled-tree which we now introduce.

The vertices of a tree are partitioned into \emph{leaves} and \emph{nodes}, of degree $=1$ and $>1$ respectively.

\begin{Definition}
    A \emph{graph-labelled-tree} (or GLT) is a set of graphs $G_x$ indexed by the nodes of a tree $T$, with bijections $\rho_x \colon V_{x} \to N(x)$ from the vertices of $G_x$ to the neighbours of $x$ in $T$. 
    
    The \emph{accessibility graph} $G$ of such a graph-labelled-tree has vertices the leaves of $T$, and distinct vertices $u,v\in V_G$ are connected when for every node $x_i$ of the injective sequence $x_1,\dots,x_l\in V_T$ joining $u=x_0$ to $v=x_{l+1}$, the vertices of $G_{x_i}$ mapped by $\rho_{x_i}$ towards $x_{i-1}$ and $x_{i+1}$ are connected.
    

    The graph-labelled-tree is \emph{reduced} when: the nodes have degree at least $3$, every $G_x$ is prime or degenerate, and for connected nodes $x,y\in V_T$ we exclude the possibility that $G_x$ and $G_y$ are both cliques, or stars linked by a center and an extremity (as in Figure~\ref{fig:clique_star_contraction}).
\end{Definition}

\begin{figure}[h]
 \centering
 \def\svgwidth{6cm}
 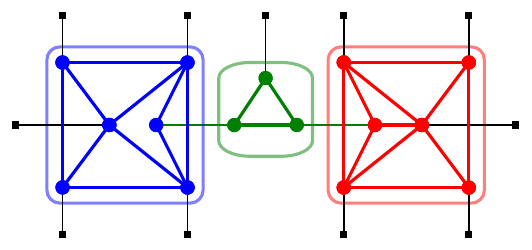
 \hfill
 \def\svgwidth{6cm}
 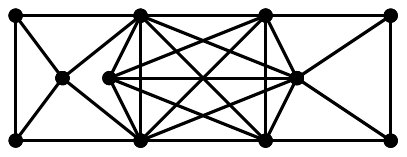
 \caption{\label{fig:acessibility} The accessibility graph of a reduced graph-labelled-tree.}
\end{figure}

The following \cite[Theorem 2.9]{Gio-Paul_split-decomp_2012} by Gioan--Paul revisits \cite[Theorem 3]{Cunningham_decomp-di-graphs_1982} of Cunningham.
\begin{Theorem}
\label{Thm:Cunningham}
A connected graph is the accessibility graph of a unique reduced graph-labelled-tree.
\end{Theorem}

\begin{Remark}[Fusion]
The reducibility conditions for a GLT ensure the uniqueness in the theorem.

Indeed, consider adjacent nodes $x,y$ in $T$.
If $G_x=K_n$ and $G_y=K_m$ are cliques then one may contract them into a clique $G_z=K_{m+n-1}$. 
If $G_x=S_n$ and $G_y=S_m$ are stars with $\rho_x^{-1}(y)$ the center of $G_x$ and $\rho_y^{-1}(x)$ an extremity of $G_y$ then one may contract them into a star $S_{n+m-1}$.
This yields a new GLT (whose tree is the contraction $T/(x,y)$) with the same accessibility graph.

\begin{figure}[h]
 \centering
 \def\svgwidth{12cm}
\begingroup%
  \makeatletter%
  \providecommand\color[2][]{%
    \errmessage{(Inkscape) Color is used for the text in Inkscape, but the package 'color.sty' is not loaded}%
    \renewcommand\color[2][]{}%
  }%
  \providecommand\transparent[1]{%
    \errmessage{(Inkscape) Transparency is used (non-zero) for the text in Inkscape, but the package 'transparent.sty' is not loaded}%
    \renewcommand\transparent[1]{}%
  }%
  \providecommand\rotatebox[2]{#2}%
  \newcommand*\fsize{\dimexpr\f@size pt\relax}%
  \newcommand*\lineheight[1]{\fontsize{\fsize}{#1\fsize}\selectfont}%
  \ifx\svgwidth\undefined%
    \setlength{\unitlength}{675bp}%
    \ifx\svgscale\undefined%
      \relax%
    \else%
      \setlength{\unitlength}{\unitlength * \real{\svgscale}}%
    \fi%
  \else%
    \setlength{\unitlength}{\svgwidth}%
  \fi%
  \global\let\svgwidth\undefined%
  \global\let\svgscale\undefined%
  \makeatother%
  \begin{picture}(1,0.13333333)%
    \lineheight{1}%
    \setlength\tabcolsep{0pt}%
    \put(0,0){\includegraphics[width=\unitlength,page=1]{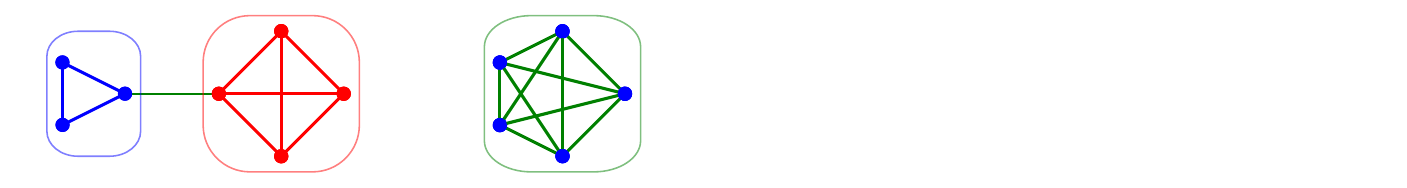}}%
    \put(0.28138027,0.05275612){\makebox(0,0)[lt]{\lineheight{1.25}\smash{\begin{tabular}[t]{l}=\end{tabular}}}}%
    \put(0,0){\includegraphics[width=\unitlength,page=2]{fusion-star-clique.pdf}}%
    \put(0.81471367,0.05275614){\makebox(0,0)[lt]{\lineheight{1.25}\smash{\begin{tabular}[t]{l}=\end{tabular}}}}%
  \end{picture}%
\endgroup%

 \caption{\label{fig:clique_star_contraction} Contraction of cliques and stars}
\end{figure}

\end{Remark}

\begin{proof}[Existence of a reduced GLT]

Every connected graph $G$ is the accessibility graph of a tree with one node decorated by $G$.
Now consider any GLT factorisation $(T,G_x)$ of $G$, and observe that every internal edge of $T$ yields a split of $G$.
If a factor $G_x$ admits a split corresponding to a factorisation $G_x=(G_0,U_0) \boxtimes (G_1,U_1)$, then define for $i=0,1$ the graphs $G'_i=(\{x_i\},\{x_i\})\boxtimes (G_i,U_i)$, and split the node decorated by $G_x$ in two nodes decorated by the $G'_i$ connected by an edge between the $u_i$, as in the left of Figure~\ref{fig:split_facto}.
Pursue these factorisation until all factors are prime or degenerate, and then fusion stars or cliques to obtain a reduced GLT factorisation.
\begin{figure}[h]
 \centering
 \def\svgwidth{4.8cm}
 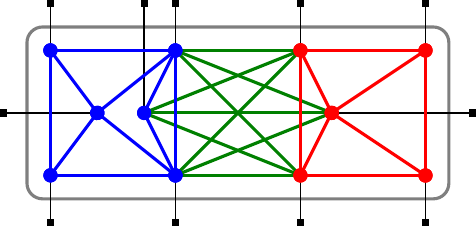
 \hfill
 \def\svgwidth{4.8cm}
 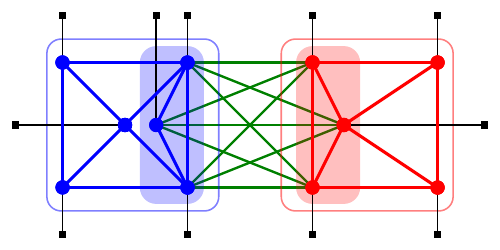
 \hfill
 \def\svgwidth{4.8cm}
 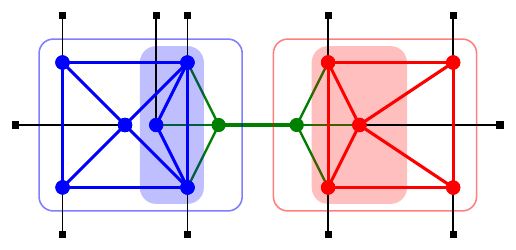
 \caption{\label{fig:split_facto} Interpreting a split-factorisation as a graph-labelled-tree with $2$ nodes.}
\end{figure}
\end{proof}

\begin{Remark}[Algorithmic]
    \label{Rem:Cunning-algo}
    The reduced GLT factorisation of a connected graph can be computed with polynomial (space \& time) complexity in its number of vertices. We refer to \cite{Charbit-Montgolfier-Raffinot_Split-decomp_2012} for the latest linear time algorithm and a short survey of the literature. 
\end{Remark}

When a connected graph has a distinguished vertex called the \emph{root}, its GLT factorisations have a distinguished leaf now called the \emph{root}. (This yields a rooted node, and a root in each factor graph.)

\begin{Definition}[Composition]
\label{Def:Composition-graphs}
For a rooted graph-labelled-tree $\TT_0$ whose leaves $f_j$ are bijected with $n$ other rooted GLTs $\TT_1, \dots,\TT_n$, we define the \emph{composition} $\TT_0(\TT_1,\dots,\TT_n)$ by grafting each root $r_j$ of $\TT_j$ on the corresponding leaf $f_j$ of $\TT_0$. 

This yields a well-defined \emph{composition} of their rooted connected accessibility graphs, which is a simultaneous toric-sum of $(G_j\setminus\{r_j\}, N(r_j))$ according to the pattern $(G_0,r_0)$.


When $\TT_0$ is a complete graph, we call this the \emph{komposition} of the $\TT_j$, and for $n=2$ we denote this operation on accessibility graphs by $(G_1,r_1)\triangle (G_2,r_2)$.
\end{Definition}

\begin{figure}[h]
 \centering
 \def\svgwidth{8cm}
 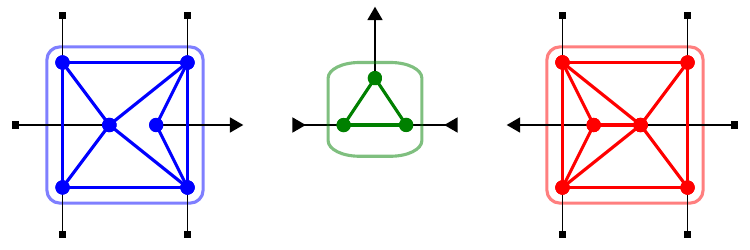
 \caption{\label{fig:operad_composition} An operad composition (which in particular a komposition).}
\end{figure}

\begin{Definition}[Complement and local complement]
    \label{Def:locomp-graph}
    Consider a set $V$ with $U\subset V$ and $c\in V$, and let us recall involutions on the set of graphs over $V$.
    Let $G$ be a graph over $V$.
    
    The \emph{complement} of $G$ on $U\subset V$ is the graph obtained by complementing the off-diagonal entries in the adjacency matrix of its induced subgraph $G[U]$.
    The \emph{local complement} of $G$ at $c\in V$ is its complement on $U=N_G(c)$ denoted $G*c$, and observe in particular that $N_{G*c}(c)=N_G(c)$.
\end{Definition}

\begin{Remark}[Complementation relations]
    \label{Rem:locomp-commute-braid}
    Consider the set $\mathcal{E}(V)$ of simple graphs over a set $V$. The local complementation involutions $*v\in \mathfrak{S}(\mathcal{E}(V))$ satisfy the variable relations:
    \begin{itemize}[noitemsep, align=left]
        \item[Commutation] If $\E(x,y)=0$ then $\E{}*x*y=\E{}*y*x$.
        \item[Braiding] If $\E(x,y)=1$ then $\E{}*x*y*x=\E{}*y*x*y$. 
    \end{itemize} 
\end{Remark}

\begin{Lemma}[Local complements and GLTs \texorpdfstring{\cite[Lemma 2.1]{Bouchet_prime-graphs-circle-graphs_1987}}{Bouchet-locomplement-Split}]
\label{Lem:GLT-Local-Comp}
    Consider a GLT $(T,G_x)$ with accessibility graph $G$.
    The local complementation of $G$ at $c\in V_G$ induces an isomorphism of $T$ and a local complementation of each factor $G_x$ accessible from $c$ at its control vertex leading to $c$.

    In particular local complementations preserve the classes of prime and degenerate graphs.
    (Note that complementing $K_n$ yields $S_n$ whereas complementing $S_n$ yields $S_n$ or $K_n$.)
\end{Lemma}

\begin{Lemma}[Subgraphs \& Connectivity]
    \label{Lem:GLT-Subgraph-Connect}
    Let $G$ be the accessibility graph of a GLT $(T, G_x)$.

    The graph $G$ is connected if and only if every factor $G_x$ is connected (\cite[Lemma 2.3]{Gio-Paul_split-decomp_2012}).

    Every factor $G_x$ appears among the induced subgraphs of $G$ (\cite[Corollary 2.6]{Gio-Paul_split-decomp_2012}).
\end{Lemma}

The last lemma enables to extend the Cunningham decomposition to disconnected graphs.

\begin{Corollary}[Forests]
We define a \emph{graph-labelled-forest} (or GLF) as a disjoint union of GLTs, and called \emph{reduced} when each of its GLT components is reduced and has connected factors.

Every graph is the accessibility graph of a unique reduced GLF.
\end{Corollary}

\subsection{Chord diagrams with a prescribed interlace graph}




\begin{Remark}[Deletion]
    For a chord diagram with chord set $V$, its subdiagram induced by $U\subset V$ is obtained by deleting the chords in $V\setminus U$.
    These operations map naturally to interlace graphs. 
\end{Remark}

\begin{Lemma}[spheric sum]
    \label{Lemma:spheric-sum-chordiag}
    For linear chord diagrams $A$ and $B$, we define their \emph{spheric sum} as the linear chord diagram $C=AB$.
    It has interlace graph $G_C = G_A \sqcup G_B$.
    
    The interlace graph of a chord diagram is disconnected if and only if it is a spheric sum of non-trivial chord diagrams.
    (The spheric sum factorisations of a chord diagram form a plane tree, whose edges correspond to the connected components of the interlace graph.)
\end{Lemma}

\begin{proof}
    %
    Suppose that a chord diagram $C$ has a disconnected interlace graph $G_C$.
    By induction, one of the connected components of $G_C$ corresponds to a sub-chord diagram $A\subset C$ consisting of a single interval. Hence denoting by $B\subset C$ the complementary interval we have $C=AB$.
    %
\end{proof}

\begin{figure}[h]
 \centering
 \def\svgwidth{15cm}
 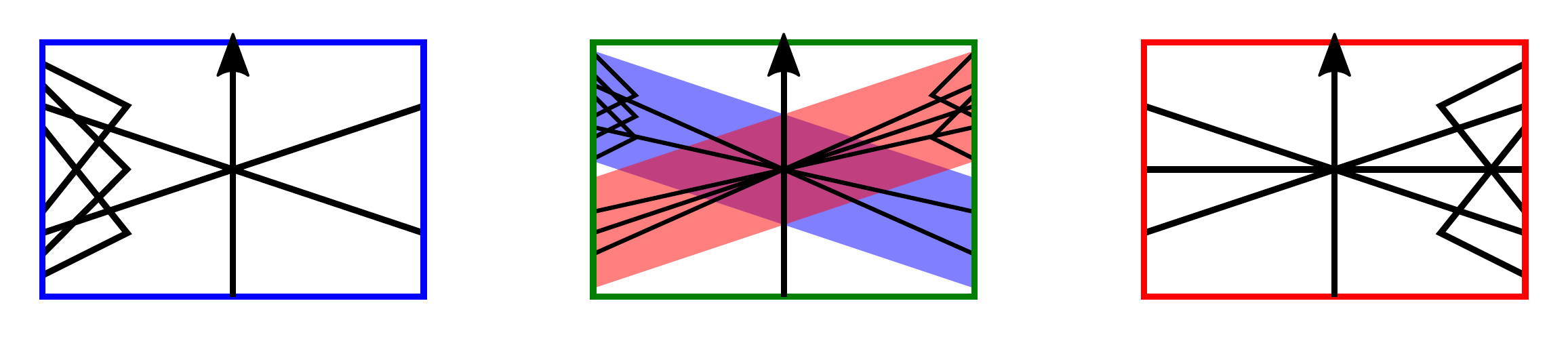
 \caption{\label{fig:cord-operad} Composition of framed rooted chord diagrams: inserting diagrams inside chords.}
\end{figure}

\begin{figure}[h]
 \centering
 \def\svgwidth{12cm}
 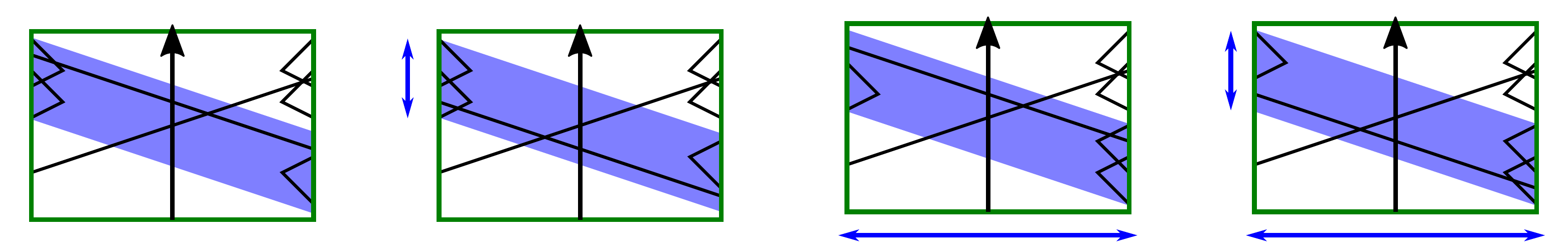
 \caption{\label{fig:cordiag_mutation} Action of $(\Z/2)^2$ by mutations of a sub-chord diagram defined by a pair of intervals.}
\end{figure}

\begin{Definition}[Composition and Mutation]
    \label{Def:Composition-chordiag}
    For a rooted framed chord diagram $C_0$ whose non-rooted chords are bijected with $n$ rooted framed chord diagrams $C_1,\dots,C_n$, we define the \emph{composition} $C_0(C_1,\dots,C_n)$ as in Figure~\ref{fig:cord-operad}, using the framings to avoid any ambiguity during insertions.
    
    Changing framings yields \emph{mutation-equivalent} framed chord diagrams (see figure \ref{fig:cordiag_mutation}).
    
    The interlace graph of such a composition is equal to the composition of the interlace graphs.
\end{Definition}

For any word $W$, we denote by $W^{-1}$ the word read backwards (regardless of a possible framing). 

\begin{Definition}[Local complement]
    \label{Def:locomp-chordiag}
    Consider unoriented chord diagrams over a set of chords $V$. The local complementation of a chord $c\in V$ acts on them by $cAcB \mapsto cAcB^{-1}$.
    
    It maps functorially to the local complementation of their interlace graphs $G\mapsto G*c$. 
\end{Definition}


\label{SubSec:chordiagraphs}

The interlace graphs of chord diagrams will be called \emph{chordiagraphs}.
The previous constructions show that they are closed by induced subgraphs, disjoint union, composition, and local complements.

\begin{Theorem}[Characterisations of chordiagraphs \cite{Bouchet_Circle-graph-obstructions_1994}]
\label{Thm:Bouchet-chordiagraph}
A graph is a chordiagraph if and only if none of its local equivalents contains an induced subgraph isomorphic to one of the following:
\begin{figure}[h]
 \centering
 \def\svgwidth{4cm}
 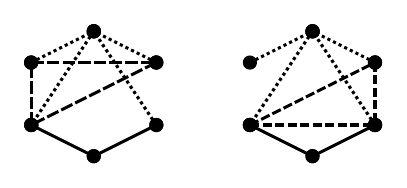
 \hspace{2cm}
 \centering
 \def\svgwidth{6cm}
 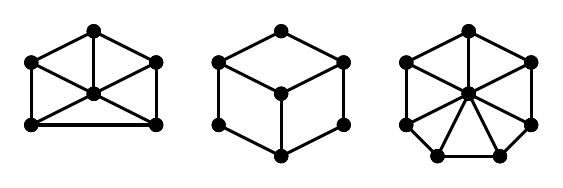\caption{\label{fig:bouchet_local_motifs}
 Local complement. Forbidden subgraphs in local equivalents of chordiagraphs.}
\end{figure}
\end{Theorem}

\begin{proof}[Comments on the proof]
By the previous sections it suffices to know which primes graphs are not chordiagraphs, and actually to give one representative in each of their local equivalence classes.
\end{proof}

Now let us explain, following \cite{Bouchet_prime-graphs-circle-graphs_1987} and \cite[section 4.8.5]{ChDuMo_Vassiliev_2012}, how Cunningham's factorisation of connected graphs provides a description of all chord diagrams with a given interlace graph.

\begin{Theorem}[Realising chordiagraphs]
    A connected chordiagraph, either prime or degenerate, is the interlace graph of a unique chord diagram up to orientation reversal. 
\end{Theorem}

\begin{proof}
    The proof is easy for degenerate graphs and the chord diagrams are depicted in Figure~\ref{fig:degenerate_cordiag_root}.
    
    For prime graphs we refer to \cite[Statement 4.4]{Bouchet_prime-graphs-circle-graphs_1987}, as well as \cite[Theorem 5]{Bouchet_Circle-graph-obstructions_1994} providing the first algorithm of polynomial (space-time) complexity which given a graph, determines if it is prime and a chordiagraph, and when so computes the unique corresponding chord diagram.
\begin{figure}[h]
 \centering
 \def\svgwidth{12cm}
 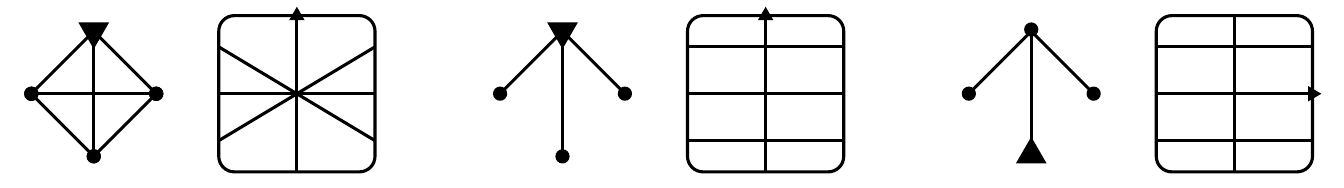
 \caption{\label{fig:degenerate_cordiag_root} Chord diagrams associated to degenerate graphs.}
\end{figure}
\end{proof}


\begin{Corollary}[Cunningham's factorisation \& Bouchet's realisation]
\label{Cor:all-filoops-with-graph-G}
Let $G$ be a connected graph. 
The chord diagrams with interlace graph $G$ are obtained as follows.

\begin{enumerate}
    \item Compute the Cunningham decomposition $(T,G_x)$ of $G$; for every factor $G_x$, check if it is a chordiagraph, and if so compute its unique unoriented chord diagram (\cite[Theorem 5]{Bouchet_Circle-graph-obstructions_1994}).
    
    \item Compose the chord diagrams according to the tree $T$ by inserting chord diagrams in the place of chords (Figure~\ref{fig:cord-operad}) in all possible ways (at most $4^{k-1}$ choices where $k=\Card(\text{nodes of }T)$).
\end{enumerate}

In particular, two (unoriented) chord diagrams have the same interlace graphs if and only if they are related by a sequence of mutations. 
    

\end{Corollary}

\subsection{Unique factorisation of filoops into spheric and toric sums}

A loop $\gamma \colon \S^1\looparrowright F_\gamma$ based at a framed self-intersection $c_j$ where $j\in \{0,\infty\}$, or equivalently pointed at a smooth point $c_j^+$ just after $c_j$, corresponds to a linear framed chord diagram $C=Ac_{1/j}Bc_{j}$.

\begin{Definition}[spheric sum]
    \label{Def:spheric-sum-loops}
    Consider loops $\alpha \colon \S^1\looparrowright F_\alpha$ and $\beta \colon \S^1\looparrowright F_\beta$ based at framed self-intersections $a_i$ and $b_j$ where $i,j\in \{0,\infty\}$, corresponding to the linear framed chord diagrams $A=Ua_{1/i}Va_{i}$ and $B=X b_{1/j} Y b_{j}$.
    Choose discs $D_\alpha\subset F_\alpha$ and $D_\beta\subset F_\beta$ intersecting the images of $\alpha$ and $\beta$ along intervals 
    just after the points $\alpha(a_i)$ and $\beta(b_j)$ in direction of $\Vec{\alpha}(a_i)$ and $\Vec{\beta}(b_j)$.

    The \emph{spheric sum} $\gamma = \alpha \otimes \beta$ is the based loop obtained from $(F_\alpha\setminus D_\alpha)\sqcup (F_\beta\setminus D_\beta)$ by identifying its oriented boundary components so as to match the cooriented pairs of points $\alpha \cap \partial D_\alpha$ and $\beta \cap \partial D_b$.
    It corresponds to the framed linear chord diagram $C=AB=Ua_{1/i}Va_{i} X b_{1/j} Y b_{j}$.

    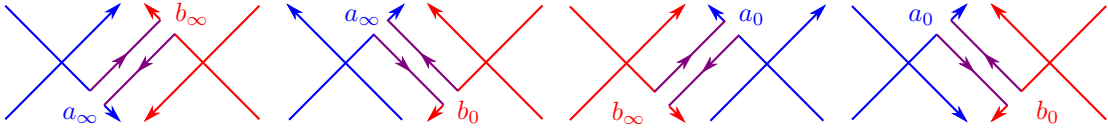
\begin{figure}[h]
        \centering
        \def\svgwidth{15cm}
\begingroup%
  \makeatletter%
  \providecommand\color[2][]{%
    \errmessage{(Inkscape) Color is used for the text in Inkscape, but the package 'color.sty' is not loaded}%
    \renewcommand\color[2][]{}%
  }%
  \providecommand\transparent[1]{%
    \errmessage{(Inkscape) Transparency is used (non-zero) for the text in Inkscape, but the package 'transparent.sty' is not loaded}%
    \renewcommand\transparent[1]{}%
  }%
  \providecommand\rotatebox[2]{#2}%
  \newcommand*\fsize{\dimexpr\f@size pt\relax}%
  \newcommand*\lineheight[1]{\fontsize{\fsize}{#1\fsize}\selectfont}%
  \ifx\svgwidth\undefined%
    \setlength{\unitlength}{1200bp}%
    \ifx\svgscale\undefined%
      \relax%
    \else%
      \setlength{\unitlength}{\unitlength * \real{\svgscale}}%
    \fi%
  \else%
    \setlength{\unitlength}{\svgwidth}%
  \fi%
  \global\let\svgwidth\undefined%
  \global\let\svgscale\undefined%
  \makeatother%
  \begin{picture}(1,0.125)%
    \lineheight{1}%
    \setlength\tabcolsep{0pt}%
    \put(0,0){\includegraphics[width=\unitlength,page=1]{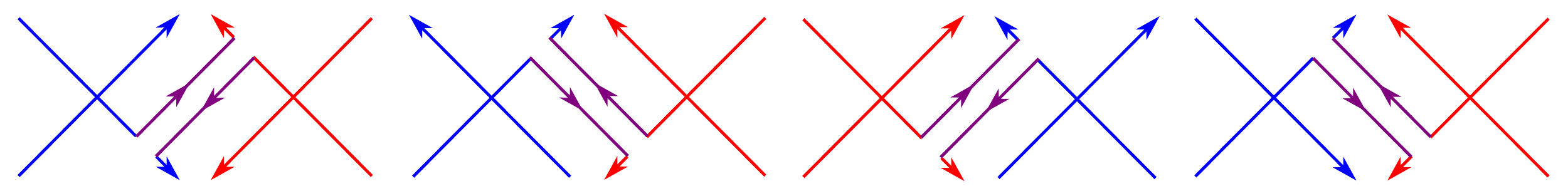}}%
    \put(0.0625,0.0125){\color[rgb]{0,0,1}\makebox(0,0)[lt]{\lineheight{1.25}\smash{\begin{tabular}[t]{l}$a_\infty$\end{tabular}}}}%
    \put(0.66200247,0.09979914){\color[rgb]{0,0,1}\makebox(0,0)[lt]{\lineheight{1.25}\smash{\begin{tabular}[t]{l}$a_0$\end{tabular}}}}%
    \put(0.81250002,0.1){\color[rgb]{0,0,1}\makebox(0,0)[lt]{\lineheight{1.25}\smash{\begin{tabular}[t]{l}$a_0$\end{tabular}}}}%
    \put(0.31250001,0.1){\color[rgb]{0,0,1}\makebox(0,0)[lt]{\lineheight{1.25}\smash{\begin{tabular}[t]{l}$a_\infty$\end{tabular}}}}%
    \put(0.92544733,0.01313902){\color[rgb]{1,0,0}\makebox(0,0)[lt]{\lineheight{1.25}\smash{\begin{tabular}[t]{l}$b_0$\end{tabular}}}}%
    \put(0.16249999,0.1){\color[rgb]{1,0,0}\makebox(0,0)[lt]{\lineheight{1.25}\smash{\begin{tabular}[t]{l}$b_\infty$\end{tabular}}}}%
    \put(0.5495025,0.01229916){\color[rgb]{1,0,0}\makebox(0,0)[lt]{\lineheight{1.25}\smash{\begin{tabular}[t]{l}$b_\infty$\end{tabular}}}}%
    \put(0.41249996,0.0125){\color[rgb]{1,0,0}\makebox(0,0)[lt]{\lineheight{1.25}\smash{\begin{tabular}[t]{l}$b_0$\end{tabular}}}}%
  \end{picture}%
\endgroup%

        \caption{The $(i,j)$-spheric sums of loops for all $(i,j)\in \{0,\infty\}^2$.}
        \label{fig:spheric-sums-loops}
    \end{figure}
\end{Definition}

\begin{Remark}[Genus]
    A spheric-sum of pointed filoops $\alpha$ and $\beta$ has $g(\alpha\otimes \beta)=g(\alpha)+g(\beta)$.
\end{Remark}

\begin{Corollary}[Spheric-sum factorisation]
    \label{Cor:Spheric-sum-facto}
    A filoop admits a unique spheric-sum factorisation into filoops with connected interlace graphs: it has the structure of a plane tree, which can be read off its chord diagram, whose vertices correspond to the connected components of its interlace graph.
\end{Corollary}
\begin{proof}
    The spheric-sum of linear framed chord diagrams descends to linear chord diagrams as in \ref{Lemma:spheric-sum-chordiag}. Hence a filoop is a non-trivial spheric-sum if and only if its interlace graph is disconnected.
\end{proof}

\begin{Definition}[Plumbing]
    \label{Def:plumbing-loops}
    Consider a pair of filoops $\alpha \colon \S^1\looparrowright F_\alpha$ and $\beta \colon \S^1\looparrowright F_\beta$ based at self-intersections $a$ and $b$, corresponding to framed rooted chord diagrams $a_\infty U a_0 V$ and $b_\infty X b_0 Y$. 
    
    
    Their \emph{plumbing} is the filoop $\alpha \asymp \beta$ defined by the framed chord diagram $UX^{-1}V^{-1}Y$, whose topological construction is depicted in Figure~\ref{fig:plumbing-loops}.
    
    It is obtained by removing disc-neighbourhoods of the based points $D_a \subset F_\alpha$ and $D_b \subset F_\beta$ and attaching the boundary components of $(F_\alpha \setminus D_\alpha) \sqcup (F_\beta \setminus D_\beta)$ by a cylinder which contains two non-intersecting segments connecting the pairs of endpoints $(\alpha \cap \partial D_\alpha)\sqcup (\beta\cap \partial D_\beta)$.

\begin{figure}[h]
    \centering
    \def\svgwidth{12cm}
    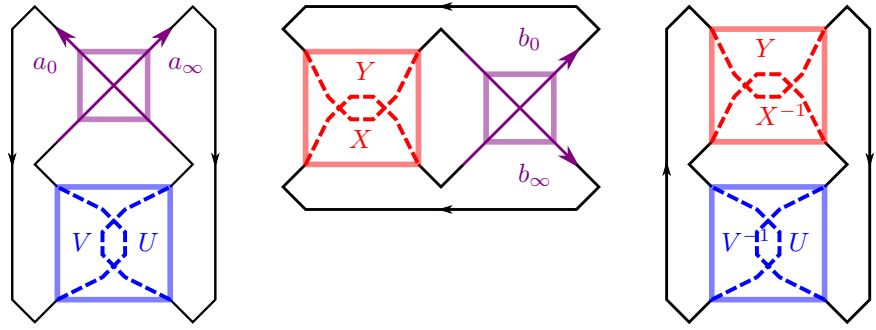
    \caption{Based filoops $\alpha$, $\beta$ and their plumbing $\alpha \asymp \beta$. (The dotted lines may have intersections.)}
    \label{fig:plumbing-loops}
\end{figure}
\end{Definition}

\begin{Definition}[Smoothing, Trivial, Essential]
    The plumbing of framed rooted chord diagrams $A=a_\infty U a_0 V$ and $B =b_\infty X b_0 Y$ will be called:
    \begin{enumerate}[noitemsep]
        \item[] \emph{smoothing} $\alpha$ at $a$ when $XY$ has $0$ chords, and denoted $\alpha \asymp a$
        \item[] \emph{trivial} when $UV$ or $XY$ has $1$ chord
        \item[] \emph{essential} when the roots $a$ and $b$ have $\ge 1$ interlaced neighbours.
    \end{enumerate}
    Note that a non-essential plumbing can be expressed using (one or two) spheric sums, and conversely a spheric sum can be expressed with a non-essential plumbing.
\end{Definition}

\begin{Remark}[Genus]
        An essential plumbing of based filoops $\alpha$ and $\beta$ has $g(\alpha \asymp \beta)=g(\alpha)+g(\beta)$.
\end{Remark}

\begin{Lemma}[Plumbing interlace graphs]
    Consider a pair $(A,a)$ and $(B,b)$ of (framed) rooted chord diagrams, whose interlace graphs $(G_A,a)$ and $(G_B,b)$ have rooted GLTs $(\TT_A,a)$ and $(\TT_B,b)$.

    The smoothing $\gamma \asymp a$ has interlace graphs $G_A\asymp a = (G_A*a)\setminus a$.
    Its GLTs are obtained by complementing the rooted GLTs $(\TT_A,a)$ and deleting the root.
    
    The plumbing $\alpha \asymp \beta$ has interlace graph $(G_A,a) \asymp (G_B,b)=(G_A\asymp a, N(a)) \boxtimes (G_B\asymp b, N(b))$. 
    Its GLTs are obtained by complementing the rooted GLTs $(\TT_A,a)$ and $(\TT_B,b)$ and grafting the results along their roots (or composing them into a rooted tripod $(K_3,c)$ and deleting the root).
\end{Lemma}

\begin{Proposition}[Plumbing and splits]
    Consider a filoop $\gamma$ with a connected interlace graph $G$.
    
    The filoop $\gamma$ is a non-trivial essential skew-plumbing if and only if $G$ admits a split.
\end{Proposition}

\begin{proof}
    If $\gamma$ is a skew plumbing of $\alpha$ and $\beta$, then a GLT for its interlace graphs is obtained from the rooted GLTs of $\alpha$ and $\beta$ by complementing and grafting them at the roots. By Lemma \ref{Lem:GLT-Local-Comp}, local complementation preserves the trees, the grafting yields a split of a connected graph provided the skew plumbing is essential.

    Conversely if $G$ has a split, then by Corollary \ref{Cor:all-filoops-with-graph-G} the chord diagram $C$ associated to $\gamma$ is of the form $UXVY$ where $UV$ and $XY$ are sub-chord diagrams having $\ge 2$ chords, with chords between $U \& V$ and between $X\& Y$ (by definition of a split for the connected graph $G$).
    Hence the framed chord diagram $C_\varphi$ is a skew-plumbing of $a_\infty U a_0 V^{-1}$ and $b_\infty X^{-1} b_0 Y$.
\end{proof}

\begin{Theorem}[Plumbing factorisation]
    \label{Thm:unifacto-filoops}
    Consider a filoop $\gamma$ with a connected interlace graph $G$. 
    
    A split factorisation of $G$ into a graph-labelled-tree $(T,G_x)$ yields an essential-plumbing factorisation of $\gamma$ into a filoop-labelled-tree $(T,\gamma_x)$.

    The interlace graph of $\gamma_x$ is the local-complement of $G_x$ at all its vertices in any order, which will prescribe the order in which one must perform the plumbings associated to its incident tree-edges. 
    
    The genus of $\gamma$ is the sum of genera of the $\gamma_x$.
\end{Theorem}

\begin{Corollary}[Mutations]
    Consider a filoop $\gamma$.
    For a plumbing factorisation $\gamma=\alpha \asymp \beta$ which is neither smoothing nor trivial, we define its \emph{mutations} as the filoops $\alpha \asymp \beta'$, $\alpha' \asymp \beta$, $\alpha' \asymp \beta'$ obtained by inverting the orientation and framing of $\alpha$ or $\beta$.
    
    The mutation is called \emph{essential} or \emph{non-essential} in accordance with the plumbing factorisation.
    Its non-essential mutations are given by the edges in its spheric-sum factorisation tree.
    Its essential mutations are given by the edges of the trees appearing in its essential-plumbing factorisation.
\end{Corollary}

\section{Generating grammars for Gaussian graphs and spheriloops}
\label{Sec:GLT-Gaussian}

\subsection{Spheriloops with a prescribed interlace graph}

Let us combine Corollary \ref{Cor:all-filoops-with-graph-G} and \ref{Cor:Spherical-framings} to describe all spheriloops corresponding to a given connected Gaussian chordiagraph.

\begin{Corollary} 
\label{Cor:all-spheriloops-with-graph-G}
Consider a connected Gaussian chordiagraph $G$. The set of spheriloops whose chord diagrams have interlace graph $G$ are obtained as follows.

\begin{enumerate}
    \item Compute the Cunningham decomposition $(T,G_x)$ of $G$, and realise every (prime or degenerate) factor by a (unique) unoriented chord diagram.
    
    \item Compose the chord diagrams according to the tree $T$ in all possible ways (at most $4^{k-1}$ choices where $k$ is the number of internal nodes of the tree).
    
    \item Compute the two framings of genus $0$, by integrating the Rosenstiehl form as in Corollary \ref{Cor:Spherical-framings} 

    \item Construct the loop from the framed chord diagram as in Proposition \ref{Prop:framed-chordiag-loop-S}.
\end{enumerate}

In particular, if a connected Gaussian chordiagraph is prime or degenerate, then it corresponds to a unique loop in the sphere up to homeomorphism.
\end{Corollary}

\begin{Remark}[Minimal genus framings of bicolourable graphs]
    We may apply the algorithm in Corollary \ref{Cor:all-spheriloops-with-graph-G} to a connected chordiagraph $G$ satisfying only \eqref{EveN1}, and adapt step 3 using \ref{Cor:minimal-genus-bicolor} to find all (bicolourable) loops of minimal genus with interlace graph $G$.
    The step 3 would then consist in computing the rank of all matrices $\CR+\CB$ for $\CB$ ranging among colouring forms of the graph.
\end{Remark}

\begin{Example}[From Gaussian chordiagraphs to spheriloops]
\label{Ex:Spheriloop-from-GaussChordiagraph}
The images \ref{fig:AbCdEfBaGeDhIcFgHi_graph} and \ref{fig:Intro:AbcDeFbagHdCiGfEhI_graph} show connected Gaussian graphs (which one may construct dynamically according to Theorem \ref{Thm:Construct-Gauss-GLT}), and their reduced GLT decompositions (one may automatise this following Remark \ref{Rem:Cunning-algo}).

The vertices of the accessibility graphs have been coloured by integrating the Rosenstiehl cocycle, and this colouring is used to find the framing of the chord diagram at the next step.

The images \ref{fig:AbCdEfBaGeDhIcFgHi_chordialoop} and \ref{fig:Intro:AbcDeFbagHdCiGfEhI_chordialoop} show their essentially unique associated framed chord diagrams computed as in Corollary \ref{Cor:all-spheriloops-with-graph-G}, and the corresponding spheriloops drawn according to Proposition \ref{Prop:framed-chordiag-loop-S}.


\begin{figure}[h]
 \centering
 \def\svgwidth{5cm}
 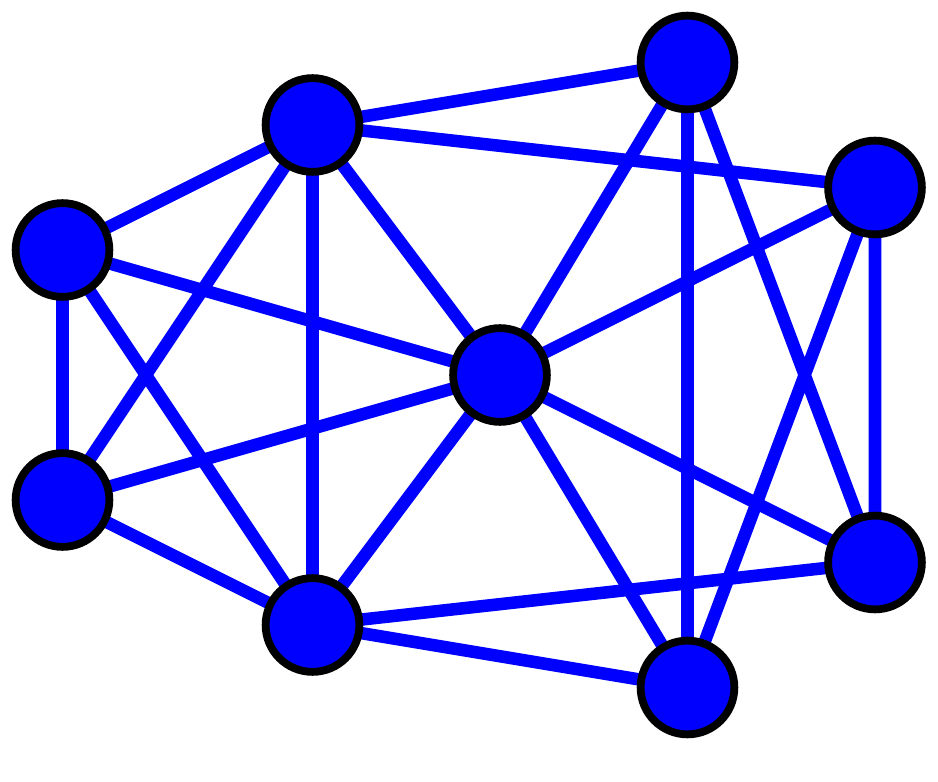
 \hfill
 \def\svgwidth{8cm}
 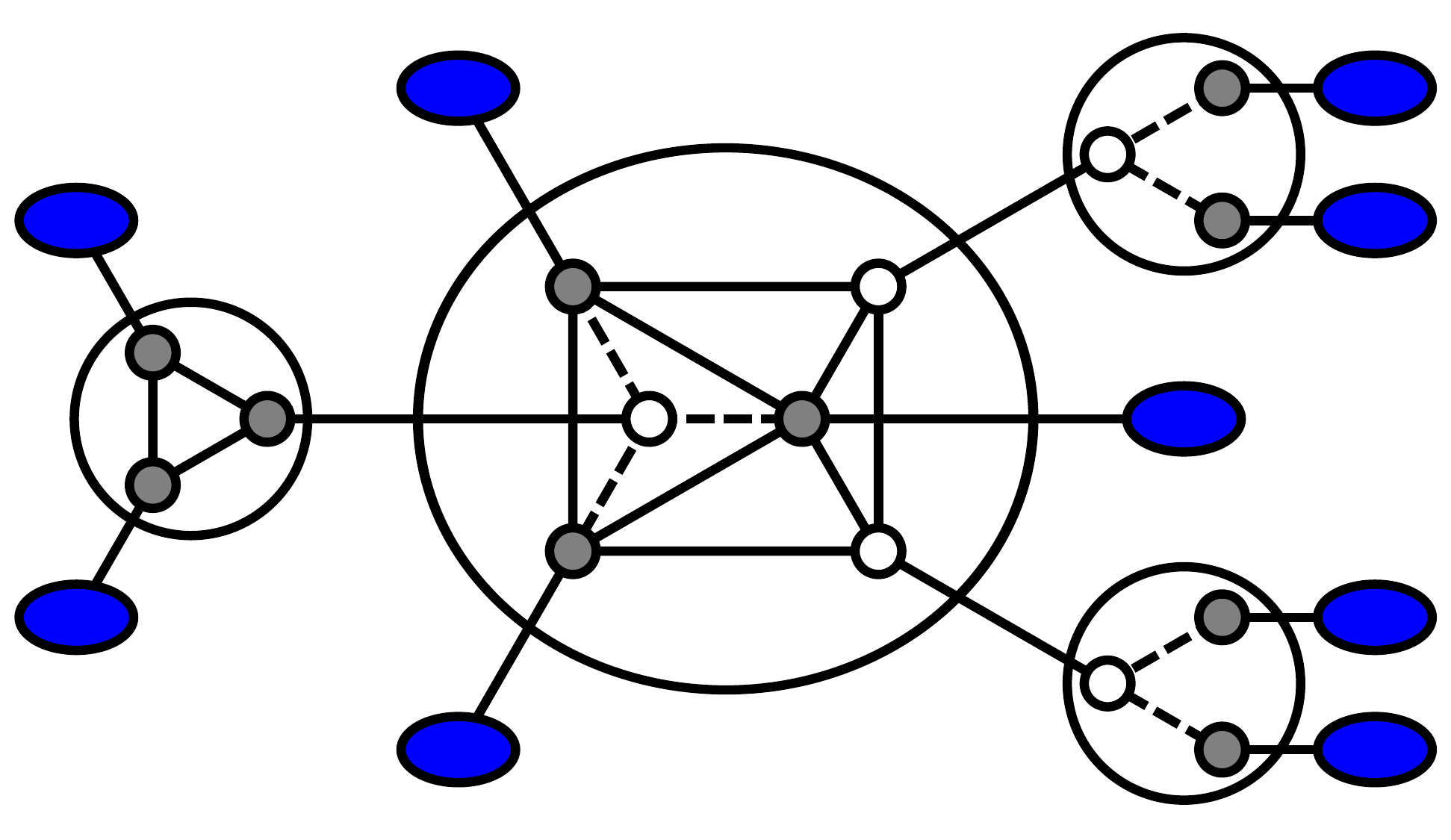\caption{\label{fig:AbCdEfBaGeDhIcFgHi_graph} A connected Gaussian chordiagraph and its reduced GLT factorisation.}
\end{figure}

\begin{figure}[h]
 \centering
 \def\svgwidth{4cm}
 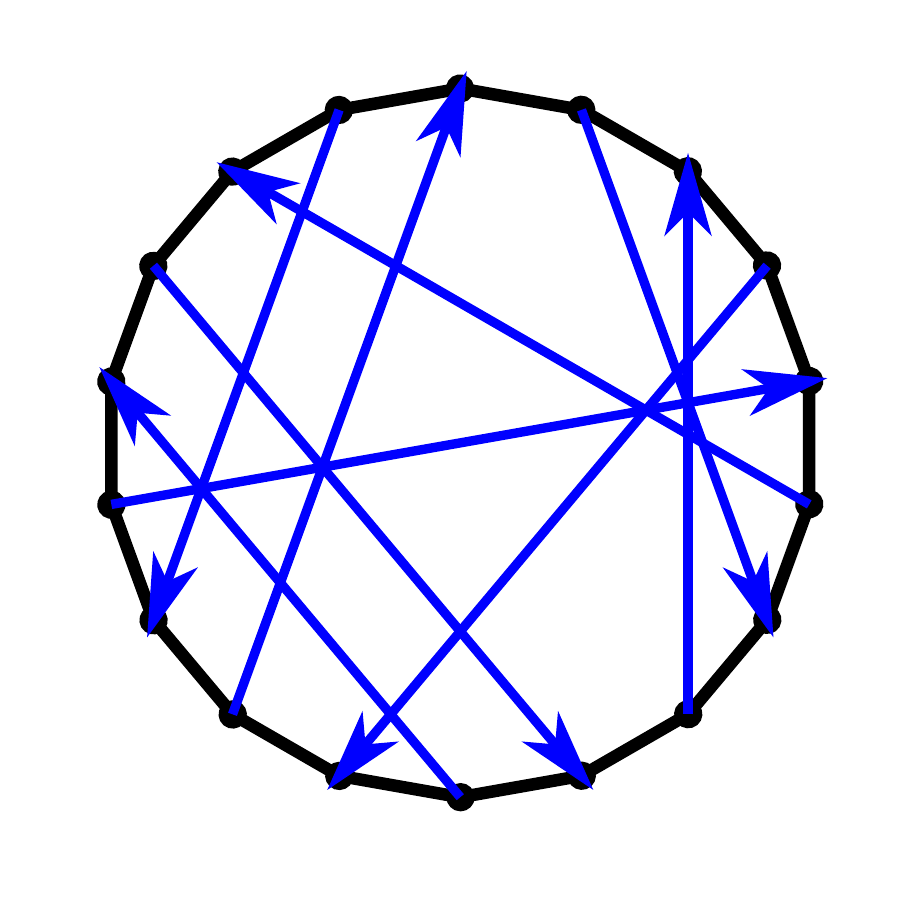
 \hspace{2cm}
 \def\svgwidth{7cm}
 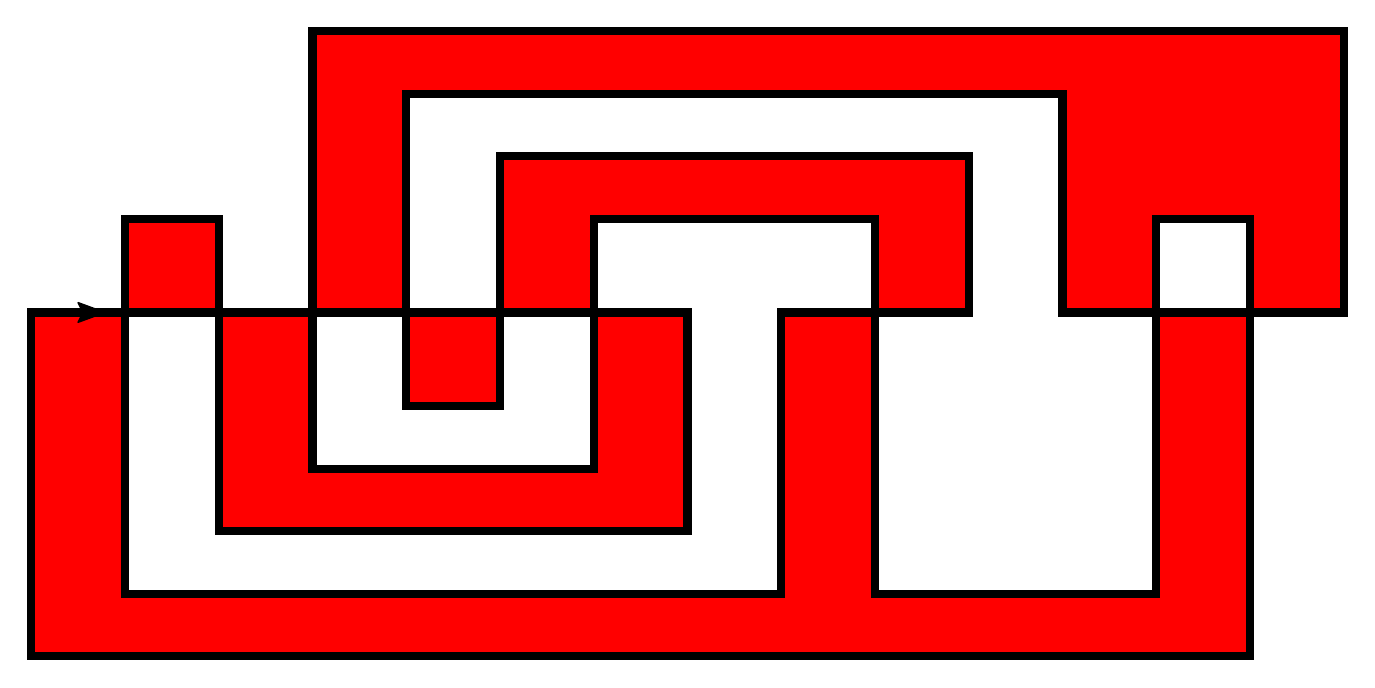\caption{\label{fig:AbCdEfBaGeDhIcFgHi_chordialoop} The unique framed chord diagram and its corresponding spheriloop.}
\end{figure}
\end{Example}


\subsection{The generating grammar of Gaussian graphs}


\begin{Definition}[Weights]
\label{Def:weights-GLT}
Consider a graph-labelled-tree $\TT=(T, G_x)$, and let us assign weights to the vertices of its factors $G_x$, whose values will mostly import $\bmod{2}$.

A leaf $v$ of $T$ is \emph{immediately accessible} from a vertex $u$ of the factor $G_{x_0}$ when for every node $x_i$ of the injective sequence $x_0,\dots,x_l\in V_T$ joining $x_0$ to $v=x_{l+1}$, the vertices of $G_{x_i}$ mapped by $\rho_{x_i}$ towards $x_{i-1}$ and $x_{i+1}$ are connected (we immediately leave $G_{x_0}$ through the edge $\rho_{x_0}(u)=(x_0,x_1)$).

The \emph{weight} of the vertex $u$ in $G_{x_0}$ is the number of leaves in $T$ immediately accessible from $u$.




\end{Definition}

\begin{Definition}[Weighted Rosenstiehl]
\label{Def:weight-Rosenstiehl}
We define the weighted conditions \ref{EveN1}, \ref{EveN2} and \ref{RoCo}.

Consider a simple graph $\E\colon V\times V \to \{0,1\}$ with a diagonal matrix $\W \colon V\times V \to \Z$ of weights (whose entries will mostly import $\bmod{2}$).
For $U\subset V$ denote $\W(U)=\sum \asrt{u\in U}\W(u)$.

We re-define the symmetric bilinear form $\CR \colon V_G\times V_G \to \Z$ as $\CR=\E+\E \W\E$ so that:
\begin{equation*}
    \forall x,y\in V_G, \quad
    \CR(x,y) = \E(x,y) + \W( N(x)\cap N(y)) 
\end{equation*}
and the conditions \ref{EveN2}, \ref{RoCo} remain the same as in Definition \ref{Def:Graph-Conditions}, 
whereas \ref{EveN1} becomes: 
\begin{equation*}
    \forall x\in V_G\colon \quad \W(x)=1 \implies \E^2(x,x)\equiv 0
\end{equation*}
where the equality $\W(x)=1$ is to hold in the working domain (namely $\Z$ or $\Z\bmod{2}$).

A graph satisfying the $\Z/2$-weighted conditions \ref{EveN2}, \ref{RoCo} will be called \emph{CL2}.

A CL2 graph satisfying the $\Z$-weighted condition \ref{EveN1} will be called \emph{CL2L1}.
\end{Definition}

\begin{Example}[Decorations]
    In the figures of example \ref{Ex:Spheriloop-from-GaussChordiagraph}, the factors labelling the nodes of the trees have decorated vertices and edges according to Definitions \ref{Def:weights-GLT} and \ref{Def:weight-Rosenstiehl}.
\end{Example}


\begin{Remark}[Algebraic variety of CL2-graphs]
\label{Rem:Vari-CL2-Graphs}
Let us work over the field $\Z/2$.
The vector space $\mathcal{WE}(n)$ of weighted simple graphs $(\W,\E)$ on $n$ vertices is a vector space of dimension $n+\binom{n}{2}=\binom{n+1}{2}$ containing the algebraic subset $\mathcal{CL}_2(n)$ of CL2 graphs described by the equations \eqref{EveN2}, \eqref{RoCo}.

The projection $\mathcal{WE}(n)\to \mathcal{E}(n)$ restricts to a linear fibration $\mathcal{CL}_2(n)\to \mathcal{E}(n)$.
Indeed, for a fixed graph $\E$, the \ref{EveN2} conditions yield $\binom{n}{2}-\Card(E)$ linear equations in $\W$, and the \ref{RoCo} conditions yield $b_1 = \dim H_1(G)$ linear equations in $\W$.
Since $b_1=b_0+\Card(E)-n$, this is a total of $\binom{n-1}{2}+b_0$ equations, which may share many relations depending on $\E$.

As in Remark \ref{Rem:Vari-LaGraphs}, for fixed weights $\W$, condition \ref{EveN2} yields $\binom{n}{2}$ singular cubic equations in $\E$, the condition \ref{RoCo} yields $\sum_{3}^{n} \binom{n}{k}(k-1)! < 2n^{n-1}$ polynomial equations $\E$ of degree $\le n+2$. 

\end{Remark}

\begin{Example}[CL2 graphs]
    Let us determine the CL2 weightings on various graphs.

    On a graph with no edges, all weightings are CL2. This applies in particular a single vertex $K_1$.

    On a tree, the CL2 weightings are those assigning weight $0$ to all internal nodes. This applies in particular to the edge $K_2$, and the star $S_{1,n}$.

    On the complete graph $K_n$ with $n>2$ vertices, the CL2 weightings $\W$ are those for which the total weight $W=\sum \W(x)$ is $\equiv 1$. 
    Indeed we have $\CR(x\sim y)=W-\W(x)-\W(y)+1$ and its integral along a cycle is $\equiv (W+1)\times \len(cycle)$. This applies in particular to the $3$-cycle $K_3$.

    A cycle of length $n>3$ admits a CL2 weighting only when $n$ is even, in which case $\W \equiv 0$. 
    Indeed a CL2 weighting must assign weight $0$ to all vertices by \ref{EveN2}, whence $\CR(x\sim y) \equiv 1$ so we must have $n\equiv 0$ by \ref{RoCo}.

    On a bipartite graph (namely whose cycles have even length), the weighting $\W \equiv 0$ satisfies CL2.

    The house graph and the gem graph do not admit any CL2 weightings. 
    The domino graph, which is bipartite, admits a unique CL2 weighting, namely $\W \equiv 0$.

    Among the graphs in Figure~\ref{fig:bouchet_local_motifs} characterising the (non) chordiagraphs, only the bipartite one admits a CL2 weighting and it is unique, namely $\W\equiv 0$.
\end{Example}


\begin{Proposition}
    \label{Prop:GLT-Gauss}
    Consider a GLT factorisation $(T, G_x)$ of a connected graph $G$. 
    \begin{enumerate}
        \item[] $G$ satisfies \ref{EveN1} if and only if all the $\Z$-weighted-factors $G_x$ satisfy \ref{EveN1}.
        \item[] $G$ satisfies \ref{EveN2} if and only if all the $\Z/2$-weighted-factors $G_x$ satisfy \ref{EveN2}.
        \item[] $G$ satisfies \ref{RoCo} if and only if all the $\Z/2$-weighted-factors $G_x$ satisfy \ref{RoCo}.
    \end{enumerate}
     Consequently if $G$ is bicolourable, then it is Gaussian if and only if all factors $G_x$ are CL2.
\end{Proposition}

\begin{proof}

First recall from Lemma \ref{Lem:GLT-Subgraph-Connect} that $G$ is connected if and only if every factor $G_x$ is connected, and that every factor $G_x$ appears among the induced subgraphs of $G$.


Now consider two leaves of $(T,G_x)$. They are connected by a unique geodesic $x_1,\dots,x_l$, which singles out a pair of vertices $u_j, v_j$ in each graph $G_{x_j}$ decorating its nodes $x_j$.

Suppose that the two leaves are disconnected.
If there is more than one disconnecting node then they have no common neighbours.
If there is a unique disconnecting node $x_j$, then their common neighbours are the leaves immediately accessible from the common neighbours of the vertices $u_j, v_j$ in $G_{x_j}$, and their cardinal is $\CR_j(u_j,v_j)$.
This proves the claim for the condition \ref{EveN2}.

Suppose that the leaves are connected.
Their common neighbours are the disjoint union over the nodes $x_j$ along the geodesic that connects them of the leaves that are immediately accessible from the common neighbours of the vertices $u_j, v_j$ in $G_{x_j}$, and their cardinal is $\sum_j \CR_j(u_j,v_j)$.

Finally consider a cycle of the accessibility graph $G$, and observe that in the GLT it decomposes as a connected sum of local cycles in the graphs $G_{x_j}$.
The condition \ref{RoCo} for the global cycle in $G$ is thus equivalent to the sum of the conditions \ref{RoCo} for the local cycles in the weighted graphs $G_{x_j}$.
\end{proof}

The following theorem, which is a Corollary of \ref{Thm:Cunningham} and \ref{Prop:GLT-Gauss}, says that connected Gaussian graphs form a language generated by an arborescent grammar which is context-sensitive and unambiguous:
    \begin{itemize}[noitemsep]
        \item[$\circ, \circ$] non-terminal symbols are leaves in CL2-weighted (prime or degenerate) graphs,
        \item[$\circ - \circ$] branching rules are constrained by matching weights and internal degrees, 
        \item[$\bullet - \bullet$] terminal symbols are (prime or degenerate) CL2L1-weighted graphs.
    \end{itemize}
    Moreover this generation is unique when all primitive and terminal graphs are prime or degenerate.

\begin{Theorem}[Grammatical generation of Gaussian GLT]
\label{Thm:Construct-Gauss-GLT}
    The following dynamical generation algorithm yields all and only connected Gaussian graphs:
    \begin{enumerate}
        \item[Atoms] For every (prime or degenerate) graph $G$, compute its CL2-weightings $\bmod{2}$ and construct a weighted GLT whose tree has only one node decorated by $G$.

        \item[Bonds] Choose weighted GLT $A$ and $B$ from those already constructed with leaves $f_a$ and $f_b$ connected to control vertices $u_a\in A_x$ and $u_b\in B_y$ such that the weight of $u_a$ is equal to the degree of $f_b$ and the weight of $f_a$ is equal to the degree of $u_b$, and form the grafting $C = K_2((A,f_a),(B,f_b))$.
        
        \item[Molecules] Among the weighted GLT so constructed, those for which all leaves are connected to control vertices of odd weight and even degree are Gaussian.
    \end{enumerate}
    When restricting to prime and degenerate graphs in the first step, every Gaussian graph is generated only once (possibly up to automorphisms depending on the implementation).
\end{Theorem}

\begin{Remark}[Generating spheriloops]
Adapting the first step in algorithm \ref{Thm:Construct-Gauss-GLT} by restricting to chordiagraphs (using Bouchet's algorithmic description \cite{Bouchet_Circle-graph-obstructions_1994}), yields all Gaussian chordiagraphs.

Combining this with Corollary \ref{Cor:all-spheriloops-with-graph-G} yields all spheriloops.

\end{Remark}

\begin{Corollary}[Constructing Gaussian graphs]
The following construction yields a GLT factorisation of a Gaussian graph from any CL2-weighted graph.
    \begin{enumerate}
        \item Choose a graph $G_0$ with a CL2 weighting and construct a graph-labelled-tree whose tree has only one node and is decorated by $G_0$
        \item For every leaf connected to a vertex of even weight, compute the total weight of its neighbours: if it is odd then attach a clique $K_{3}$ and if it is even then attach a star $S_{1,2}$ by its center.
    \end{enumerate}
If $G_0$ has $n$ vertices and the CL2-weighting has $k$ even weights then the result has $n+2K$ vertices.

If $G_0$ is already Gaussian then we recover the graph $G_0$, so all Gaussian graphs appear trivially from this construction.
However if we impose $G_0$ to be prime then only the prime Gaussian graphs appear in this way.

One may apply this procedure to families of graphs $G_0$ with CL2-weightings having even weights (such as trees or bipartite graphs), to construct families of genuine Gaussian graphs.
\end{Corollary}



\begin{Example}[Smallest Gaussian non-chordiag graph]
    Consider the graph $(\{c\},\{c\})\boxtimes (C_3,C_6) $ in the middle of Figure~\ref{fig:bouchet_local_motifs}, which cannot appear as an induced subgraph of a chordiagraph. 
    
    It is bipartite so the constant $\equiv 0$ weighting is Gaussian, and the previous algorithm yields a graph on $14$ vertices which is Gaussian, but not a chordiagraph.

\begin{figure}[h]
    \centering
    \def\svgwidth{6cm}
    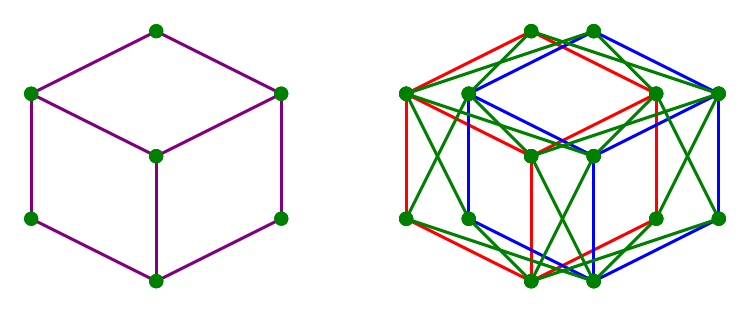
    \caption{The smallest Gaussian non-chordiag graph}
    \label{fig:Gauss-graph-non-chordiag}
\end{figure}
\end{Example}

\subsection{Smoothing and plumbing Gaussian graphs}

\begin{Definition}[Smoothing graphs]
For a graph $G$, its \emph{smoothing} at $v\in V_G$ is $G\asymp v=(G*v)\setminus v$.    
\end{Definition}


We know that Gaussian chordiagraphs are stable under smoothings.
We now show that all Gaussian graphs are stable under smoothings.

\begin{figure}[h]
    \centering
    \def\svgwidth{5cm}
    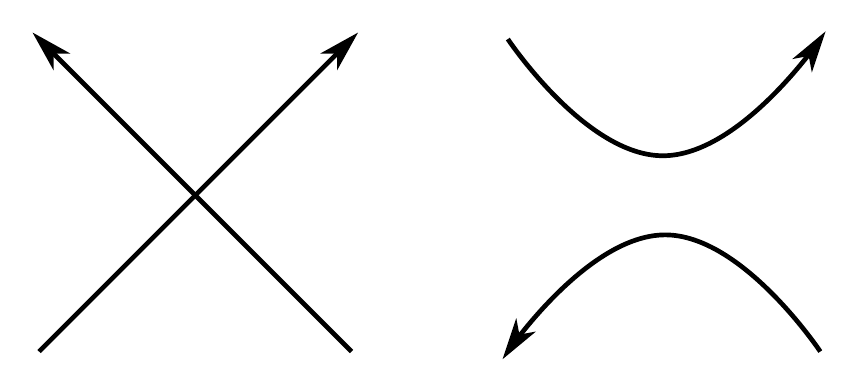
    \hfill
    \centering
    \def\svgwidth{3cm}
    \import{images/inkscape/}{AbCdEfBaGeDhIcFgHi_interlace-graph.pdf_tex}
    \hspace{0.5cm}
    \def\svgwidth{3cm}
    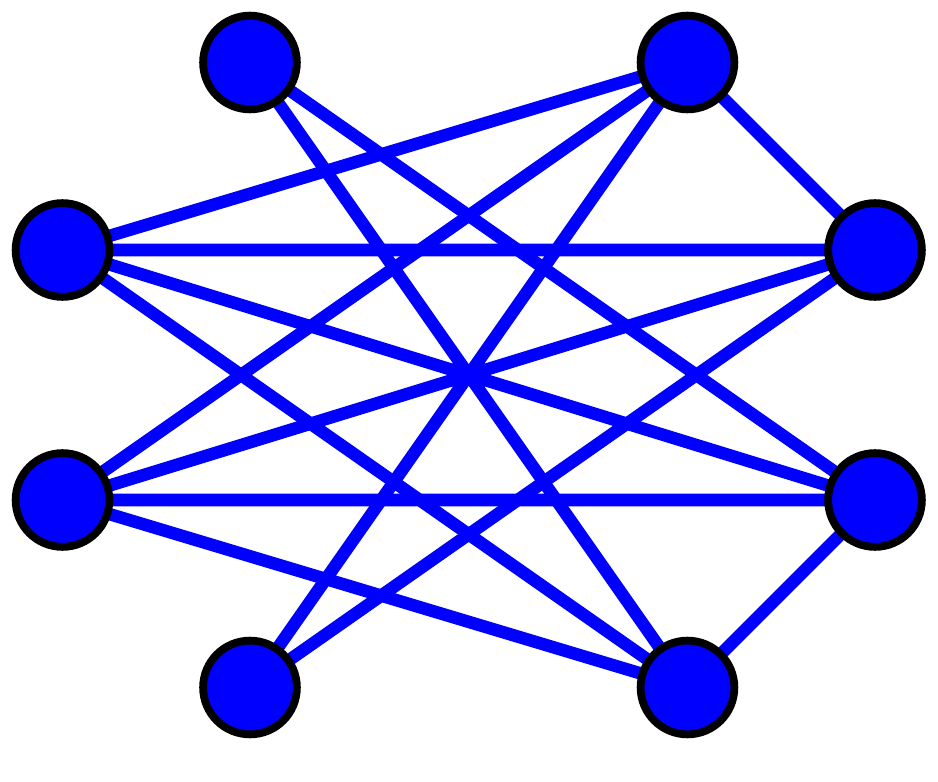
    \caption{Smoothing a loop at an intersection, smoothing a graph at a vertex}
    \label{fig:smoothing-loop-graph}
\end{figure}

\begin{Definition}[Relative evenness properties]
Consider a simple graph $G$ with adjacency matrix $\E$ and Rosenstiehl form $\CR=\E+\E^2$, so that $\forall x,y\in V_G$ we have $\CR(x,y) = \E(x,y)+\Card N(x)\cap N(y)$.
The restriction of $\CR$ to the edges $E_G$ defines the Rosenstiehl cocycle $\CR\in Z^1(G;\Z/2)$.

We say that graph $G$ satisfies the property:
\begin{enumerate}
    \item[\ref{EveN1}] at $v\in V_G$ when $\CR(v,v)\equiv 0$.
    \item[\ref{EveN1}] on $W\subset V_G$ when $\forall v\in W$ we have $\CR(v,v)\equiv 0$.
    \item[\ref{EveN2}] at $v$ when $\forall x\in V_G$, we have $\E(v,x)=0\implies \CR(v,x) \equiv 0$.
    \item[\ref{EveN2}] on $W\subset V_G^2$ when $\forall (x,y)\in W$ with $x\ne y$, we have $\E(x,y)=0\implies \CR(x,y) \equiv 0$.
    \item[\ref{RoCo}] at $v\in V_G$ when $\forall x,y\in N_G(v)$, we have $\E(x,y)\ne 0 \implies \CR(x,y)\equiv \CR(v,x)+\CR(y,v)$. 
    \item[\ref{RoCo}] relative to $W\subset V_G$ when the integral of $\CR$ vanishes along paths $(x_0,\dots,x_l)\in V_G^{l+1}$ with $\prod_1^l \E(x_{i-1},x_{i})\ne 0$ and either $x_0=x_l$ or else $x_0,x_l\in W$ and $\CR(x_0,x_l)=0$.
    
    (Note that if $U \subset W$ then \ref{RoCo} relative to $W$ implies \ref{RoCo} relative to $U$.)
\end{enumerate}
\end{Definition}

\begin{Lemma}
    \label{lem:Gaussian-smoothing}
    Consider a graph $G$ with a vertex $v\in V_G$, and let $G'=(G*v)\setminus v$ be its \emph{smoothing}.

    If $G$ is Gaussian then $G'$ is Gaussian and furthermore satisfies \ref{RoCo} relative to $N_G(v)$.
    
    If $G'$ is Gaussian and satisfies \ref{RoCo} relative to $N_G(v)$ then $G$ is Gaussian provided it satisfies \ref{EveN1} and \ref{RoCo} at $v$ as well as \ref{EveN2} on $\{v\}\cup N_G(v)$.
    
    More precisely:
    \begin{enumerate}
        \item If $G$ satisfies \ref{EveN1} at $v$, then $G$ satisfies \ref{EveN1} if and only if $G'$ satisfies \ref{EveN1}. 
        
        \item If $G$ satisfies 
        \ref{RoCo} at $v$, then $G$ satisfies \ref{EveN2} implies $G'$ satisfies \ref{EveN2}.

        If $G$ satisfies \ref{EveN2} on $\{v\}\cup N_G(v)$, then $G'$ satisfies \ref{EveN2} implies $G$ satisfies \ref{EveN2}.
        
        \item If $G$ satisfies \ref{EveN2} on $\{v\}\cup N_G(v)$, 
        then $G$ satisfies \ref{RoCo} implies $G'$ satisfies \ref{RoCo}.

        If $G$ satisfies \ref{EveN2} on $\{v\}\cup N_G(v)$ and \ref{RoCo} at $v$, then $G'$ satisfies \ref{RoCo} relative to $N_G(v)$ implies $G$ satisfies \ref{RoCo}.
    \end{enumerate}
    
\end{Lemma}

\begin{proof}
    We may assume $G$ to be connected by restricting it to the connected component of $v$.

    \emph{1.} Suppose $G$ satisfies \ref{EveN1} at $v$ and let $x\in V_{G'}$.
    
    If $x\notin N_G(v)$ then $N_{G'}(x) = N_G(x)$.
    
    If $x\in N_G(v)$ then $N_{G'}(x) = (N_G(x)\setminus\{v\}) \triangle (N_G(v)\setminus\{x\})= N_G(x)\triangle \{x,v\} \triangle N_G(v)$ so by additivity of the cardinal modulo $2$ we have $\Card N_{G'}(x) \equiv \Card N_G(x)$.

    \emph{2.} Suppose $G$ satisfies \ref{EveN2} at $v$ and \ref{RoCo} at $v$. Let $x,y\in V_{G'}$ be distinct.
    
    If $x,y\notin N_G(v)$ then $\E(x,y)=\E'(x,y)$ and $N_{G'}(x)\cap N_{G'}(y)=N_{G}(x)\cap N_{G}(y)$.
    
    If $x\in N_G(v)$ and $y\notin N_G(v)$ then on the one hand $\E'(x,y)=\E(x,y)$ and on the other hand, using the previous computations and the distributivity of $\cap$ over $\triangle$, we find that:
    \begin{align*}
        N_{G'}(x)\cap N_{G'}(y)
        &= (N_{G}(x)\triangle \{x,v\} \triangle N_G(v))
        \cap N_G(y) \\
        &= \big(N_{G}(x)\cap N_G(y)\big) \triangle \big(\{x,v\} \cap N_G(y)\big) \triangle \big(N_G(v)\cap N_G(y)\big)
    \end{align*}
    so in this case $\CR'(x,y)-\E'(x,y) = \CR(x,y)+\CR(v,y)$, and as $G$ satisfies \ref{EveN2} at $v$ we have: 
    \begin{equation*} 
        \CR'(x,y)-\E'(x,y)=\CR(x,y).
    \end{equation*}
    In particular, if $\E'(x,y)=\E(x,y)=0$ then $\CR'(x,y)=\CR(x,y)=0$, proving the second point on $(x,y)$.
    However, if $\E'(x,y)=\E(x,y)=1$ then $\CR'(x,y)=1+\CR(x,y)$ and this will serve later to prove 3.
    
    If $x,y\in N_G(v)$, then on the one hand $\E'(x,y)=1-\E(x,y)$ and on the other hand we compute:
    \begin{align*}
        N_{G'}(x)\cap N_{G'}(y) 
        & = 
        \big((N_G(x)\setminus\{v\}) \triangle (N_G(v)\setminus\{x\})\big) \cap \big((N_G(y)\setminus\{v\}) \triangle (N_G(v)\setminus\{y\})\big) 
        \\
        &= \big(N_{G}(x)\cap N_G(y) \setminus\{v\}\big)
        \triangle \big(N_G(v)\setminus\{x\} \cap N_G(v)\setminus\{y\} \big) 
        \\
        &\triangle \big(N_G(x) \cap N_G(v)\setminus\{y\}\big)
        \triangle \big(N_G(y)\cap N_G(v)\setminus\{x\}\big)
    \end{align*}
    so $\CR'(x,y)-\E'(x,y) \equiv \big(\CR(x,y)-\E(x,y)-1\big)+\big(\CR(x,v)-\E(x,v)-1\big)+\big(\CR(y,v)-\E(y,v)-1\big)$ whereby:
    \begin{equation*}
        \CR'(x,y) \equiv \CR(v,x) + \CR(x,y) + \CR(y,v).
    \end{equation*}
    Thus when $\E'(x,y)=1-\E(x,y)=0$, as $G$ satisfies \ref{RoCo} at $v$ we have $\CR(v,x)+\CR(x,y)+\CR(y,v)\equiv 0$, whence $\CR'(x,y)\equiv 0$, and this concludes the proof of 2.
    
    However when $\E'(x,y)=1-\E(x,y)=1$, if $G$ satisfies \ref{EveN2} on $\{v\}\cup N_G(v)$ then $\CR(x,y)\equiv0$, whereby $\CR'(x,y) \equiv \CR(x,v) + \CR(v,y)$, and this will serve to prove 3. 

    \emph{3.} Suppose $G$ satisfies \ref{EveN2} on $\{v\}\cup N_G(v)$ and \ref{RoCo} at $v$ (so all previous calculations hold).
    
    Consider a relative cycle $\gamma'\in Z_1(G', N_G(v))$, represented as a sequence $x_0,\dots, x_l \in V_{G'}^{l+1}$ with $\prod_1^l \E'(x_{i-1},x_{i})=1$, and either $x_0=x_l$ or else $x_0,x_l\in N_G(v)$ and $\CR(x_0,x_l)=0$.
    %
    %
    We construct a cycle $\gamma\in Z_1(G)$ obtained from $\gamma'$ 
    by replacing every segment $(x_{i-1},x_{i}) \in N_G(v)^2$ with the segment $(x_{i-1},v,x_{i})$, and appending $v$ if $x_0 \ne x_l$.

    Let us show $\CR'(\gamma')=\CR(\gamma)$.
    If $x,y\in \gamma'$ avoid $N_G(v)$ then $\CR(x,y)=\CR'(x,y)$.
    The edges $(x,y)$ of $\gamma'$ with exactly one endpoint in $N_G(v)$ are in even number, and they all satisfy $\CR'(x,y)=1+\CR(x,y)$.
    If $x,y\in \gamma'$ are both in $N_G(v)$ then $\CR'(x,y)=\CR(x,v)+\CR(v,y)$.
    The first statement in 3 is proven.

    Note that the inverse map $\gamma \mapsto \gamma'$, obtained by deleting all appearances of $v$, is well-defined (only) for cycles in which every occurrence of $v$ appears between $x,y$ with $\E(x,y)=0$.
    Hence every cycle in $Z_1(G;\Z/2)$ is homologous to a sum of cycles in $Z_1(\{v\}\cup N_{G}(v);\Z/2)$ and relative cycles in $Z_1(G',N_G(V))$ whose endpoints are not connected in $G'$.
    The converse statement in 3 follows. 
    %
    %
\end{proof}

\begin{Remark}[Komposition]
    For framed rooted chord diagrams $\alpha = a_\infty U a_0 V$ and $\beta = b_\infty X b_0 Y$, we define their \emph{komposition} as the framed rooted chord diagram $\alpha \triangle  \beta = c_\infty UX c_0 VY$.
    The interlace graph of a komposition is the komposition of interlace graphs (defined in \ref{Def:Composition-graphs}).
    
    The topological construction for the corresponding based filoops is depicted in Figure~\ref{fig:skew-plumbing-loops}.
    The genus is additive under komposition, so Gaussian chordiagraphs are stable under komposition.

    This operation was defined in \cite{Ghys_promenade_2017} to construct an operad on spheriloops.
    
\begin{figure}[h]
    \centering
    \def\svgwidth{12cm}
    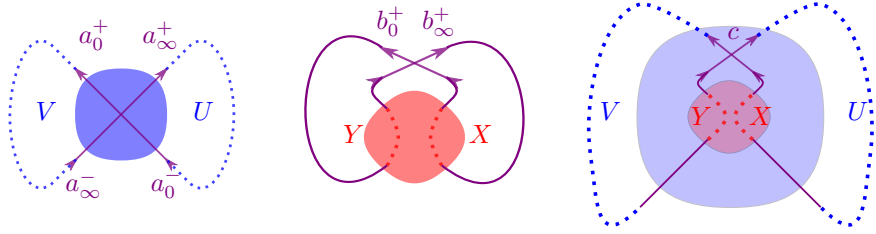
    \caption{Based filoops $\alpha$ and $\beta$, and their komposition $\alpha \triangle \beta$.}
    \label{fig:skew-plumbing-loops}
\end{figure}
\end{Remark}

\begin{Corollary}
    If $(G_A,a)$ and $(G_B,b)$ are Gaussian then so is $(G_A,a) \triangle (G_B,b)$.
    
    If $(G_A,a) \triangle (G_B,b)$ is Gaussian then $\Card N_A(a)$ and $\Card N_B(b)$ have the same parity, and if it is even then $(G_A,a)$ and $(G_B,b)$ are Gaussian.
\end{Corollary}

\begin{proof}
    Assume $(G_A,a)$ and $(G_B,b)$ are Gaussian and let $(G_C,c)=(G_A,a)\triangle(G_B,b)$. By Lemma \ref{lem:Gaussian-smoothing}, the smoothings $(G_A*a)\setminus a$ and $(G_B*b)\setminus b$ are Gaussian and satisfy \ref{RoCo} relative to $N_A(a)$ and $N_B(b)$. Hence the smoothing $G_C'= (G_C*c)\setminus c = ((G_A*a)\setminus a)\sqcup((G_B*b)\setminus b)$ is Gaussian and satisfies \ref{RoCo} relative to $N_C(c)=N_A(a)\sqcup N_B(b)$. Moreover $\Card N_C(c)=\Card N_A(a)+\Card N_B(b)$ is even (by \ref{EveN1} for $G_A$ and $G_B$), so $G_C$ satisfies \ref{EveN1} at $c$.
    We show that $G_C$ satisfies \ref{EveN2} on ${c}\cup N_C(c)={c}\sqcup N_A(a)\sqcup N_B(b)$: for $x,y\in N_A(a)$ we have $\CR_C(x,y)\equiv \CR_A(x,y)$ (because $N_B(b)$ is even) and for $x,y\in N_B(b)$ we have $\CR_C(x,y)\equiv\CR_B(x,y)$ (because $N_A(a)$ is even), and that is all we must check (as all other pairs of distinct vertices in ${c}\cup N_C(c)$ are connected).
    We finally show that $G_C$ satisfies \ref{RoCo} at $c$: this follows from the corresponding property for $G_A$ at $a$ and for $G_B$ at $b$, while for $x\in N_A(a)$ and $y\in N_B(b)$ we have by construction $\CR_C(x,y)-1\equiv (\CR_A(x,a)-1)+(\CR_B(b,y)-1)+1$ but also $\CR_C(x,c)\equiv \CR_A(x,a)$ and $\CR_B(b,y)\equiv \CR_C(c,y)$, hence $\CR_C(x,y) \equiv \CR_C(x,c)+\CR_C(c,y)$ as desired.
    Thus Lemma \ref{lem:Gaussian-smoothing} implies that $G_C$ is Gaussian.

    If $(G_C,c)=(G_A,a) \triangle (G_B,b)$ is Gaussian then its GLT factorisations satisfy the weighted Gaussian conditions in Proposition \ref{Prop:GLT-Gauss}.
    Consider a (possibly non-reduced) GLT factorisation of $(G_C,c)$ in which we may identify GLT factorisations of $(G_A,a)$ and $(G_B,b)$ composed into a rooted triangle.
    
    The vertices of the triangle have weights $1$, $\Card N_A(a)$, $\Card N_B(b)$, so the Gaussian hypothesis implies that $\Card N_A(a)\equiv \Card N_B(b)\bmod{2}$. 
    The weights of the vertices in the factors of $G_C$ are obtained from those of $G_A$ and $G_B$ by adding $\Card N_B(b)$ to the vertices directly accessible from $a$ in the GLT of $G_A$ and $\Card N_A(a)$ to the vertices directly accessible from $b$ in the GLT of $G_B$. 

    Consequently if these cardinals are even, the GLT factorisations of $(G_A,a)$ and $(G_B,b)$ satisfy the weighted Gaussian conditions in Proposition \ref{Prop:GLT-Gauss}.
\end{proof}

\begin{Question}[Converse]
    For the converse, the assumption that $\Card N_A(a)\equiv \Card N_B(b)$ is even cannot, in general, be lifted. Indeed $K_m\triangle K_n=K_{1+m+n}$, but a clique is Gaussian if and only if it has an odd number of vertices.
    Can one find a reduced counterexample?
    
    Can one find graphs $(G_A,a)$ and $(G_B,b)$ such that $(G_A,a)\triangle (G_B,b)$ is Gaussian but such that $\Card N_A(a)\equiv \Card N_B(b)$ is odd?
    Can one find a reduced counterexample?
\end{Question}

\begin{Question}[Plumbing Gaussian graphs]
We showed that Gaussian graphs are stable under smoothing and komposition, both of which can be constructed from plumbings.

We know that Gaussian chordiagraphs are stable under plumbing, and by Theorem \ref{Thm:unifacto-filoops} that they admit unique reduced factorisations into plumbings of CL2-chordiagraphs.
However it is not clear how to use Proposition \ref{Prop:GLT-Gauss} to show that all Gaussian graphs stable under plumbing: are they?
\end{Question}

\bibliographystyle{alpha}
\bibliography{biblio_CDGP.bib}

@article{FPST_morsifications-mutations_2022,
    AUTHOR = {Fomin, Sergey and Pylyavskyy, Pavlo and Shustin, Eugenii and
              Thurston, Dylan},
     TITLE = {Morsifications and mutations},
   JOURNAL = {J. Lond. Math. Soc. (2)},
  FJOURNAL = {Journal of the London Mathematical Society. Second Series},
    VOLUME = {105},
      YEAR = {2022},
    NUMBER = {4},
     PAGES = {2478--2554},
      ISSN = {0024-6107,1469-7750},
   MRCLASS = {13F60 (16G20 20F36 57K10 58K65)},
  MRNUMBER = {4440540},
       DOI = {10.1112/jlms.12566},
       URL = {https://doi.org/10.1112/jlms.12566},
}

@article{Arratia-Bollobas-Sorkin_interlace-polynomial_2004,
    AUTHOR = {Arratia, Richard and Bollob\'{a}s, B\'{e}la and Sorkin,
              Gregory B.},
     TITLE = {The interlace polynomial of a graph},
   JOURNAL = {J. Combin. Theory Ser. B},
  FJOURNAL = {Journal of Combinatorial Theory. Series B},
    VOLUME = {92},
      YEAR = {2004},
    NUMBER = {2},
     PAGES = {199--233},
      ISSN = {0095-8956,1096-0902},
   MRCLASS = {05C45},
  MRNUMBER = {2099142},
MRREVIEWER = {Gary\ Gordon},
       DOI = {10.1016/j.jctb.2004.03.003},
       URL = {https://doi.org/10.1016/j.jctb.2004.03.003},
}

@book{Lando-Zvonkin_graphs-on-surfaces_2004,
    AUTHOR = {Lando, Sergei K. and Zvonkin, Alexander K.},
     TITLE = {Graphs on surfaces and their applications},
    SERIES = {Encyclopaedia of Mathematical Sciences},
    VOLUME = {141},
      NOTE = {With an appendix by Don B. Zagier,
              Low-Dimensional Topology, II},
 PUBLISHER = {Springer-Verlag, Berlin},
      YEAR = {2004},
     PAGES = {xvi+455},
      ISBN = {3-540-00203-0},
   MRCLASS = {14H55 (05-02 05C10 05C30 05C50 14H10 14H30 32G15)},
  MRNUMBER = {2036721},
MRREVIEWER = {Athanase\ Papadopoulos},
       DOI = {10.1007/978-3-540-38361-1},
       URL = {https://doi.org/10.1007/978-3-540-38361-1},
}

@book{Arnold_Topo-plane-curves-caustics_2000,
    AUTHOR = {Arnol'd, V. I.},
     TITLE = {Topological invariants of plane curves and caustics},
    SERIES = {University Lecture Series},
    VOLUME = {5},
      NOTE = {Dean Jacqueline B. Lewis Memorial Lectures presented at
              Rutgers University, New Brunswick, New Jersey},
 PUBLISHER = {American Mathematical Society, Providence, RI},
      YEAR = {1994},
     PAGES = {viii+60},
      ISBN = {0-8218-0308-5},
   MRCLASS = {57M25 (53A04 57-02 57R42 58C28)},
  MRNUMBER = {1286249},
MRREVIEWER = {Louis\ Funar},
       DOI = {10.1090/ulect/005},
       URL = {https://doi.org/10.1090/ulect/005},
}

@book{ChDuMo_Vassiliev_2012,
    AUTHOR = {Chmutov, S. and Duzhin, S. and Mostovoy, J.},
     TITLE = {Introduction to {V}assiliev knot invariants},
 PUBLISHER = {Cambridge University Press, Cambridge},
      YEAR = {2012},
     PAGES = {xvi+504},
      ISBN = {978-1-107-02083-2},
   MRCLASS = {57M27},
  MRNUMBER = {2962302},
MRREVIEWER = {Hitoshi\ Murakami},
       DOI = {10.1017/CBO9781139107846},
       URL = {https://doi.org/10.1017/CBO9781139107846},
}

@book{Ghys_promenade_2017,
    AUTHOR = {Ghys, \'{E}tienne},
     TITLE = {A singular mathematical promenade},
 PUBLISHER = {ENS \'{E}ditions, Lyon},
      YEAR = {2017},
     PAGES = {viii+302},
      ISBN = {978-2-84788-939-0},
   MRCLASS = {14H20 (00A09 01A05 14P25 18D50 55P48 57M25 58K05)},
  MRNUMBER = {3702027},
MRREVIEWER = {Evelyn James Lamb},
}

@article{Tao_polynomial-method_2014,
    AUTHOR = {Tao, Terence},
     TITLE = {Algebraic combinatorial geometry: the polynomial method in
              arithmetic combinatorics, incidence combinatorics, and number
              theory},
   JOURNAL = {EMS Surv. Math. Sci.},
  FJOURNAL = {EMS Surveys in Mathematical Sciences},
    VOLUME = {1},
      YEAR = {2014},
    NUMBER = {1},
     PAGES = {1--46},
      ISSN = {2308-2151,2308-216X},
   MRCLASS = {14-02 (05-02 11B30 14G15 14Q99 51E30)},
  MRNUMBER = {3200226},
MRREVIEWER = {Victoria\ Powers},
       DOI = {10.4171/EMSS/1},
       URL = {https://doi.org/10.4171/EMSS/1},
}

@article{Chmutov-Pak_Kauffman-virtualink-Tutte_2007,
    AUTHOR = {Chmutov, Sergei and Pak, Igor},
     TITLE = {The {K}auffman bracket of virtual links and the
              {B}ollob\'{a}s-{R}iordan polynomial},
   JOURNAL = {Mosc. Math. J.},
  FJOURNAL = {Moscow Mathematical Journal},
    VOLUME = {7},
      YEAR = {2007},
    NUMBER = {3},
     PAGES = {409--418, 573},
      ISSN = {1609-3321,1609-4514},
   MRCLASS = {57M15 (05C10 05C22 57M25 57M27)},
  MRNUMBER = {2343139},
MRREVIEWER = {Gw\'{e}na\"{e}l\ Massuyeau},
       DOI = {10.17323/1609-4514-2007-7-3-409-418},
       URL = {https://doi.org/10.17323/1609-4514-2007-7-3-409-418},
}

@incollection{Conway_enumeration-knots-links-tangles_1970,
    AUTHOR = {Conway, J. H.},
     TITLE = {An enumeration of knots and links, and some of their algebraic
              properties},
 BOOKTITLE = {Computational {P}roblems in {A}bstract {A}lgebra ({P}roc.
              {C}onf., {O}xford, 1967)},
     PAGES = {329--358},
 PUBLISHER = {Pergamon, Oxford-New York-Toronto, Ont.},
      YEAR = {1970},
   MRCLASS = {55.20},
  MRNUMBER = {258014},
MRREVIEWER = {H.\ E.\ Debrunner},
}

@article{Scott_geometry-3-manifolds_1983,
    AUTHOR = {Scott, Peter},
     TITLE = {The geometries of {$3$}-manifolds},
   JOURNAL = {Bull. London Math. Soc.},
  FJOURNAL = {The Bulletin of the London Mathematical Society},
    VOLUME = {15},
      YEAR = {1983},
    NUMBER = {5},
     PAGES = {401--487},
      ISSN = {0024-6093,1469-2120},
   MRCLASS = {57N10 (22E10 53C20)},
  MRNUMBER = {705527},
MRREVIEWER = {John\ Hempel},
       DOI = {10.1112/blms/15.5.401},
       URL = {https://doi.org/10.1112/blms/15.5.401},
}

@article{Thurston_geometrisation_1982,
    AUTHOR = {Thurston, William P.},
     TITLE = {Three-dimensional manifolds, {K}leinian groups and hyperbolic
              geometry},
   JOURNAL = {Bull. Amer. Math. Soc. (N.S.)},
  FJOURNAL = {American Mathematical Society. Bulletin. New Series},
    VOLUME = {6},
      YEAR = {1982},
    NUMBER = {3},
     PAGES = {357--381},
      ISSN = {0273-0979,1088-9485},
   MRCLASS = {57N10 (20H15 30F40 57M35 57M40 57S17)},
  MRNUMBER = {648524},
MRREVIEWER = {Klaus\ Johannson},
       DOI = {10.1090/S0273-0979-1982-15003-0},
       URL = {https://doi.org/10.1090/S0273-0979-1982-15003-0},
}

@inproceedings{Gauss_words-plane-loops_1807,
    author = {Carl Frederich Gauss},
    title = {Zur Geometria Situs},
    booktitle = {Werke (Nachlass)},
    volume = {VIII},
    year = {1807},
    pages = {271--286},
    editor = {Schilling, C.},
    publisher={Georg Olms Verlag, Hildesheim-New York}
}

@article{Dowker-Thistlethwaite_classification-knot-projections_1983,
    AUTHOR = {Dowker, C. H. and Thistlethwaite, Morwen B.},
     TITLE = {Classification of knot projections},
   JOURNAL = {Topology Appl.},
  FJOURNAL = {Topology and its Applications},
    VOLUME = {16},
      YEAR = {1983},
    NUMBER = {1},
     PAGES = {19--31},
      ISSN = {0166-8641,1879-3207},
   MRCLASS = {57M25},
  MRNUMBER = {702617},
       DOI = {10.1016/0166-8641(83)90004-4},
       URL = {https://doi.org/10.1016/0166-8641(83)90004-4},
}

@article{Fraysseix-Mendez_Gauss-codes_1999,
    AUTHOR = {de Fraysseix, H. and Ossona de Mendez, P.},
     TITLE = {On a characterization of {G}auss codes},
   JOURNAL = {Discrete Comput. Geom.},
  FJOURNAL = {Discrete \& Computational Geometry. An International Journal
              of Mathematics and Computer Science},
    VOLUME = {22},
      YEAR = {1999},
    NUMBER = {2},
     PAGES = {287--295},
      ISSN = {0179-5376},
   MRCLASS = {05C10 (94B60)},
  MRNUMBER = {1698548},
MRREVIEWER = {M. L. Marx},
       DOI = {10.1007/PL00009461},
       URL = {https://doi.org/10.1007/PL00009461},
}

@article{Rosenstiehl_Gauss-interlace_1999,
    AUTHOR = {Rosenstiehl, Pierre},
     TITLE = {A new proof of the {G}auss interlace conjecture},
   JOURNAL = {Adv. in Appl. Math.},
  FJOURNAL = {Advances in Applied Mathematics},
    VOLUME = {23},
      YEAR = {1999},
    NUMBER = {1},
     PAGES = {3--13},
      ISSN = {0196-8858},
   MRCLASS = {05C10},
  MRNUMBER = {1692992},
MRREVIEWER = {M. L. Marx},
       DOI = {10.1006/aama.1999.0643},
       URL = {https://doi.org/10.1006/aama.1999.0643},
}

@article{Bouchet_Circle-graph-obstructions_1994,
    AUTHOR = {Bouchet, Andr\'{e}},
     TITLE = {Circle graph obstructions},
   JOURNAL = {J. Combin. Theory Ser. B},
  FJOURNAL = {Journal of Combinatorial Theory. Series B},
    VOLUME = {60},
      YEAR = {1994},
    NUMBER = {1},
     PAGES = {107--144},
      ISSN = {0095-8956,1096-0902},
   MRCLASS = {05C05},
  MRNUMBER = {1256586},
MRREVIEWER = {Hubert\ de Fraysseix},
       DOI = {10.1006/jctb.1994.1008},
       URL = {https://doi.org/10.1006/jctb.1994.1008},
}

@article{Bouchet_prime-graphs-circle-graphs_1987,
    AUTHOR = {Bouchet, Andr\'{e}},
     TITLE = {Reducing prime graphs and recognizing circle graphs},
   JOURNAL = {Combinatorica},
  FJOURNAL = {Combinatorica. An International Journal of the J\'{a}nos
              Bolyai Mathematical Society},
    VOLUME = {7},
      YEAR = {1987},
    NUMBER = {3},
     PAGES = {243--254},
      ISSN = {0209-9683},
   MRCLASS = {05C75},
  MRNUMBER = {918395},
MRREVIEWER = {R.\ C.\ Read},
       DOI = {10.1007/BF02579301},
       URL = {https://doi.org/10.1007/BF02579301},
}

@article{Chaves-Weber_plombages-rubans-mots-Gauss_1994,
    AUTHOR = {Chaves, Nathalie and Weber, Claude},
     TITLE = {Plombages de rubans et probl\`eme des mots de {G}auss},
   JOURNAL = {Exposition. Math.},
  FJOURNAL = {Expositiones Mathematicae. International Journal for Pure and
              Applied Mathematics},
    VOLUME = {12},
      YEAR = {1994},
    NUMBER = {1},
     PAGES = {53--77},
      ISSN = {0723-0869},
   MRCLASS = {57M25},
  MRNUMBER = {1267628},
MRREVIEWER = {Fran\c{c}ois\ Jaeger},
}

@incollection{Weber_classical-knot-theory_2001,
    AUTHOR = {Weber, Claude},
     TITLE = {Elements of classical knot theory},
 BOOKTITLE = {An introduction to the geometry and topology of fluid flows
              ({C}ambridge, 2000)},
    SERIES = {NATO Sci. Ser. II Math. Phys. Chem.},
    VOLUME = {47},
     PAGES = {57--75},
 PUBLISHER = {Kluwer Acad. Publ., Dordrecht},
      YEAR = {2001},
      ISBN = {1-4020-0206-8},
   MRCLASS = {57M25},
  MRNUMBER = {1999081},
       DOI = {10.1007/978-94-010-0446-6\_4},
       URL = {https://doi.org/10.1007/978-94-010-0446-6_4},
}

@article{Valette_classification-spheriloops-6-gauss-diagram_2016,
    AUTHOR = {Valette, Guy},
     TITLE = {A classification of spherical curves based on {G}auss
              diagrams},
   JOURNAL = {Arnold Math. J.},
  FJOURNAL = {Arnold Mathematical Journal},
    VOLUME = {2},
      YEAR = {2016},
    NUMBER = {3},
     PAGES = {383--405},
      ISSN = {2199-6792,2199-6806},
   MRCLASS = {57M50 (53A04)},
  MRNUMBER = {3550432},
MRREVIEWER = {Shy\={u}ichi\ Izumiya},
       DOI = {10.1007/s40598-016-0049-3},
       URL = {https://doi.org/10.1007/s40598-016-0049-3},
}

@article{Gio-Paul_split-decomp_2012,
    AUTHOR = {Gioan, Emeric and Paul, Christophe},
     TITLE = {Split decomposition and graph-labelled trees:
              characterizations and fully dynamic algorithms for totally
              decomposable graphs},
   JOURNAL = {Discrete Appl. Math.},
  FJOURNAL = {Discrete Applied Mathematics. The Journal of Combinatorial
              Algorithms, Informatics and Computational Sciences},
    VOLUME = {160},
      YEAR = {2012},
    NUMBER = {6},
     PAGES = {708--733},
      ISSN = {0166-218X,1872-6771},
   MRCLASS = {05C85 (05C05 05C12)},
  MRNUMBER = {2901084},
       DOI = {10.1016/j.dam.2011.05.007},
       URL = {https://doi.org/10.1016/j.dam.2011.05.007},
}

@article{Cunningham_decomp-di-graphs_1982,
    AUTHOR = {Cunningham, William H.},
     TITLE = {Decomposition of directed graphs},
   JOURNAL = {SIAM J. Algebraic Discrete Methods},
  FJOURNAL = {Society for Industrial and Applied Mathematics. Journal on
              Algebraic and Discrete Methods},
    VOLUME = {3},
      YEAR = {1982},
    NUMBER = {2},
     PAGES = {214--228},
      ISSN = {0196-5212},
   MRCLASS = {05C20 (05C70 68E10)},
  MRNUMBER = {655562},
       DOI = {10.1137/0603021},
       URL = {https://doi.org/10.1137/0603021},
}

@article{Lovasz-Marx_forbidden-substructure-Gauss-codes_1976,
    AUTHOR = {Lov\'{a}sz, L\'{a}szl\'{o} and Marx, Morris L.},
     TITLE = {A forbidden substructure characterization of {G}auss codes},
   JOURNAL = {Acta Sci. Math. (Szeged)},
  FJOURNAL = {Acta Universitatis Szegediensis. Acta Scientiarum
              Mathematicarum},
    VOLUME = {38},
      YEAR = {1976},
    NUMBER = {1-2},
     PAGES = {115--119},
      ISSN = {0001-6969},
   MRCLASS = {94A10 (05C10)},
  MRNUMBER = {441540},
MRREVIEWER = {Andr\'{a}s\ Recski},
}

@article{Charbit-Montgolfier-Raffinot_Split-decomp_2012,
    AUTHOR = {Charbit, Pierre and de Montgolfier, Fabien and Raffinot,
              Mathieu},
     TITLE = {Linear time split decomposition revisited},
   JOURNAL = {SIAM J. Discrete Math.},
  FJOURNAL = {SIAM Journal on Discrete Mathematics},
    VOLUME = {26},
      YEAR = {2012},
    NUMBER = {2},
     PAGES = {499--514},
      ISSN = {0895-4801,1095-7146},
   MRCLASS = {68R10 (05C51 05C85)},
  MRNUMBER = {2967479},
MRREVIEWER = {Matthew\ Yancey},
       DOI = {10.1137/10080052X},
       URL = {https://doi.org/10.1137/10080052X},
}

@incollection{Crapo-Rosenstiehl_bicoloops_2001,
    AUTHOR = {Crapo, Henry and Rosenstiehl, Pierre},
     TITLE = {On lacets and their manifolds},
      NOTE = {Graph theory (Prague, 1998)},
   JOURNAL = {Discrete Math.},
  FJOURNAL = {Discrete Mathematics},
    VOLUME = {233},
      YEAR = {2001},
    NUMBER = {1-3},
     PAGES = {299--320},
      ISSN = {0012-365X,1872-681X},
   MRCLASS = {05C10 (30G12 57M15)},
  MRNUMBER = {1825623},
MRREVIEWER = {M.\ L.\ Marx},
       DOI = {10.1016/S0012-365X(00)00248-X},
       URL = {https://doi.org/10.1016/S0012-365X(00)00248-X},
}

@article{Turaev_CurveSurfacesChartsWords_2005,
    AUTHOR = {Turaev, Vladimir},
     TITLE = {Curves on surfaces, charts, and words},
   JOURNAL = {Geom. Dedicata},
  FJOURNAL = {Geometriae Dedicata},
    VOLUME = {116},
      YEAR = {2005},
     PAGES = {203--236},
      ISSN = {0046-5755,1572-9168},
   MRCLASS = {57M99 (68R15)},
  MRNUMBER = {2195447},
MRREVIEWER = {Igor\ Rivin},
       DOI = {10.1007/s10711-005-9013-4},
       URL = {https://doi.org/10.1007/s10711-005-9013-4},
}

@article{Turaev_Virtual-strings_2015,
    AUTHOR = {Turaev, Vladimir},
     TITLE = {Virtual strings},
   JOURNAL = {Ann. Inst. Fourier (Grenoble)},
  FJOURNAL = {Universit\'{e} de Grenoble. Annales de l'Institut Fourier},
    VOLUME = {54},
      YEAR = {2004},
    NUMBER = {7},
     PAGES = {2455--2525},
      ISSN = {0373-0956,1777-5310},
   MRCLASS = {57M27 (16W30 57M25 57N05)},
  MRNUMBER = {2139700},
MRREVIEWER = {Vladimir\ V.\ Tchernov},
       URL = {http://aif.cedram.org/item?id=AIF_2004__54_7_2455_0},
}

@article{Gibson_Tabulating-Virtual-strings_2008,
    AUTHOR = {Gibson, Andrew},
     TITLE = {On tabulating virtual strings},
   JOURNAL = {Acta Math. Vietnam.},
  FJOURNAL = {Acta Mathematica Vietnamica},
    VOLUME = {33},
      YEAR = {2008},
    NUMBER = {3},
     PAGES = {493--518},
      ISSN = {0251-4184},
   MRCLASS = {57M27 (57M25)},
  MRNUMBER = {2501855},
MRREVIEWER = {Heather\ A.\ Dye},
}

@article{Ito_finite-type-inv-curveSurfaces_2009,
    AUTHOR = {Ito, Noboru},
     TITLE = {Finite-type invariants for curves on surfaces},
   JOURNAL = {Proc. Japan Acad. Ser. A Math. Sci.},
  FJOURNAL = {Japan Academy. Proceedings. Series A. Mathematical Sciences},
    VOLUME = {85},
      YEAR = {2009},
    NUMBER = {9},
     PAGES = {129--134},
      ISSN = {0386-2194},
   MRCLASS = {57M27},
  MRNUMBER = {2573961},
MRREVIEWER = {Teruhisa\ Kadokami},
       URL = {http://projecteuclid.org/euclid.pja/1257430680},
}

@article{Cahn-Levi_VassilInv-VirtuaLegendKnots_2015,
    AUTHOR = {Cahn, Patricia and Levi, Asa},
     TITLE = {Vassiliev invariants of virtual {L}egendrian knots},
   JOURNAL = {Pacific J. Math.},
  FJOURNAL = {Pacific Journal of Mathematics},
    VOLUME = {273},
      YEAR = {2015},
    NUMBER = {1},
     PAGES = {21--46},
      ISSN = {0030-8730,1945-5844},
   MRCLASS = {57M27},
  MRNUMBER = {3290443},
MRREVIEWER = {Micah\ Whitney\ Chrisman},
       DOI = {10.2140/pjm.2015.273.21},
       URL = {https://doi.org/10.2140/pjm.2015.273.21},
}

@article{Schaeffer-Justin_asymptotic-plane-curves_2004,
    AUTHOR = {Schaeffer, Gilles and Zinn-Justin, Paul},
     TITLE = {On the asymptotic number of plane curves and alternating
              knots},
   JOURNAL = {Experiment. Math.},
  FJOURNAL = {Experimental Mathematics},
    VOLUME = {13},
      YEAR = {2004},
    NUMBER = {4},
     PAGES = {483--493},
      ISSN = {1058-6458,1944-950X},
   MRCLASS = {05A16 (05C30 57M25 83C80)},
  MRNUMBER = {2118273},
MRREVIEWER = {Valery\ A.\ Liskovets},
       URL = {http://projecteuclid.org/euclid.em/1109106440},
}

@Article{CLS_TopoDenoCourbAlgReel_2018,
     author = {Christopher-Lloyd Simon},
     title = {Topologie et d\'enombrement des courbes alg\'ebriques r\'eelles},
     journal = {Annales de la Facult\'e des sciences de Toulouse : Math\'ematiques},
     pages = {383--422},
     publisher = {Universit\'e Paul Sabatier, Toulouse},
     volume = {6e s{\'e}rie, 31},
     number = {2},
     year = {2022},
     doi = {10.5802/afst.1698},
     zbl = {07549944},
     language = {fr},
     url = {https://afst.centre-mersenne.org/articles/10.5802/afst.1698/}
}

\end{document}